\documentclass[sn-mathphys,Numbered]{sn-jnl}% Math and Physical Sciences Reference Style
\usepackage{graphicx}%
\usepackage{multirow}%
\usepackage{amsmath,amssymb,amsfonts}%
\usepackage{amsthm}%
\usepackage{mathrsfs}%
\usepackage[title]{appendix}%
\usepackage{xcolor}%
\usepackage{textcomp}%
\usepackage{manyfoot}%
\usepackage{booktabs}%
\usepackage{algorithm}%
\usepackage{algorithmicx}%
\usepackage{algpseudocode}%
\usepackage{listings}%
\usepackage{lmodern}
\usepackage{mathtools}
\usepackage{amssymb}

\usepackage[shortlabels]{enumitem}
\usepackage{dsfont}
\usepackage{verbatim}
\usepackage[autostyle]{csquotes}
\usepackage{hyperref}
\usepackage{cleveref}
\usepackage{xcolor}
\usepackage{array}
\usepackage{nicefrac}

\newcommand*{\Q}{\mathbb{Q}}
\newcommand*{\R}{\mathbb{R}}

\newcommand*{\N}{\mathbb{N}}
\newcommand*{\eps}{\varepsilon}

\newcommand*{\spdx}{(\cdot,\cdot)}
\newcommand*{\normx}{\|\cdot\|}

\newcommand*{\acore}{\mathcal{D}}
\newcommand*{\sccs}{(T_t)_{t\geq 0}}

\newcommand*{\sccsz}{(T_t^0)_{t\geq 0}}

\newcommand*{\dm}{\,\mathrm{d}\mu}
\providecommand{\norm}[1]{\lVert#1\rVert} %Norm

\providecommand{\abs}[1]{\lvert#1\rvert}

\DeclareMathOperator{\tr}{tr}

\DeclareMathOperator{\diag}{diag}

\newcommand*{\defeq}{\mathrel{\vcenter{\baselineskip0.5ex \lineskiplimit0pt
			\hbox{\scriptsize.}\hbox{\scriptsize.}}}%
	=}
\newcommand*{\eqdef}{=\mathrel{\vcenter{\baselineskip0.5ex \lineskiplimit0pt
			\hbox{\scriptsize.}\hbox{\scriptsize.}}}}

\newcommand*{\bs}[1]{\mathcal{B}(#1)}
\newcommand*{\lop}[1]{\mathcal{L}(#1)}
\newcommand*{\lopp}[1]{\mathcal{L}^+(#1)}
\newcommand*{\loppt}[1]{\mathcal{L}_1^+(#1)}
\newcommand*{\lophs}[1]{\mathcal{L}_2(#1)}
\newcommand*{\s}[1]{(#1_n)_{n\in\N}}
\newcommand*{\sk}[1]{(#1_k)_{k\in\N}}

\newcommand*{\lin}[1]{\operatorname{span}\{ #1\}}
\newcommand*{\fbs}[1]{\mathcal{F}C_b^\infty(#1)}

\theoremstyle{thmstyleone}%
\newtheorem{thm}{Theorem}[section]%  meant for continuous numbers
\newtheorem{prop}[thm]{Proposition}% 
\newtheorem{corollary}[thm]{Corollary}
\newtheorem{lem}[thm]{Lemma}

\newtheorem{summary}[thm]{Summary}
\newtheorem*{ass}{Assumption}

\newtheorem{remark}[thm]{Remark}%

\newtheorem{defn}[thm]{Definition}%

\newenvironment{cond}[1]{\begin{ass}[#1] \def\@currentlabelname{(\upshape\boldmath\textbf{#1}}}{\end{ass}}
\begin{document}

\title[Article Title]{Hypocoercivity for infinite-dimensional non-linear degenerate stochastic differential equations with multiplicative noise}


\author[1]{\fnm{Alexander} \sur{Bertram}}\email{alexander.bertram@tutanota.de}

\author*[1]{\fnm{Benedikt} \sur{Eisenhuth}}\email{eisenhuth@mathematik.uni-kl.de}

\author[1]{\fnm{Martin} \sur{Grothaus}}\email{grothaus@mathematik.uni-kl.de}

\affil*[1]{\orgdiv{Department of Mathematics}, \orgname{University of Kaiserslautern-Landau}, \orgaddress{\street{Erwin-Schr{\"o}dinger-Stra{\ss}e 48}, \city{Kaiserslautern}, \postcode{67663}, \country{Germany}}}

%\affil[2]{\orgdiv{Department}, \orgname{Organization}, \orgaddress{\street{Street}, \city{City}, \postcode{10587}, \state{State}, \country{Country}}}
%
%\affil[3]{\orgdiv{Department}, \orgname{Organization}, \orgaddress{\street{Street}, \city{City}, \postcode{610101}, \state{State}, \country{Country}}}


\abstract{We analyze infinite-dimensional non-linear degenerate stochastic differential equations with multiplicative noise. First, essential m-dissipativity of their associated Kolmogorov backward generators on $L^2(\mu^{\Phi})$ defined on smooth finitely based functions is established. Here $\Phi$ is an appropriate potential and $\mu^{\Phi}$ is the invariant measure with density $e^{-\Phi}$ w.r.t.~an infinite-dimensional non-degenerate Gaussian measure. Second, we use resolvent methods to construct corresponding right processes with infinite lifetime, solving the martingale problem for the Kolmogorov backward generators. They provide weak solutions, with weakly continuous paths, to the non-linear degenerate stochastic partial differential equations. Moreover, we identify the transition semigroup of such a process with the strongly continuous contraction semigroup $(T_t)_{t\geq 0}$ generated by the corresponding Kolmogorov backwards generator. 
Afterwards, we apply a refinement of the abstract Hilbert space hypocoercivity method, developed by Dolbeaut, Mouhot and Schmeiser, to derive hypocoercivity of $(T_t)_{t\geq 0}$. I.e.~we take domain issues into account and use the formulation in the Kolmogorov backwards setting worked out by Grothaus and Stilgenbauer. The method enables us to explicitly compute the constants determining the exponential convergence rate to equilibrium of $(T_t)_{t\geq 0}$. The identification between $(T_t)_{t\geq 0}$ and the transition semigroup of the process yields exponential ergodicity of the latter. Finally, we apply our results to second order in time stochastic reaction-diffusion and Cahn-Hilliard type equations with multiplicative noise. More generally, we analyze corresponding couplings of infinite-dimensional deterministic with stochastic differential equations.}



%35B40 (Asymptotic behavior of solutions to PDEs,
%\pacs[MSC Classification]{37A25 (Ergodicity, mixing, rates of mixing), 37k45 (Stability problems for infinite-dimensional Hamiltonian and Lagrangian systems), 60H15(Stochastic partial differential equations (aspects of stochastic analysis), 35K57 (reaction-diffusion equations (For diffusion processes and reaction effects)}
\maketitle

\section{Introduction}\label{Intro}
In this article we consider infinite-dimensional Langevin dynamics with multiplicative noise. Sometimes they are also called infinite-dimensional stochastic Hamiltonian systems. They are described via the following infinite-dimensional degenerate stochastic 				differential equation
\begin{equation}\label{eq:sde}
	\begin{aligned}
		\mathrm{d} X_t &=\,K_{21}Q_2^{-1}Y_t\,\mathrm{d}t \\
		\mathrm{d}Y_t &= \sum_{i=1}^\infty \partial_iK_{22}(Y_t)e_i - K_{22}(Y_t)Q_2^{-1}Y_t\,\mathrm{d}t-K_{12}Q_1^{-1}X_t-K_{12}					D\Phi(X_t)\,\mathrm{d}t \\&\qquad+ \sqrt{2K_{22}(Y_t)}\,\mathrm{d}B_t.
	\end{aligned}
\end{equation}
Their state space is the Cartesian product $W=U\times V$ of two real separable Hilbert space $(U,\spdx_U)$ and $(V,\spdx_V)$, 				respectively. Above, $K_{21}$ and $K_{12}$ are linear bounded operators from $V$ to $U$ and from $U$ to $V$, respectively. $K_{22}$ is 	the variable diffusion operator, i.e.~a map into the space of symmetric linear bounded operator from $V$ to $V$. Moreover, $D\Phi$ is 		the gradient of a potential $\Phi:U\rightarrow (-\infty,\infty]$, which is bounded from below. $(B_t)_{t\geq 0}$ is a cylindrical 			Brownian motion on $V$. Finally, $Q_1$ and $Q_2$ are the covariance operators of two non-degenerate infinite-dimensional Gaussian 			measures $\mu_1$ and $\mu_2$ on $U$ and $V$, respectively. The covariance operators and the potential determine the measure $				\mu^{\Phi}$ on the Borel $\sigma$-algebra $\mathcal{B}(W)$ by setting
\[
\mu^{\Phi}\defeq \frac{1}{\int_W e^{-\Phi}\,\mathrm{d}\mu_1}e^{-\Phi}\,\mu_1\otimes\mu_2.
\]
The degeneracy of the equations refers to the fact, that the Brownian motion only appears in the second component. On the space of smooth bounded cylinder functions $\fbs{B_W}$, the associated Kolmogorov backward operator $L^{\Phi}$ admits the following representation
\[
	L^{\Phi}\defeq S-A^{\Phi},
\]
where for $f\in \fbs{B_W}$ and $(u,v)\in W$:
\[
	\begin{aligned}
		Sf(u,v)
		&\defeq \tr\left[K_{22}(v)\circ D_2^2f(u,v)\right]
		+ \sum_{j\in\N} (\partial_jK_{22}(v)D_2f(u,v),e_j)_V \\
		&\qquad- (v,Q_2^{-1}K_{22}(v)D_2f(u,v))_V\quad\text{and}
			\end{aligned}
		\]
		\[
		\begin{aligned}
		A^{\Phi}f(u,v)
		&\defeq (u,Q_1^{-1}K_{21}D_2f(u,v))_U+(D\Phi(u),K_{21}D_2f(u,v))_U\\
		&\qquad-(v,Q_2^{-1}K_{12}D_1f(u,v))_V.
	\end{aligned}
\]
Above, $D_1$ and $D_2$ denote the gradient in the first and second component, respectively, compare \Cref{def:grad_comp}. Hereafter, $(L^{\Phi},\fbs{B_W})$ is also referred to as the infinite-dimensional Langevin operator.
	
$(A^{\Phi},\fbs{B_W})$ and $(S,\fbs{B_W})$ correspond to the antisymmetric and symmetric part of $(L^{\Phi},\fbs{B_W})$ on $L^2(W;\mu^{\Phi})$, respectively. The degeneracy of the equation corresponds to the degeneracy of $L^{\Phi}$ in the sense that the second order differential operator in the definition of $L^{\Phi}$ only acts in the second component. As $A^{\Phi}$ contains first order differential operators in the first component and $S$ only differential operators in the second component, the operator $L^{\Phi}$ is non-sectorial.

The structure of this article is as follows. \Cref{preliminaries} introduces Sobolev spaces w.r.t.~infinite-dimensional  Gaussian measures, corresponding integration by parts formulas and an important Poincar\'{e} inequality for measures of type $\mu^{\Phi}$. 

In \Cref{setting_and_operaots} we state conditions on $K_{21}$, $K_{12}$ and $K_{22}$ in order to show that the operators $(L^{\Phi},\fbs{B_W})$, $(S,\fbs{B_W})$ and $(A^{\Phi},\fbs{B_W})$ are dissipative, symmetric and antisymmetric on $L^2(W;\mu^{\Phi})$, respectively.

Afterwards, in \Cref{Essential m-dissipativity}, we show that the closure of the infinite-dimensional Langevin operator $(L^{\Phi},\fbs{B_W})$ in $L^2(W;\mu^{\Phi})$ generates a strongly continuous contraction semigroup $\sccs$, i.e.~is essentially m-dissipative. We overcome the hindrances of non-sectorality of $(L^{\Phi},\fbs{B_W})$ and the $V$-dependence of the diffusion operator $K_{22}$ by using techniques from \cite{EG21} and  \cite{BG21-wp}. To be precise, we start with $\Phi=0$ and use the invariance properties of $K_{21}$, $K_{12}$ and $K_{22}$ stated in \Cref{def:inf-dim-operators} as well as Assumptions \nameref{ass:inf-dim-elliptic} and \nameref{ass:inf-dim-growth}. We transfer the finite-dimensional essential m-dissipativity result from \cite[Theorem 1.1]{BG21-wp} to our infinite-dimensional setting. For $\Phi\neq 0$ we assume \nameref{ass:pert-pot} and extend the essential m-dissipativity result from \cite[Theorem 4.11]{EG21} to infinite-dimensional Langevin operators with variable diffusion matrix $K_{22}$, see \Cref{thm:ess_diss_Lphi}.
	
Invoking general results from \cite{BBR06_Right_process}, we show in \Cref{sec:process}, that there exists a right process
\[
	\mathbf{M}=(\Omega,\mathcal{F},(\mathcal{F}_t)_{t\geq 0}, (X_t,Y_t)_{t\geq 0},(P_w)_{w\in W})
\]
providing a martingale solution for the infinite-dimensional Langevin operator $L^{\Phi}$. Furthermore, the process has infinite life-time and its transition semigroup $(p_t)_ {t>0}$ coincides in $L^2(W;\mu^{\Phi})$ with the strongly continuous contraction semigroup $\sccs$ generated by $L^{\Phi}$. 
We assume \nameref{ass:inf-dim-diff-bound} and use resolvent methods described in \cite{BBR06} by Lucian Beznea, Nicu Boboc and Michael R\"{o}ckner to show that $\mathbf{M}$ is a $\mu^{\Phi}$-invariant Hunt process with weakly continuous paths.
Under Assumption \nameref{ass:weak_sol}, we show that $\mathbf{M}$ provides an analytic and stochastic weak solution to \eqref{eq:sde} in the sense of \Cref{thm:sto_ana_weak_sol}. 

In Section \Cref{sec:hypo_appl} we use the essential m-dissipativity of $(L^{\Phi},\fbs{B_W})$ on $L^2(W;\mu^{\Phi})$ to show exponential convergence to equilibrium of $\sccs$. The idea is to use the abstract hypocoercivity method as described by Grothaus and Stilgenbauer in \cite{GS14} and \cite{GS16}, summarized in \Cref{sec:abstract-hypoc}. The concept of hypocoercivity was developed by C\'{e}dric Villani in \cite{Villani}. Abstract hypocoercivity concepts for a quantitative description of convergence rates were introduced by Dolbeault, Mouhot and Schmeiser, compare \cite{DMS09} and \cite{DMS15}. In order to treat stochastic (partial) differential equations, Grothaus and Stilgenbauer translated these concepts to the Kolmogorov backwards setting taking domain issues into account, see \cite{GS14} and \cite{GS16}. Using generalized Dirichlet forms and martingale techniques, the Ergodicity and the rate of convergence to equilibrium of finite-dimensional Langevin dynamics with singular potentials were studied in \cite{Gs15sing}. Singular potentials could be treated also by Lyapunov techniques, see \cite{HM17}, \cite{Camrud_2022} and \cite{BGH21}. To study finite-dimensional Langevin dynamics with multiplicative noise, a Lyapunov function approach is also used in \cite{Wu2001}. 
Using coupling methods, finite-dimensional Langevin dynamics and their quantitative contraction rates were studied in \cite{EbAn19}. In the context of sub-exponential convergence rates to equilibrium of finite-dimensional Langevin dynamics with multiplicative noise the aforementioned results were generalized in \cite{BG21-wp}.
In \cite{W17}, infinite-dimensional stochastic Hamiltonian systems were treated. The results obtained there are complementary to ours. There the author used coupling methods and dimension-free Harnack inequalities in order to prove hypercontractivity in a degenerate situation without the necessity to explicitly state the invariant measure.
All of the above mentioned references deal either with a constant diffusion matrix $K_{22}$ or with finite-dimensional spaces. We combine the ideas and methods from \cite{BG21} and \cite{EG21_Pr}, where exponential convergence to equilibrium of finite-dimensional Langevin dynamics with multiplicative noise, respectively infinite-dimensional Langevin dynamics without multiplicative noise was shown. The major ingredients to handle both the infinite dimensionality and the variability of the diffusion matrix are besides the essential m-dissipativity of $L^{\Phi}$, the Poincar\'{e} inequality from \Cref{lem:inf-dim-poincare} and the Assumptions \nameref{ass:inf-dim-v-poincare} and \nameref{ass:inf-dim-u-poincare}. In addition, assuming \nameref{ass:qc-negative-type}, the regularity estimates for a perturbed infinite-dimensional Ornstein-Uhlenbeck operator shown in \Cref{lem:inf-dim-regularity-estimates} are essential. In contrast to \cite{EG21_Pr} we can consider non-convex potentials as described in Assumption \nameref{ass:ess-N_Phi} and \nameref{ass:L_Phi-convex-reg-est}. Moreover, note that the variability of $K_{22}$ drastically affects the symmetric operator $S$ and the infinite-dimensionality prevents us to rely on \cite[Proposition 3.10]{BG21} to bound the auxiliary operator $BS$, compare \Cref{def:b-operator}. Therefore boundedness of $BS$ needs to be shown in a different way, which we do in Proposition \Cref{prop:est_T} taking into account either Assumptions \nameref{ass:inf-dim-matrix-bounded-ev} or \nameref{ass:inf-dim-matrix-bounded-ev_star}.
In summary, assuming \nameref{ass:inf-dim-elliptic}-\nameref{ass:pert-pot}, \nameref{ass:qc-negative-type}, \nameref{ass:inf-dim-matrix-bounded-ev}  or \nameref{ass:inf-dim-matrix-bounded-ev_star}, \nameref{ass:inf-dim-v-poincare}, \nameref{ass:inf-dim-u-poincare} \nameref{ass:ess-N_Phi} and \nameref{ass:L_Phi-convex-reg-est} hold, our main hypocoercivity results reads as follows:
For each $\theta_1\in (1,\infty)$, there is some $\theta_2\in (0,\infty)$, which can be explicitly computed (cp. \Cref{thm:inf-dim-hypoc-applied}), such that
\[
	\left\| T_tf- \mu^{\Phi}(f) \right\|_{L^2(W;\mu^{\Phi})} \leq \theta_1 \mathrm{e}^{-\theta_2 t}\left\| f-\mu^{\Phi}(f)  \right\|			_{L^2(W;\mu^{\Phi})}
\]
for all $f\in L^2(W;\mu^{\Phi})$ and all $t\geq 0$.
Via the identification of $\sccs$ and the transition semigroup of the right process $\textbf{M}$ associated to the equation and the inequality above, we are able to show an $L^2$-exponential ergodicity result for $\textbf{M}$, see \Cref{ergodic}.

In the last section we focus on degenerate stochastic reaction-diffusion and Cahn-Hilliard type equations with multiplicative noise. We emphasis how the results from above can be applied. We translate the classic non-degenerate reaction-diffusion and Cahn-Hilliard type equations in our framework of infinite-dimensional Langevin equations with multiplicative noise. Existence of invariant measures and m-dissipativity in an $L^1$-setting for classic stochastic Cahn-Hilliard type equations were established in \cite{Stannat_Cahn}. Alessandra Lunardi and Guiseppe Da Prato studied Sobolev regularity and m-dissipativity of the classic equations in an $L^2$-setting in \cite{Lunardi14} and for Neumann problems in \cite{Lunardi15}. 
Here, we show essential-m dissipativity in $L^2(W;\mu^{\Phi})$, hypocoercivity of the associated semigroup and existence of stochastically and analytically weak solutions of second order in time stochastic reaction-diffusion and Cahn-Hilliard type equations with multiplicative noise. 
In our degenerate situation we are able to deal with potentials of type
\[
	\Phi=\Phi_1+\Phi_2,\quad\text{where}\quad\Phi_1(u)=\int_{(0,1)}\phi( u(\xi))\,\mathrm{d}\xi,\quad u\in U
\]
and $\Phi_2$ is bounded and two times continuously Fr\'{e}chet-differentiable with bounded first and second order derivative. To fulfill Assumption \nameref{ass:ess-N_Phi} and \nameref{ass:L_Phi-convex-reg-est}  we assume that $\phi$ is convex, continuous differentiable with bounded derivative. Moreover, for both $U$ and $V$ we choose $L^2((0,1);\mathrm{d}\xi)$ for the reaction-diffusion equation and the dual space of the Sobolev space $\lbrace x\in W^{1,2}(0,1)\mid \int_0^1x(\xi)\,\mathrm{d}\xi=0\rbrace$ for the Cahn-Hilliard equation. The operators $K_{22}$, $K_{21}$, $K_{12}$, $Q_1$ and $Q_2$ are determined by suitable powers of minus the second order derivative with Dirichlet boundary condition in the reaction-diffusion and by powers of fourth order derivative with zero boundary condition for the first and third order derivative in the Cahn-Hilliard setting. The Assumptions \nameref{ass:inf-dim-elliptic}--\nameref{ass:inf-dim-u-poincare} (including \nameref{ass:inf-dim-matrix-bounded-ev_star}) translate into inequalities of these powers.

\smallbreak

Our main results are summarized in the following list:
\begin{enumerate}
	\item Essential m-dissipativity of $(L^{\Phi},\fbs{B_W})$ on $L^2(W;\mu^{\Phi})$ is proved, see \Cref{thm:ess_diss_Lphi}.
	\item Existence of a $\mu^{\Phi}$-invariant Hunt process $\textbf{M}$ with weakly continuous paths, providing a stochastically and analytically weak solution to the infinite-dimensional degenerate stochastic differential equation \eqref{eq:sde} is shown in \Cref{thm:sto_ana_weak_sol}. Moreover, we represent the strongly continuous contraction semigroup $\sccs$, generated by $(L^{\Phi},\fbs{B_W})$ on $L^2(W;\mu^{\Phi})$, as the 		transition semigroup of $\textbf{M}$.
	\item Hypocoercivity of infinite-dimensional Langevin dynamisc, with explicitly computable rates of convergence is shown in 				\Cref{thm:inf-dim-hypoc-applied}. The identification of $\sccs$ and the transition semigroup of $\textbf{M}$ enables us to derive 			$L^2$-exponential ergodicity of the $\textbf{M}$, see \Cref{ergodic}.
	\item A formulation of degenerate stochastic reaction-diffusion and Cahn-Hilliard type equations in the context of \eqref{eq:sde} is given in \Cref{sec:examples}. In addition, we establish conditions under which we have a corresponding m-dissipative generator, a associated $\mu^{\Phi}$-invariant Hunt process with weakly continuous paths, providing a stochastically and analytically weak solution and hypocoercivity.
\end{enumerate}

\section{Preliminaries}\label{preliminaries}
Let $H$ be a real separable Hilbert space with inner product $\spdx_H$, $\bs{H}$ the corresponding Borel-$\sigma$-algebra, and $\mu$ a centered non-degenerate Gaussian measure on the measurable space $(H,\bs{H})$. 
We denote the set of all linear bounded operators on $H$ by $\lop{H}$, the subset of positive semi-definite symmetric operators by $\lopp{H}$, and the set of operators additionally being of trace class by $\loppt{H}$.
The set of Hilbert-Schmidt operators on $H$ is denoted by $\lophs{H}$. If $\tilde{H}$ is another Hilbert space we denote by $\lop{H,\tilde{H}}$
the set of all linear bounded operators from $H$ to $\tilde{H}$.

The covariance operator corresponding to $\mu$ is denoted by $Q$ and is an element of $\loppt{H}$. It is injective and positive-definite as $\mu$ is non-degenerate. Since it is symmetric and of trace class, there is a complete orthonormal basis $B_H=\s{e}$ in $H$ consisting of eigenvectors of $Q$ with corresponding positive eigenvalues $\s{\lambda}\in\ell^1(\N)$, where $\ell^{\tau}(\N)\defeq \lbrace(a_n)_{n\in\N}\subseteq \R\mid \sum_{n=1}^{\infty}\abs{a_n}^{\tau}<\infty \rbrace$ for $\tau\in (0,\infty)$. W.l.o.g.~we can assume that the eigenvalues are decreasing to zero. Due to positivity of $Q$, there exists an inverse operator defined on $Q(H)$ that satisfies $Q^{-1}e_n=\frac1{\lambda_n}e_n$, for all $n\in\N$. Moreover, we can define $Q^{\theta}$, $\theta\in \R$ on $\lin{e_n:n\in\N}$ characterized by $Q^{\theta}e_n={\lambda^{\theta}_n}e_n$, for all $n\in\N$.
\begin{defn}\label{def:finitely_based}
	For each $n\in\N$, define $H_n\defeq\lin{e_1,\dots,e_n}$ and denote the orthogonal projection from $H$ to $H_n$ by $P_n$, with the 			corresponding coordinate map $p_n:H\to\R^n$. This means for all $x\in H$:
	\[
		P_n(x)\defeq \sum_{k=1}^n (x,e_k)_H e_k	\quad\text{ and }\quad
		p_n(x) \defeq ((x,e_1)_H,\dots,(x,e_n)_H).
	\]	
	Let $k\in\N\cup\lbrace\infty\rbrace$. Denote by $C_b^k(\R^{n})$ the space of bounded k-times differentiable real-valued functions on $		\R^n$ with bounded derivatives and by $\norm{\cdot}_{\infty}$ the sup norm. Let $\mu^n$ be the image measure of $\mu$ under $p_n$. Then define
	\[
	\begin{aligned}
		\mathcal{F}C_b^k(B_H,n)
		&\defeq \{ f:H\to\R \mid f= \varphi\circ p_n
			\text{ for some } \varphi\in C_b^k(\R^{n}) \}\quad\text{and} \\
		\mathcal{F}C_b^k(B_H)
		&\defeq \bigcup_{n\in\N} \mathcal{F}C_b^k(B_H,n). 
	\end{aligned}
	\]
	It is well-known, see for example \cite[Lemma~2.2]{Lunardi14}, that $\fbs{B_H}$ is dense in $L^2(H,\mu)$. 
\end{defn}

Via \cite[Corollary~1.19]{dP06}, we obtain the following characterization for the measure $\mu^n$:
\begin{lem}\label{lem:n-dim-gaussian}
	Let $n\in\N$, and $x_1,\dots, x_n\in H$. The image measure of $\mu$ under the map
	\[
		H\ni x\mapsto ((x,x_1)_H,\dots,(x,x_n)_H)\in \R^n
	\]
	is the centered $n$-dimensional Gaussian measure with covariance matrix
	 \[
	 	((Qx_i,x_j)_H)_{1\leq i,j\leq n}.
	 \]
	 In particular, $\mu^n$ is the centered $n$-dimensional Gaussian measure with covariance matrix $Q_{n}\defeq\diag(\lambda_1,\dots,\lambda_n)$.
\end{lem}

This directly implies the following:
\begin{corollary}\label{cor:gaussian-integral}
	For any $x_1,x_2\in H$, it holds that
	\[
		\int_H (x,x_1)_H (x,x_2)_H\,\mu(\mathrm{d}x) = (Q x_1,x_2)_H.
	\]
	Moreover for all $k\in\N$, $H\ni x\mapsto\|x\|_H^k\in\R$ is $\mu$-integrable. In particular it holds 
	\[
		\int_H \|x\|_H^2\,\mathrm{d}\mu
		= \int_H \sum_{n\in\N} (x,e_n)_H^2 \,\mathrm{d}\mu
		= \sum_{n\in\N} (Qe_n,e_n)_H
		= \sum_{n\in\N} \lambda_n <\infty
	\]
	due to monotone convergence and the fact that $Q$ is trace class.
\end{corollary}

\begin{defn}
	Let $f:H\to\R$ be a Fr\'{e}chet-differentiable function, then we denote its Fr\'{e}chet derivative at the point $x$ by $Df(x)\in H'$, where $H'$ denotes the 				topological dual space of $H$.
	If $f$ is a twice Fr\'{e}chet-differentiable function, then its second order Fr\'{e}chet derivative at the point $x$ is denoted by $D^2f(x)\in\lop{H;H'}$. By 	identifying $H$ with its dual via the Riesz isomorphism, we can interpret $Df(x)$ as an element of $H$ and $D^2f(x)$ as an element of 		$\lop{H}$. For $i,j\in\N$, we denote the partial derivative of $f$ in direction of $e_i$ at the point $x$ by $\partial_if(x)=(Df(x),e_i)$ and the second order partial derivative of $f$ in direction of $e_i$ and $e_j$ by $\partial_{ij}f(x)=\partial_{ji}f(x)=(D^2f(x)e_i,e_j)$, since 			$D^2f(x)$ is symmetric, see for example \cite[(8.12.2)]{D69}.
\end{defn}

\begin{remark}
	By definition, it holds that $Df(x)=\sum_{n\in\N}\partial_if(x) e_i$ for all Fr\'{e}chet-differentiable $f:H\to\R$. If $f=\varphi\circ 			p_n$ for some $n\in\N$, $\varphi\in C_b^\infty(\R^n)$, then the chain rule implies $Df(x)=\sum_{i=1}^n \partial_i\varphi(p_n(x))e_i\in 	H_n$ for all $x\in H$.
\end{remark}

\begin{remark}\label{remark:sobolev}
	Given $p\in [1,\infty)$ and a linear operator $(F,D(F))$, with $H_n\subseteq D(F)$ for all $n\in\N$. Suppose 
	\[
		FD:\fbs{B_H}\rightarrow L^p(H;\mu;H)
	\] 	
	is closable. In this situation we denote by $W^{1,p}_F(H;\mu)$ the domain of the closure of $FD$. In abuse of notation we denote the 		closure of $FD$ again by $FD$. If $F=\operatorname{Id}$ we simply write $W^{1,p}(H;\mu)$.\\
	One can show, compare \cite[Section 2]{dPL14} and \cite[Section 2]{EG21}, that for $\theta\in \R$ the operator $Q^{\theta}D:		\fbs{B_H}\rightarrow L^p(H;\mu;H)$ is closable. In particular for $f\in W^{1,p}_{Q^{\theta}}(H;\mu)$ it is reasonable to set
	\[
		\partial_if\defeq(Q^{\theta}Df,d_i)\frac{1}{\lambda_i^{\theta}}\in L^2(H;\mu).
	\]
	Moreover, it is easy to show that $\tilde{\theta}\in [\theta,\infty)$ implies
	\[
		 W^{1,p}_{Q^{\theta}}(H;\mu) \subseteq W^{1,p}_{Q^{\tilde{\theta}}}(H;\mu).
	\]
\end{remark}

Using that notation, we can state an integration by parts formula for the measure $\mu$. It follows by using \Cref{lem:n-dim-gaussian}, the classical integration by parts formula in $\R^n$ for the Lebesgue measure as well as the remarks in \cite[Section~2]{Lunardi14}.

\begin{lem}\label{lem:inf-dim-ibp-fbs}
	Let $\theta\in [0,1]$, $p,q\in (1,\infty)$ with $\frac{1}{p}+\frac{1}{q}=1$ and $f\in W^{1,p}_{Q^{\theta}}(H;\mu)$, $g\in W^{1,q}			_{Q^{\theta}}(H;\mu)$.
	Then
	\[
		\int_H \partial_i f g\,\mathrm{d}\mu= -\int_H f\partial_i g\,\mathrm{d}\mu+ \int_H (x,Q^{-1}e_i)_H fg\,\mathrm{d}\mu.
	\]
\end{lem}
\begin{defn}
	Suppose $\Phi: H\rightarrow (-\infty,\infty]$ is measurable and bounded from below, s.t.~$\int_He^{-\Phi}\,\mathrm{d}\mu>0$. For such 		$\Phi$ we consider the measure $\mu^{\Phi}\defeq \frac{1}{\int_He^{-\Phi}\,\mathrm{d}\mu}e^{-\Phi}\mu$ and set $\mu^0\defeq \mu$ as 		well as
	\[
		\mu^{\Phi}(f)\defeq \int_Hf\,\mathrm{d}\mu^{\Phi}\quad\text{for}\quad f\in L^1(H;\mu^{\Phi}).
	\]
\end{defn}
As in \cite[Chapter~2.2]{Lunardi14} we can extend the integration by parts formula for measures of type $\mu^{\Phi}$.
\begin{lem}\label{lem:inf-dim-ibp-fbs-Pot}
	Let $\theta\in [0,\infty)$ and consider a potential $\Phi: H\rightarrow (-\infty,\infty]$ which is bounded from below and in $W^{1,2}				_{Q^{\theta}}(H;\mu)$. Then for $f,g\in \mathcal{F}C_b^1(B_H)$ and $i\in\N$, it holds the integration by parts formula
	\[
		\int_H \partial_ifg\,\mu^{\Phi}=-\int_H f\partial_ig\,\mathrm{d}\mu^{\Phi}+\int_H (x,Q^{-1}e_i)_U fg\,\mathrm{d}\mu^{\Phi}+\int_H\partial_i\Phi fg\,\mathrm{d}\mu^{\Phi}.
	\]
\end{lem}

\begin{cond}{$\Phi 0$}\label{hypoapproximationPotential}
	There exist two sequences $(\Phi_{n})_{n\in\N}$ and $(\Phi_{(m,n)})_{(m,n)\in\N^2}$ of convex functions from $H$ to $\R$ such that for $\mu$-almost all $x\in H$ and for all $(m,n)\in\N^2$
	\begin{enumerate}[(i)]
		\item  $ \lim_{m\rightarrow \infty}\Phi_{(m,n)}(x)=\Phi_{n}(x)$  and $-\infty<\inf_{y\in H}\Phi(y)\leq\Phi_{(m,n)}(x)$ as well as\\ $\lim_{(m,n)\rightarrow \infty}\Phi_{(m,n)}(x)\defeq \lim_{n\rightarrow \infty} \lim_{m\rightarrow \infty}\Phi_{(m,n)}(x)=\Phi(x)$.
		\item $\Phi_{(m,n)}$ is Fr\'{e}chet-differentiable with Lipschitz continuous gradient.
		\item If additionally $\Phi\in W^{1,2}_{Q^{\theta}}(H,\mu)$ for some $\theta\in [0,\infty)$, then for all $i\in\N$
		\[
		\lim_{(m,n)\rightarrow 0}\norm{(Q^{\theta}D\Phi_{(n,m)}-Q^{\theta}D\Phi,e_i)_H}_{L^2(H;\mu)}=0.
		\]
	\end{enumerate}
\end{cond}
\begin{remark}\label{rem:weak_conv_measure}
	Suppose Item (i) of Assumption \nameref{hypoapproximationPotential} holds true. Then $e^{-\Phi_{(m,n)}}$ and $e^{-\Phi_{n}}$ are 			bounded by $e^{-\inf_{y\in H}\Phi(y)}$. Since we also have pointwise convergence we obtain by the theorem of dominated convergence  	\[
	\lim_{(m,n)\rightarrow \infty}\mu^{\Phi_{(m,n)}}(f)=\mu^{\Phi}(f)\quad \text{for all}\quad f\in L^1(H;\mu^{\Phi}).
	\]
\end{remark}
The next lemma provides a general Poincar\'{e} inequality for $\mu^{\Phi}$. The result is essential for our further analysis. 

\begin{lem}\label{lem:inf-dim-poincare}
	Suppose $\Phi=\Phi_1+\Phi_2$, where $\Phi_1:H\rightarrow (-\infty,\infty]$ fulfills Item (i) and (ii) from Assumption \nameref{hypoapproximationPotential} 		and $\Phi_2:H\rightarrow \R$ is measurable with $\norm{\Phi_2}_{osc}\defeq\sup_{x\in H}\Phi_2(x)-\inf_{x\in H}\Phi_2(x)<\infty$. Then 
	\[
		\int_H (QDf,Df)_H\,\mathrm{d}\mu^{\Phi}
		\geq \frac{\lambda_1}{e^{\norm{\Phi_2}_{osc}}}\int_H \left(f-\mu^{\Phi}(f) \right)^2\,\mathrm{d}\mu^{\Phi},\quad\text{for all}\quad f\in\fbs{B_H}.
	\]
\end{lem}
\begin{proof}
	For $\Phi_2=0$ we approximate $\Phi$ by $(\Phi_{(m,n)})_{(m,n)\in\N^2}$ and use \cite[Proposition 3.12]{EG21_Pr} to get the 					desired Poincar\'{e} inequality for $\mu^{\Phi_{(m,n)}}$ with Poincar\'{e} constant $\lambda_1$ independent of $(m,n)\in\N^2$. The observation		of the previous remark shows the promised Poincar\'{e} inequality for the measure $\mu^{\Phi_1}		$. Using that $\int_H \left(f-\mu^{\Phi}(f) \right)^2\,\mathrm{d}\mu^{\Phi}\leq \int_H \left(f-c \right)^2\,\mathrm{d}\mu^{\Phi}$ for all $c\in\R$ we can conclude with the Poincar\'{e} inequality for $\mu^{\Phi_1}$
	\[
		\begin{aligned}
		\int_H \left(f-\mu^{\Phi}(f) \right)^2\,\mathrm{d}\mu^{\Phi}&\leq \int_H \left(f-\mu^{\Phi_1}(f) \right)^2\,\mathrm{d}\mu^{\Phi}\\
		&\leq e^{-\inf_{x\in H}\Phi_2(x)}\frac{\int_He^{-\Phi_1}\,\mathrm{d}\mu}{\int_He^{-\Phi}\,\mathrm{d}\mu}\int_H \left(f-\mu^{\Phi_1}(f) \right)^2 \,\mathrm{d}\mu^{\Phi_1}\\
		&\leq \frac{e^{-\inf_{x\in H}\Phi_2(x)}}{\lambda_1}\frac{\int_He^{-\Phi_1}\,\mathrm{d}\mu}{\int_He^{-\Phi}\,\mathrm{d}\mu}\int_H (QDf,Df)_H\,\mathrm{d}\mu^{\Phi_1}\\
		&\leq \frac{e^{\norm{\Phi_2}_{osc}}}{\lambda_1}\int_H (QDf,Df)_H\,\mathrm{d}\mu^{\Phi}.
			\end{aligned}
	\]
\end{proof}

\section{Setting and Operators}\label{setting_and_operaots}

Let $U$ and $V$ be two real separable Hilbert spaces with inner products $\spdx_U$ and $\spdx_V$, respectively. Denote by $\mu_1$ and $\mu_2$ two non-degenerate Gaussian measures $U$ and $V$, respectively.
Let $Q_i$ be the covariance operator of $\mu_i$, $i=1,2$ with corresponding basis of eigenvectors $B_U=\s{d}$ and $B_V=\s{e}$ and positive eigenvalues $(\lambda_{1,k})_{k\in\N}$ and $(\lambda_{2,k})_{k\in\N}$, respectively. W.l.o.g.~we assume that $(\lambda_{1,k})_{k\in\N}$ and $(\lambda_{2,k})_{k\in\N}$ are decreasing to zero. The respective orthogonal projections to the induced subspaces $U_n$ and $V_n$ are denoted by $P_n^U$, $p_n^U$, $P_n^V$ and $p_n^V$, respectively.  Furthermore we fix a potential $\Phi:U\rightarrow (-\infty,\infty]$ and assume that for $\alpha\in [0,\infty)$ the following assumption is valid. 

\begin{cond}{$\Phi^{\alpha}$}\label{ass:Phi_general}
	The potential $\Phi:U\rightarrow (-\infty,\infty]$ fulfills $\int_U e^{-\Phi}\,\mathrm{d}\mu=1$, is bounded from below by zero, and such that $\Phi\in W_{Q_1^{\alpha}}^{1,2}(U;\mu_1)$. 
\end{cond}
Note that all considerations below also work for potentials which are bounded from below and in $W_{Q_1^{\alpha}}^{1,2}(U;\mu_1)$.
\begin{defn}
	Define the real separable Hilbert space $W=U\times V$ with the canonically defined inner product $\spdx_W\defeq \spdx_U+\spdx_V$. Moreover, define the measure $\mu_1^{\Phi}\defeq e^{-\Phi}\mu_1$ on $\mathcal{B}(U)$ and let $\mu^{\Phi}\defeq \mu_1^{\Phi}\otimes\mu_2$ be the product measure on $\bs{W}=\bs{U}\otimes\bs{V}$.
	We set $\mu=\mu^0$. Due to \cite[Theorem~1.12]{dP06}, $\mu$ is a centered Gaussian measure with covariance $Q$ defined by 					$Q(u,v)=(Q_1u,Q_2v)$.
	We set
	\[
		B_W\defeq \{ (d_n,0)\mid n\in\N \} \cup \{ (0,e_n)\mid n\in\N \} \subseteq W
	\]
	as well as
	\[
		\mathcal{F}C_b^k(B_W,n)\defeq \{ f:W\to\R \mid f(u,v)= \varphi(p_n^U(u),p_n^V(v)) 
		\text{ with } \varphi\in C_b^k(\R^{n}\times\R^n) \}
	\]
	for $k\in\N \cup \lbrace\infty\rbrace$. The space $\mathcal{F}C_b^k(B_W)$ is defined analogously to \Cref{def:finitely_based}. Note that $\fbs{B_W}$ is dense in $L^2(W;\mu^{\Phi})$.\\
	Furthermore we set $\mu^n\defeq \mu_1^n\otimes\mu_2^n$, with $\mu_i$ being centered Gaussian measures on $\R^n$ with diagonal 				covariance matrices $Q_{i,n}$.
\end{defn}

\begin{defn}\label{def:grad_comp}
	For Fr\'{e}chet-differentiable $f:W\to\R$ and all $w=(u,v)\in W$ set
	\begin{align*}
		D_1f(w)&\defeq \sum_{n\in\N} (Df(w),(d_n,0))_W d_n\in U,\quad  \partial_{i,1} f(w) \defeq (D_1f(w),d_i)_U\quad\text{and}\\
		D_2f(w)&\defeq \sum_{n\in\N} (Df(w),(0,e_n))_W e_n\in V,\quad\partial_{i,2} f(w) \defeq (D_2f(w),e_i)_V.
	\end{align*}
	The second order derivatives and partial derivatives are defined analogously.
\end{defn}

Now we can introduce the differential operators $S$, $A^{\Phi}$ and $L^{\Phi}$ on $\fbs{B_W}$.

\begin{defn}\label{def:inf-dim-operators}
	Let $K_{12}\in\lop{U;V}$ and set $K_{21}\defeq K_{12}^*\in\lop{V;U}$.
	Let the map $K_{22}:V\to\lopp{V}$ be Fr\'{e}chet-differentiable with $DK_{22}(v)\in\lop{V;\lop{V}}$ and partial derivatives $\partial_iK_{22}(v)=(DK_{22}(v))(e_i)\in\lop{V}$ for each $v\in V$.		
	Moreover, assume there is a strictly increasing sequence $\sk{m}$ in $\N$ such that for each $n\leq m_k$, it holds that	
	\[
		\begin{aligned}
			&K_{12}(U_n)\subseteq V_{m_k},\quad K_{21}(V_n)\subseteq U_{m_k},\quad K_{22}(v)(V_n)\subseteq V_{m_k}\quad\text{and}\\ 					&K_{22}(v)|_{V_n}= K_{22}(P_{m_k}^V(v))|_{V_n}\quad \text{ for all 	}v\in V.
		\end{aligned}
	\]
	Moreover assume that for each $k\in\N$, there is a constant $M_k\in(0,\infty)$ such that
	\[
		\begin{aligned}
			\sup_{v\in V_{m_k}} \|K_{22}(v)\|_{\lop{V_{m_k}}} &\leq M_k
			\quad\text{ and }\\
			\|\partial_iK_{22}(v)\|_{\lop{V_{m_k}}} &\leq M_k(1+\|v\|_{V_{m_k}})
			\quad\text{ for all }v\in V_{m_k}, 1\leq i\leq m_k.
		\end{aligned}
	\]	
	In the following we set $m^K(n)\defeq \min_{k\in\N} \{m_k: m_k\geq n \}$.	
	Then define the differential operators $(S,\fbs{B_W})$ and $(A^{\Phi},\fbs{B_W})$ on $L^2(W;\mu^{\Phi})$ by
	\[
		\begin{aligned}
			Sf(u,v)
			\defeq &\tr\left[K_{22}(v)\circ D_2^2f(u,v)\right]
			+ \sum_{j\in\N} (\partial_jK_{22}(v)D_2f(u,v),e_j)_V \\
			&\qquad- (v,Q_2^{-1}K_{22}(v)D_2f(u,v))_V
		\end{aligned}
	\]
	and
	\[
		\begin{aligned}
			A^{\Phi}f(u,v)
			&\defeq (u,Q_1^{-1}K_{21}D_2f(u,v))_U+(D\Phi(u),K_{21}D_2f(u,v))_U\\
			&\qquad-(v,Q_2^{-1}K_{12}D_1f(u,v))_V,
		\end{aligned}
	\]
	respectively.
	Finally, $(L^{\Phi},\fbs{B_W})$ is defined via $L^{\Phi}\defeq S-A^{\Phi}$. We set $L\defeq L^0$.
\end{defn}

\begin{remark}\label{rem:inf-dim-op-well-defined}
	The invariance assumptions made on $K_{12}$, $K_{21}$ and $K_{22}$ ensure that $S$ and $A^{\Phi}$ are well-defined on $\fbs{B_W}$. 			Indeed, let $f\in\fbs{B_W,n}$ for some $n\in\N$ with corresponding $\varphi\in C_b^\infty(\R^n\times\R^n)$. By trivially extending $		\varphi$ if necessary, we can assume $m^K(n)=n$. Then $Q_1^{-1}K_{21}D_2f(u,v)\in U_n$, $Q_2^{-1}K_{12}D_1f(u,v)\in V_n$, and $Q_2^{-1}K_{22}		(v)D_2f(u,v)\in V_n$ for all $(u,v)\in W$. Therefore, these maps are uniformly bounded in $(u,v)$ due to uniform boundedness of 			$K_{22}:V_{n}\rightarrow \mathcal{L}(V_{n})$ and the fact that all derivatives of $f$ are bounded. Moreover by these invariance 			properties we can interpret  
	\[
		(D\Phi(u),K_{21}D_2f(u,v))_U=\sum_{i=1}^n \partial_i\Phi(u)(d_i,K_{21}D_2f(u,v))_U
	\]even though we don't know if $\Phi\in W^{1,2}(U;\mu_1)$.
	By the observation that all sums appearing in the definition of $S$ are finite, the fact that $\partial_i\Phi\in L^2(\mu^{\Phi})$, as 		well as \Cref{cor:gaussian-integral}, it follows that $Sf(u,v)=Sf(P_n^U u, P_n^V v)$, $A^{\Phi}f(u,v)=A^0f(P_n^U u, P_n^V)+					(D\Phi(u),K_{21}D_2f(P_n^U u, P_n^V))_U$ and $Sf,\,A^{\Phi}f\in L^2(\mu^{\Phi})$. Since $\norm{\cdot}_U^2$ and $\norm{\cdot}_V^2$ are 		integrable with respect to $\mu$, the arguments above also show that $Lf$ is finitely based and in $L^2(\mu)$ for all $f\in \fbs{B_W}$
\end{remark}

\begin{lem}\label{lem:inf-dim-op-decomp}
	On $L^2(W;\mu^{\Phi})$, $(S,\fbs{B_W})$ is symmetric and negative semi-definite, $(A^{\Phi},\fbs{B_W})$ is antisymmetric and 				therefore $(L^{\Phi},\fbs{B_W})$ is dissipative. For $f,g\in\fbs{B_W}$, we have the representation
	\[
		(L^{\Phi}f,g)_{L^2(\mu^{\Phi})} = - \int_W (D_2 f, K_{22} D_2g)_V - (D_1 f, K_{21}D_2 g)_U + (D_2 f, K_{12}D_1 g)_V \,\mathrm{d}			\mu^{\Phi}.
	\]
\end{lem}
\begin{proof}
	Let $f,g\in \fbs{B_W}$. As in \Cref{rem:inf-dim-op-well-defined}, we can assume $f,g\in\fbs{B_W,n}$ for some $n\in\N$ with $n=m^K(n)$. 	For any $(u,v)\in W$, it holds that
	\[
		Q_1^{-1}K_{21}D_2f(u,v)
		= \sum_{k=1}^n \partial_{k,2}f(u,v) Q_1^{-1}K_{21}e_k
		= \sum_{k,\ell=1}^n \partial_{k,2}f(u,v) (K_{21}e_k,d_\ell)_U Q_1^{-1}d_\ell.
	\]
	Using \Cref{lem:inf-dim-ibp-fbs-Pot}, we obtain
	\[
		\int_W \left((u,Q_1^{-1}d_\ell)_U +\partial_l\Phi\right)\partial_{k,2}fg\,\mathrm{d}\mu^{\Phi}
		= \int_W (g\partial_{\ell,1}\partial_{k,2}f+ \partial_{k,2}f \partial_{\ell,1}g) \,\mathrm{d}\mu^{\Phi},
	\]
	which shows that
	\[
	\begin{aligned}
		&\qquad\left((u,Q_1^{-1}K_{21}D_2f)_U+(D\Phi(u),K_{21}D_2f)_U,g \right)_{L^2(\mu^{\Phi})}\\
		 &= \int_W \left(K_{21}D_2f,D_1g \right)_U\,\mathrm{d}\mu^{\Phi}
			+ \sum_{k,\ell=1}^n(K_{21}e_k,d_\ell)_U(g,\partial_{\ell,1}\partial_{k,2}f)_{L^2(\mu^{\Phi})}.		
	\end{aligned}
	\]
	Similarly, it holds that
	\[
	\begin{aligned}
		\quad\left((v,Q_2^{-1}K_{12}D_1f)_V,g \right)_{L^2(\mu^{\Phi})}
		&= \int_W \left(K_{12}D_1f,D_2g \right)_V\,\mathrm{d}\mu^{\Phi}\\
			&\qquad+ \sum_{k,\ell=1}^n(K_{12}d_\ell,e_k)_V(g,\partial_{k,2}\partial_{\ell,1}f)_{L^2(\mu^{\Phi})}.
			\end{aligned}
	\]
	$K_{12}^*=K_{21}$ implies that
	\[
		(A^{\Phi}f,g)_{L^2(\mu^{\Phi})} = \int_W (D_2 f, K_{12}D_1 g)_V-(D_1 f, K_{21}D_2 g)_U\,\mathrm{d}\mu^{\Phi}.
	\]
	In particular $(A^{\Phi}f,f)_{L^2(\mu^{\Phi})}=0$. Now consider the operator $S$. As before, we have
	\[
		Q_2^{-1}K_{22}(v)D_2f(u,v)= \sum_{i,j=1}^n \partial_{i,2}f(u,v)(K_{22}(v)e_i,e_j)_V Q_2^{-1}e_j
	\]
	for all $(u,v)\in W$. Due to the assumptions on $K_{22}$, the maps
		\[
			(u,v)\mapsto \partial_{i,2}f(u,v)(K_{22}(v)e_i,e_j)_V
		\]
	are elements of $\mathcal{F}C_b^1(B_W,n)$ for all $1\leq i,j\leq n$, see \Cref{rem:inf-dim-op-well-defined}. Therefore, the 				application of the integration by parts formula is possible and yields
	\[
	\begin{aligned}
		\int_W (v,Q_2^{-1}e_j) \partial_{i,2}f(K_{22}e_i,e_j)_V g\,\mathrm{d}\mu^{\Phi}
		&= \int_W \partial_{j,2}\partial_{i,2}f(K_{22}e_i,e_j)_V g \,\mathrm{d}\mu^{\Phi} \\
		&\qquad + \int_W \partial_{i,2}f(\partial_j K_{22}e_i,e_j)_V g \,\mathrm{d}\mu^{\Phi} \\
		&\qquad + \int_W \partial_{i,2}f(K_{22}e_i,e_j)_V\partial_{j,2}g\,\mathrm{d}\mu^{\Phi}.
	\end{aligned}
	\]
	Summing up over $i$ and $j$, we get
	\[
		\sum_{i,j=1}^n \int_W \partial_{i,2}f(K_{22}e_i,e_j)_V\partial_{j,2}g\,\mathrm{d}\mu^{\Phi}
		= \int_W (K_{22}D_2f,D_2g)_V\,\mathrm{d}\mu^{\Phi}
	\]
	and
	\[
		\sum_{i,j=1}^n \int_W \partial_{j,2}\partial_{i,2}f(K_{22}e_i,e_j)_V g \,\mathrm{d}\mu^{\Phi}= \int_W \tr\left[K_{22} D_2^2f\right] g\,\mathrm{d}\mu^{\Phi}
	\]
	due to pointwise symmetry of $K_{22}$. Therefore, we indeed get
	\[
		(Sf,g)_{L^2(\mu^{\Phi})} = -\sum_{i,j=1}^n \int_W \partial_{i,2}f(K_{22}e_i,e_j)_V\partial_{j,2}g\,\mathrm{d}\mu^{\Phi}
		= - \int_W (D_2 f, K_{22} D_2g)_V\,\mathrm{d}\mu^{\Phi}.
	\]
	Hence, $S$ is symmetric and negative semi-definite since $K_{22}$ is positive semi-definite.
	Hence, all three operators are dissipative on $\fbs{B_W}$.
\end{proof}
\section{Essential m-dissipativity}\label{Essential m-dissipativity}
In the first part of this section we assume $\Phi=0$ and prove that $(L,\fbs{B_W})$ is essentially m-dissipative on $L^2(W;\mu)$. Since dissipativity was shown in \Cref{lem:inf-dim-op-decomp}, it remains to show that $(\operatorname{Id}-L)(\fbs{B_W})$ is dense in $L^2(W;\mu)$.
Since $\fbs{B_W}$ is dense in $L^2(W;\mu)$, it suffices to approximate all such functions. The main idea here is to interpret $L$ for all $m_k$, $k\in\N$, as an operator on the finite-dimensional subspace determined by $\fbs{B_W,m_k}$, which is possible due to \Cref{rem:inf-dim-op-well-defined}.
In that case, we can apply the recent finite-dimensional m-dissipativity result from \cite[Thm.~1.1]{BG21-wp}.

\begin{defn}\label{def:inf-dim-fin-dim}
	Fix $n\in\N$ such that $n=m^K(n)$. Then we define
	\[
	\begin{aligned}
		K_{12,n}&\defeq ((K_{12}d_i,e_j)_V)_{ij},
		\quad
		K_{21,n}\defeq (K_{12,n})^*
		,\quad \\
		\Sigma_n(y)&\defeq \left(\left(K_{22}\left(\sum_{k=1}^n y_ke_k \right)e_i,e_j\right)_V\right)_{ij}
		\quad\text{ for all }y=(y_1,...,y_n)\in\R^n
		\end{aligned}
	\]
	and denote the entry of $\Sigma_n$ at position $i,j$ by $a_{ij,n}$.
	
	Moreover, define the operators $S_n$, $A_n$ and $L_n$ on $L^2(\mu^n)$ with domain $C_b^\infty(\R^n\times\R^n)$ by
	\[
	\begin{aligned}
	S_nf(x,y) &\defeq \tr[\Sigma_n D_2^2 f](x,y) + \sum_{i,j=1}^{n} \partial_j a_{ij}(y)\partial_{i,2}f(x,y)
				- \langle \Sigma_n(y)Q_{2,n}^{-1}y,D_2 f(x,y) \rangle, \\
	A_nf(x,y) &\defeq \langle K_{12,n}Q_{1,n}^{-1}x,D_2 f(x,y)\rangle - \langle K_{21,n}Q_{2,n}^{-1}y,D_1 f(x,y)\rangle \quad\text{ and }\\
	L_nf &\defeq (S_n-A_n)f,
	\end{aligned}
	\]
	where $\langle \cdot,\cdot\rangle$ and $\abs{\cdot}$ denotes the Euclidean inner product and norm on $\R^n$, respectively.
	 
	Note that these definitions coincide with the structure of operators considered in \cite{BG21-wp}, with the choices $\Theta(x)=\frac12 \langle x,Q_{1,n}^{-1}x \rangle$ and $\Psi(y)=\frac12\langle y,Q_{2,n}^{-1}y\rangle$.
\end{defn}
\begin{lem}
	Let $f\in\fbs{B_W,n}$ for some $n\in\N$ with $n=m^K(n)$, with $f(u,v)=\varphi(p_n^U u,p_n^V v)$ for $\varphi\in C_b^\infty(\R^n\times\R^n)$. Then $Sf(u,v)=S_n\varphi(p_n^U u,p_n^V v)$, and analogue statements hold for $A$ and $L$.
\end{lem}
\begin{proof}
	Follows directly from \Cref{rem:inf-dim-op-well-defined} together with $\nabla\Theta(x)=Q_{1,n}^{-1}x$ and $\nabla\Psi(y)=Q_{2,n}^{-1}		y$.
\end{proof}

Next, we fix some assumptions on $K_{22}$ such that they imply sufficient conditions on each $\Sigma_n$ to ensure that \cite[Thm.~1.1]{BG21-wp} can be applied. In particular, $\Sigma_n$ should satisfy:
\begin{cond}{$\Sigma$}\leavevmode 
\begin{enumerate}[($\Sigma$1)]
	\item $\Sigma_n$ is symmetric and uniformly strictly elliptic, i.e.~there is some $0<c_{\Sigma_n}<\infty$
		such that
		\[
			\langle v,\Sigma_n(y) v\rangle \geq c_{\Sigma_n}^{-1}\cdot |v|^2\quad\text{ for all }v\in\R^{n} \text{ and $\mu_2^n$-almost all }y\in\R^{n}.
		\]
	\item For each $1\leq i,j\leq n$, $a_{ij,n}$ is bounded and locally Lipschitz-continuous.
	\item There are constants $0\leq M<\infty$, $0\leq \beta <1$ such that for all $1\leq i,j,k\leq n$
		\[
			|\partial_k a_{ij,n}(y)|\leq M(1+|y|^\beta)\quad\text{ for $\mu_2^n$-almost all }y\in\R^{n}.
		\]
\end{enumerate}
\end{cond}
\begin{cond}{K1}\label{ass:inf-dim-elliptic}
	Assume that there is some positive-definite $K_{22}^0\in \lopp{V}$ which leaves each $V_{m_k}$ invariant, $k\in\N$, such that
	\[
		(v,K_{22}(y)v)_V \geq (v,K_{22}^0 v)_V
		\qquad\text{ for all }v,y\in V.
	\]
\end{cond}
Above $\sk{m}$ is the sequence from \Cref{def:inf-dim-operators},
\begin{cond}{K2}\label{ass:inf-dim-growth}
	For each $n\in\N$, let $k(n)$ be such that $m_{k(n)}=m^K(n)$.
	Assume that there are sequences $\sk{\beta}$ in $[0,1)$ and $\sk{N}$ in $\R$ such that, for all $n\in\N$,
	\[
	\begin{aligned}
		|(\partial_iK_{22}(v)e_n,e_j)_V|
		&\leq N_{k(n)}(1+\|v\|_{V_{m^K(n)}}^{\beta_{k(n)}})
	\end{aligned}
	\]
	for all $v\in V_{m^K(n)}$, $1\leq i\leq m^K(n)$ and $1\leq j\leq n$.
	
	For $n\in\N$, set $N^K(n)\defeq 2\max\{ N_{k(j)}: 1\leq j\leq n\}$ and
	$\beta^K(n)\defeq \max\{ \beta_{k(n)}:1\leq j\leq n\}$.
\end{cond}

\begin{remark}\label{rem:too-diagonal}
	Assume that $K_{22}(v)$ leaves $V_n$ invariant for all $n\in\N$ and $v\in V$. Using the strengthened invariance properties of $K_{22}$
	it follows quickly that $K_{22}(v)$ is diagonal. I.e.~$K_{22}(v)e_i=\lambda_{22,i}(v)e_i$ for some non-negative continuously differentiable $\lambda_{22,i}:V\to\R$. In that case, \nameref{ass:inf-dim-elliptic} just means 			that each $\lambda_{22,i}$ is bounded from below by a positive constant $\lambda_i^0\in\R$, and \nameref{ass:inf-dim-growth} reduces 		to $|\partial_i\lambda_{22,n}(v)|=|\partial_i\lambda_{22,n}(P_{m^K(n)}^V v)|\leq N_{k(n)}(1+\|P_{m^K(n)}^V v\|_V^{\beta_{k(n)}})$ for 		all $1\leq i\leq m^K(n)$ and $n\in\N$.
\end{remark}
\begin{lem}\label{lem:inf-dim-matrix-properties}
	Let $n\in\N$ with $n=m^K(n)$ and define $\Sigma_n$ as in \Cref{def:inf-dim-fin-dim}.
	If $K_{22}$ satisfies \nameref{ass:inf-dim-elliptic} and \nameref{ass:inf-dim-growth},
	then $\Sigma_n$ fulfills ($\Sigma$1)--($\Sigma$3), i.e.~it is admissible for \cite[Thm.~1.1]{BG21-wp}.
\end{lem}
\begin{proof}
	Set $\Sigma_n^0\in\R^{n\times n}$ analogously to $\Sigma_n$ for $K_{22}^0$. Since $K_{22}^0$ is positive-definite,
	all eigenvalues $\lambda_1^0,\dots,\lambda_n^0$ of $\Sigma_n^0$ are positive, and we define
	$c_n\defeq \min_{i\in\{1,\dots,n\}}\lambda_i>0$. For any $y\in\R^n$, set $\tilde{y}\defeq\sum_{i=1}^n y_ie_i\in V_n$. Then
	\[
		\langle y, \Sigma_n(v) y\rangle
		= (\tilde{y},K_{22}(\tilde{v}) \tilde{y} )_V
		\geq (\tilde{y},K_{22}^0 \tilde{y} )_V
		= \langle y, \Sigma_n^0 y\rangle
		\geq c_n |y|^2
	\]
	for all $y,v\in \R^n$. Hence, ($\Sigma$1) holds with $c_{\Sigma_n}\defeq c_n^{-1}$.
	Due to definition of $K_{22}$, we already have that all entries of $\Sigma_n$ are bounded and continuously differentiable, hence in particular locally Lipschitz. Assume that $j\leq i$, and let $k\in\{1,\dots,n\}$. Then
	\[
	\begin{aligned}
		|\partial_k a_{ij,n}(y)|
		&= |\partial_k (K_{22}(\tilde{y})e_i,e_j)_V|
		= |(\partial_kK_{22}(\tilde{y})e_i,e_j)_V|
		\leq N_{k(i)}(1+\|\tilde{y}\|_{V_n}^{\beta_n}) \\
		&\leq 2N_{k(i)}(1+\|\tilde{y}\|_{V_n}^{\beta^K(n)})
		\leq N^K(n)(1+\|\tilde{y}\|_{V_n}^{\beta^K(n)})
	\end{aligned}
	\]
	by \nameref{ass:inf-dim-growth}, so $\Sigma_n$ satisfies ($\Sigma$3) with constants $M=N^K(n)$ and $\beta=\beta^K(n)$.
\end{proof}

\begin{prop}\label{prop:inf-dim-finite-m-diss}
	Let $n\in\N$ such that $n=m^K(n)$ and let $K_{22}$ satisfy \nameref{ass:inf-dim-elliptic} and \nameref{ass:inf-dim-growth}. Then $(L_n,C_b^\infty(\R^n\times\R^n))$ is essentially m-dissipative on $L^2(\mu^n)$.
\end{prop}
\begin{proof}
	We finally apply \cite[Thm.~1.1]{BG21-wp}.
	Due to \Cref{lem:inf-dim-matrix-properties}, the conditions on $\Sigma_n$ are satisfied,
	and most properties required of $\Theta$ and $\Psi$ as chosen in \Cref{def:inf-dim-fin-dim} are immediate.
	Note in particular that due to the simple structure of $\Theta$, we can apply \cite[Thm.~2.3]{BG21-wp} to obtain $\Theta\in L^2(\R^n;\mu_1^n)$, as well as the analogue for $\Psi$.
	For any $x\in\R^n$, it holds that
	\[
		|\nabla\Theta(x)|^2=\sum_{i=1}^n \frac1{\lambda_{1,i}^2}x_i^2
		\leq \frac1{\lambda_{1,n}^2} |x|^2
	\]
	since $Q_{1,n}=\diag(\lambda_{1,1},\dots,\lambda_{1,n})$, where $(\lambda_{1,i})_{i\in\N}$ is the decreasing sequence of eigenvalues of $Q_1$. Therefore, $\Theta$ satisfies the required growth condition ($\Theta$2) for $N=\lambda_{1,n}^{-1}$ and $\gamma=1<(\beta^K(n))^{-1}$.
	
	As a result, it follows that $(L_n,C_c^\infty(\R^n\times\R^n))$ is essentially m-dissipative on $L^2(\mu^n)$. Since $C_b^\infty(\R^n\times\R^n)$ extends that domain, and $(L_n,C_b^\infty(\R^n\times\R^n))$ is dissipative on $L^2(\mu^n)$ due to \Cref{lem:inf-dim-op-decomp},
	the claim follows.
\end{proof}

Now finally, we are able to prove the first major result.
\begin{thm}\label{thm:inf-dim-ess-m-diss}
	Let $K_{22}$ satisfy \nameref{ass:inf-dim-elliptic} and \nameref{ass:inf-dim-growth}.
	Then $(L,\fbs{B_W})$ is essentially m-dissipative on $L^2(W;\mu)$.
	Furthermore, the strongly continuous contraction semigroup and the resolvent $(R_{\alpha}^L)_{\alpha >0}$ generated by the closure of $(L,\fbs{B_W})$ is sub-Markovian, conservative and $\mu$-invariant. 
\end{thm}
\begin{proof}
	As mentioned above, we only need to show that $(\operatorname{Id}-L)(\fbs{B_W})$ is dense in $L^2(W;\mu)$, since \Cref{lem:inf-dim-op-decomp} provides 	dissipativity of $(L,\fbs{B_W})$. Let $g\in \fbs{B_W}$, then there is some $n\in\N$ such that $g\in \fbs{B_W,n}$. As before, we extend 	$g$ trivially to $\fbs{B_W,m^K(n)}$, so that we can assume $n=m^K(n)$. Let $\varphi_g\in C_b^\infty(\R^n\times\R^n)$ be such that 			$g(u,v)=\varphi(p_n^U u,p_n^V v)$ for all $(u,v)\in W$ and let $\eps>0$. Then
	\[
	\begin{aligned}
		\left\|(\operatorname{Id}-L)f - g \right\|_{L^2(\mu)}^2
		&= \int_W \left((\operatorname{Id}-L)f(P_n^U u, P_n^V v)- g(P_n^U u, P_n^V v) \right)^2 \,\mu(\mathrm{d}(u,v)) \\
		&= \int_{\R^n\times\R^n} \left((\operatorname{Id}-L_n)\varphi_f(x,y)-\varphi_g(x,y)\right)^2\,\mu^n(\mathrm{d}(x,y)) \\
		&= \left\|(\operatorname{Id}-L_n)\varphi_f - \varphi_g \right\|_{L^2(\mu^n)}^2
	\end{aligned}
	\]
	for all $f\in\fbs{B_W,n}$ with corresponding $\varphi_f\in C_b^\infty(\R^n\times\R^n)$.
	Due to \Cref{prop:inf-dim-finite-m-diss}, there is some $h\in C_b^\infty(\R^n\times\R^n)$ such that
	$\|(\operatorname{Id}-L_n)h - \varphi_g \|_{L^2(\mu^n)}<\eps$. Setting $f_h(u,v)\defeq h(p_n^U u,p_n^V v)$ yields $f_h\in\fbs{B_W,n}$ with $\|(\operatorname{Id}-L)f_h - g\|_{L^2(\mu)}<\eps$. Since $\fbs{B_W}$ is dense in $L^2(W;\mu)$,
	this proves that $(\operatorname{Id}-L)(\fbs{B_W})$ is dense in $L^2(W;\mu)$.
	
	By explicit calculation, it can be shown that $(L,\fbs{B_W})$ is an abstract diffusion operator in the sense of \cite[Def.~1.3.5]{Conrad2011}, which implies that it is a Dirichlet operator and that $\sccs$ is sub-Markovian.
	By \Cref{lem:inf-dim-op-decomp}, it follows that $L1=0$ and $\mu(Lf)=0$ for all $f\in\fbs{B_W}$.
	The former implies $T_t1=1$ in $L^2(W;\mu)$ for all $t\geq 0$, while the latter shows that $L$ and therefore $\sccs$ is $\mu$-invariant.
\end{proof}

For the rest of the chapter we consider the situation where $\Phi\neq 0$, i.e.~we work on $L^2(W;\mu^{\Phi})$ and want to show essential m-dissipativity of $(L^{\Phi},\fbs{B_W})$ in $L^2(W;\mu^{\Phi})$. To do that we need to derive some regularity estimates, first.

\begin{lem}\label{lem:reg_K_22}
Let $K_{22}$ satisfy \nameref{ass:inf-dim-elliptic} and \nameref{ass:inf-dim-growth}. For $f\in D(L)$ and $\alpha \in (0,\infty)$, set $g=\alpha f-Lf$. Then $K_{22}^{\frac{1}{2}}D_2f$ exists in $L^2(W;\mu)$ and the following equation hold
	\[
		\int_W\alpha f^2+\norm{K_{22}^{\frac{1}{2}}D_2f}^2_V\dm=\int_W fg\dm.
	\]
In particular
	\begin{equation}\label{varepsilonone}
		\int_W\norm{K_{22}^{\frac{1}{2}}D_2f}^2_V\dm\leq \frac{1}{2} \int_Wf^2+(Lf)^2\dm\quad\text{and}
	\end{equation}
	\begin{equation}\label{varepsilonn2}
		\int_W\norm{K_{22}^{\frac{1}{2}}D_2f}^2_V\dm\leq \frac{1}{4\alpha} \int_Wg^2\dm.
	\end{equation}
\end{lem}
\begin{proof}
	Assume $f\in \fbs{B_W}$ and $g=\alpha f-Lf$. Now Multiply $g=\alpha f-Lf$ with $f$, integrate over $W$ w.r.t.~$\mu$ and use Lemma 			\ref{lem:inf-dim-op-decomp} to obtain the first identity. Rearranging the terms we obtain
	\[
		\int_W\norm{K_{22}^{\frac{1}{2}}D_2f}^2_V\dm=\int_W f(g-\alpha f)\dm=-\int_W fLf\dm\leq \frac{1}{2} \int_Wf^2+(Lf)^2\dm.
	\]
	Moreover, by completing the square we have
	\[
		\int_W\norm{K_{22}^{\frac{1}{2}}D_2f}^2_V\dm=-\int_W \alpha f^2-fg\dm\leq \frac{1}{4\alpha} \int_Wg^2\dm. 
	\]
	By construction $\fbs{B_W}$ is dense in the $(L,D(L))$ graph norm. Hence, $K_{22}^{\frac{1}{2}}D_2f$ exists for $f\in D(L)$ as the $L^2(W;\mu)$-limit of $K_{22}^{\frac{1}{2}}D_2f_n$, where $(f_n)_{n\in\N}$ is the the approximating sequence of $f$ w.r.t~$(L,D(L))$ graph norm. Note that $K_{22}^{\frac{1}{2}}D_2f$ does not depend on the approximating sequence. In particular the (in)equalities above are also valid for $f\in D(L)$.
\end{proof}

Having these regularity estimates at hand we need one last assumption to derive the final essential m-dissipativity result. The assumption is necessary to apply the perturbation argument described in the next theorem.
\begin{cond}{K3}\label{ass:pert-pot}$ $
	 There is an $\alpha\in [0,\infty)$ and a constant $c_{\alpha}\in (0,\infty)$ such that
	\begin{enumerate}[(i)]
		\item $(Q_1^{-\alpha}K_{21}v,Q_1^{-\alpha}K_{21}v)_V\leq c_{\alpha} (K_{22}(\tilde{v})v,v)_V$ for all $\tilde{v}\in V$ and 				$v\in V_n$, $n\in\N$.
		\item Assumption \nameref{ass:Phi_general} is valid and $Q_1^{\alpha}D\Phi$ is in $L^{\infty}(U;\mu_1)$.
	\end{enumerate}
\end{cond}
\begin{thm}\label{thm:ess_diss_Lphi}
	Let $K_{22}$ satisfy \nameref{ass:inf-dim-elliptic} and \nameref{ass:inf-dim-growth} and suppose Assumption \nameref{ass:pert-pot} is valid. Then $D(L)\subseteq D(L^{\Phi})$ with
	\[
		L^{\Phi}f=Lf-(D\Phi,K_{21}D_2f)_U, \quad f\in D(L).
	\]
	Moreover, the infinite-dimensional Langevin operator $(L^{\Phi},\fbs{B_W})$ is essentially m-dissipative on $L^2(W;\mu^{\Phi})$. 			Furthermore, the strongly continuous semigroup/resolvent $\sccs/(R_{\alpha}^{L^{\Phi}})_{\alpha >0}$ generated by the closure of $			(L^{\Phi},\fbs{B_W})$ is sub-Markovian, conservative and $\mu^{\Phi}$-invariant. 
\end{thm}
\begin{proof}
	Using Item $(i)$ from Assumption \nameref{ass:pert-pot} and Inequality \eqref{varepsilonone} there exists $\alpha\in [0,\infty)$ such 		that for all $f\in\fbs{B_W}$
	\[
	\begin{aligned}
		\frac{1}{\lambda_{1,1}^{2\alpha}}\int_W\norm{K_{21}D_2f}^2_U\dm&\leq\int_W\norm{Q_1^{-\alpha}K_{21}D_2f}^2_U\dm\\&\leq\int_W 				c_{\alpha} \norm{K^{\frac{1}{2}}_{22}D_2f}_V\dm\\&\leq \frac{c_{\alpha}}{2} \int_Wf^2+(Lf)^2\dm.
	\end{aligned}
	\]
	As in \Cref{remark:sobolev} we can show that $K_{21}D_2, Q_1^{-\alpha}K_{21}D_2:\fbs{B_W} \to L^2(\mu^{\Phi})$ are closable in 				$L^2(\mu^{\Phi})$ with closures again denoted by $K_{21}D_2$ and $Q_1^{-\alpha}K_{21}D_2$. Since $\fbs{B_W}$ is dense $D(L)$ 				w.r.t.~the $(L,D(L))$ graph norm the estimate above also holds for $f\in D(L)$.
	Let $(f_n)_{n\in\N}\subseteq \fbs{B_W}$ be a sequence converging to $f\in D(L)$ with respect to the $(L,D(L))$ graph norm. Since $\Phi$ 		is bounded 	from below it is easy to check that $(f_n)_{n\in\N}$ converges to $f$ in $L^2(W;\mu^{\Phi})$. Moreover we can estimate
	\[
	\begin{aligned}
		&\int_W(L^{\Phi}f_n-Lf+(D\Phi,K_{21}D_2f)_U)^2\dm^{{\Phi}}\\
		&\leq 2\int_W (Lf_n-Lf)^2\dm^{\Phi} +2 \int_W (Q_1^{\alpha}D\Phi,Q_1^{-\alpha}K_{21}D_2(f_n-f))^2_U\dm^{\Phi}\\
		&\leq 2\big(\int_W (Lf_n-Lf)^2\dm+\frac{c_{\alpha}}{2}\norm{Q_1^{\alpha}D\Phi}^2_{L^{\infty}(\mu_1)}\int_W(f_n-f)^2+(Lf_n-					Lf)^2\dm\big).
	\end{aligned}
	\]
	Hence the sequence $(L^{\Phi}f_n)_{n\in\N}$ converges to $Lf-(D\Phi,K_{21}D_2f)_U$ in $L^2(\mu^{\Phi})$. As the operator $					(L^{\Phi},D(L^{\Phi}))$ is closed we get $D(L)\subseteq D(L^{\Phi})$ and for all $f\in D(L)$ 
	\[
		L^{\Phi}f=Lf-(D\Phi,K_{21}D_2f)_U.
	\]
	By Lemma \ref{lem:inf-dim-op-decomp} we already know that $(L^{\Phi},\fbs{B_W})$ is dissipative. In view of the Lumer-Phillips theorem 	we are left to show the dense range condition.
	For $f\in L^2(W;
	\mu)$ and $\lambda>0$ set
	\[
		T_{\lambda}f=-(D\Phi,K_{21}D_2R_{\lambda}^Lf)_U.
	\]
	We calculate using the Cauchy-Schwarz inequality, the assumption on $\Phi$, Assumption \nameref{ass:pert-pot} and Inequality 				\eqref{varepsilonn2}
	\[
	\begin{aligned}
		\int_W (T_{\lambda}f)^2\dm &=\int_W (D\Phi,K_{21}D_2R_{\lambda}^Lf)_U^2\dm\\
		&\leq \norm{Q_1^{\alpha}D\Phi}^2_{L^{\infty}(\mu_1)}\int_W (Q_1^{-\alpha}K_{21}D_2R_{\lambda}^Lf,Q_1^{-\alpha}K_{21}D_2R_{\lambda}^Lf)_U\dm\\
		&\leq \norm{Q_1^{\alpha}D\Phi}^2_{L^{\infty}(\mu_1)}\int_W c_{\alpha}(K_{22}D_2R_{\lambda}^Lf,D_2R_{\lambda}^Lf)_U\dm\\
		&\leq \norm{Q_1^{\alpha}D\Phi}^2_{L^{\infty}(\mu_1)}\frac{c_{\alpha}}{4\lambda}\int_W f^2\dm.
	\end{aligned}
	\]
	Hence the operator $T_{\lambda}:L^2(\mu)\rightarrow L^2(\mu)$ is well-defined. Moreover, if 
	\[
		\norm{Q_1^{\alpha}D\Phi}^2_{L^{\infty}(\mu_1)}\frac{c_{\alpha}}{4\lambda}<1
	\] we can apply Neumann-Series theorem to get $(\operatorname{Id}-T_{\lambda})^{-1}\in \mathcal{L}(L^2(W;\mu))$. In particular, for such $\lambda$, for all $g\in L^2(W;\mu)$ we find $f\in L^2(\mu)$ with $f-T_{\lambda}f=g$ in $L^2(W;\mu)$. Furthermore there is $h\in D(L)$ with $(\lambda-L)h=f$. Therefore,
	\[
		(\lambda-L^{\Phi})h=(\lambda-L)h+(D\Phi,K_{21}D_2h)_U=f+(D\Phi,K_{21}D_2R_{\lambda}^Lf)_U=f-T_{\lambda}f=g.
	\]
	This yields $L^2(W;\mu)\subseteq(\lambda-L^{\Phi})(D(L))$. Since $L^2(W;\mu)$ is dense in $L^2(W;\mu^{\Phi})$ and $D(L)\subseteq D(L^{\Phi})$ the dense range condition is shown and the essential m-dissipativity is shown. Sub-Markovianity, conservativity, and $\mu^{\Phi}$-invariance is shown as in \Cref{thm:inf-dim-ess-m-diss}. 
\end{proof}

\section{The associated stochastic process}\label{sec:process}

Assume that Assumptions \nameref{ass:inf-dim-elliptic}-\nameref{ass:pert-pot} hold true. In view of \Cref{thm:ess_diss_Lphi} it follows that $(L^{\Phi},\fbs{B_W})$ is essentially m-dissipative on $L^2(W;\mu^{\Phi})$. Furthermore, the strongly continuous contraction semigroup/resolvent $\sccs/(R_{\alpha}^{L^{\Phi}})_{\alpha >0}$ generated by the closure of $(L^{\Phi},\fbs{B_W})$ is sub-Markovian, conservative, and $\mu^{\Phi}$-invariant. By \cite[Theorem~2.2]{BBR06_Right_process} there exists a right process 
\[
	\mathbf{M}=(\Omega,\mathcal{F},(\mathcal{F}_t)_{t\geq 0}, (X_t,Y_t)_{t\geq 0},(P_w)_{w\in W})
\]
such that its transition $(p_{t})_{t >0}$ semigroup and resolvent coincides with $\sccs$ and $(R_{\alpha}^{L^{\Phi}})_{\alpha >0}$ in $L^2(W;\mu^{\Phi})$.
The transition semigroup is defined for $t>0$ and $g\in B_b(W)$ via
\[
 p_tg(w)=\int_{\Omega} g((X_t,Y_t)(\omega))\,\mathrm{d}{P}^{w}(\omega)
\]
and identified with its extension to $L^2(W;\mu^{\Phi})$.

Since $(L^{\Phi},D(L^{\Phi}))$ is conservative, we can apply \cite[Lemma 2.1.14]{Conrad2011} to show that $(X_t,Y_t)_{t\geq 0}$ has infinite life time $P_{\mu^{\Phi}}\defeq \int_{W} P_w\,\mu^{\Phi}(\mathrm{d}w)$ a.s..
\begin{remark}\label{rem:inf-dim-mart-prob}
	As a direct consequence of \cite[Proposition~1.4]{BBR06}, we see that for any $f\in D(L^{\Phi})$,
	\[
		M_t^{[f],L^{\Phi}}\defeq f(X_t,Y_t)-f(X_0,Y_0)-\int_0^ tL^{\Phi}f(X_s,Y_s)\,\mathrm{d}s,\quad t\geq 0
	\]
	is $P_{\mu^{\Phi}}$-integrable and describes an $(\mathcal{F}_t)_{t\geq 0}$-martingale.
	Moreover, if $f^2\in D(L^{\Phi})$ with $L^{\Phi}f\in L^4(\mu^{\Phi})$, then
	\[
		N_t^{[f],L^{\Phi}}\defeq (M_t^{[f],L^{\Phi}})^2 - \int_0^t L^{\Phi}(f^2)(X_s,Y_s)-(2fL^{\Phi}f)(X_s,Y_s)\,\mathrm{d}s,
		\quad t\geq 0
	\]
	describes an $(\mathcal{F}_t)_{t\geq 0}$-martingale as well.
\end{remark}

Now we want to show that $\mathbf{M}$ is a $\mu^{\Phi}$-invariant Hunt process with state space $(W,\mathcal{T})$, where $\mathcal{T}$ denotes the weak topology on $W$. In order to do that we prove that \[F_n\defeq \{w\in W: \|w\|_W\leq n\}\] defines an $\mu^{\Phi}$-nest of $\mathcal{T}$-compact sets, employing the strategy used in \cite[Lemma~3.3]{EG21}. Consequently we assume the following:

\begin{cond}{K4}\label{ass:inf-dim-diff-bound}\leavevmode 
	\begin{enumerate}[(i)]
		\item There exists $\rho\in L^1(W;\mu^{\Phi})$ such that for each $n\in\N$, $\rho_n$ defined via $\rho_n(u,v)\defeq 					\rho(P_n^Uu,P_n^Vv)$ is also in $L^1(W;\mu^{\Phi})$. Moreover the sequence $(\rho_{m^K(n)})_{n\in\N}$ converges to $\rho$ in $L^1(W;						\mu^{\Phi})$ as $n\to\infty$ and 
		\[
			(P_{m^K(n)}^U u,Q_1^{-1}K_{21}P_{m^K(n)}^V v)_U
			- (Q_2^{-1}K_{12}P_{m^K(n)}^U u,P_{m^K(n)}^V v)_V
			\leq \rho_{m^K(n)}(u,v)
		\]
		for all $n\in\N$ and $(u,v)\in W$.
		\item There exist $a,b\in\N$ such that for all $n\in\N$ and $v\in V$
		\[
			\sum_{j=1}^n\norm{\partial_jK_{22}(v)e_j}_V\leq a(1+ \norm{v}^b_V).
		\]
		\item The sequence $(\overline{h}_n)\subseteq L^1(V;\mu_2)$ defined by \[\overline{h}_n(v)\defeq (P_{m^K(n)}^V v,Q_2^{-1}K_{22}(v)P_{m^K(n)}^V 				v)_V\in\R\] converges in $L^1(V;\mu_2)$.
	\end{enumerate}
\end{cond}
%		\item The sequence $(N^K(n))_{n\in\N}$ from \nameref{ass:inf-dim-growth} is in $l^2(\N)$.
\begin{prop}\label{prop:inf-dim-assoc-process}
	Let \nameref{ass:inf-dim-diff-bound} hold. Then $(F_n)_{n\in\N}$ is $\mu^{\Phi}$-nest of $\mathcal{T}$-compact sets for $					(L^{\Phi},D(L^{\Phi}))$.
\end{prop}
\begin{proof}
	By the Theorem of Banach-Alaoglu $F_n$ is $\mathcal{T}$-compact for all $n\in\N$.
	It remains to prove that $\s{F}$ is a $\mu^{\Phi}$-nest.
	For notation purposes, we write $N(u,v)\defeq \|(u,v)\|_W^2$ and $N_n(u,v)=N(P_n^U u,P_n^V v)$ for $n\in\N$.
	Moreover, we only consider those $n\in\N$ that satisfy $n=m^K(n)$, which provide an increasing sequence.
	Using a sequence of cutoff functions, we see that $N_n\in D(L^{\Phi})$ with
	\[
	\begin{aligned}
		&\quad\frac12 L^{\Phi}N_n(u,v)\\
		&=\tr[K_{22}(P_n^V v)]
			+\sum_{j=1}^n (\partial_jK_{22}(v)e_j,P_n^V v)_V
			-(P_n^V v,Q_2^{-1}K_{22}(v)P_n^V v)_V \\
		&\qquad-(D\Phi(u),K_{21}P_n^V v)_U- (P_n^U u,Q_1^{-1}K_{21}P_n^V v)_U
			+(P_n^V v,Q_2^{-1}K_{12}P_n^U u)_V \\
		&\geq \sum_{j=1}^n (\partial_jK_{22}(v)e_j,P_n^V v)_V
			-(P_n^V v,Q_2^{-1}K_{22}(v)P_n^V v)_V
			\\
			&\qquad-(D\Phi(u),K_{21}P_n^V v)_U-\rho_n(u,v).
	\end{aligned}
	\]
	Using the second item from \nameref{ass:inf-dim-diff-bound} we find $a,b\in\N$ with
%	\[ 
%		\left|\sum_{j=1}^n (\partial_jK_{22}(v)e_j,P_n^V v)_V \right|
%		&\leq \sum_{i,j=1}^n |(\partial_jK_{22}(v)e_j,e_i)_V (e_i,v)_V| \\
%		&\leq \sum_{i=1}^n \left(|(\partial_jK_{22}(v)e_j,e_i)_V^2\right)^{\frac{1}{2}} (e_i,v)_V| \\
%		&\leq   2(1+\norm{P_n^V v}_V)\norm{(N_{k(j)})_{j\in\N}}_{l^2(\N)} \norm{P_n^V v}_V\eqdef  h_n(v).
%	\end{aligned}
%	\]
	\[ 
	\begin{aligned}
		\left|\sum_{j=1}^n (\partial_jK_{22}(v)e_j,P_n^V v)_V \right|\leq a(1+ \norm{v}^b_V)\norm{P_n^V v}_V\eqdef  h_n(v).
	\end{aligned}
	\]
	Note that $(h_n)_{n\in\N}$ converge in $L^1(V;\mu_2)$.
	Moreover, by \nameref{ass:pert-pot} it holds for $n,m\in\N$ with $n<m$
	\[
	\begin{aligned}
		&\quad\int_W\left(D\Phi(u),K_{21} (P_m^V v-P_n^V v)\right)_V^2\,\mu^{\Phi}(\mathrm{d}(u,v))\\
		&\leq\norm{Q_1^{\alpha}D\Phi}_{L^{\infty}(\mu_1)}^2\int_V \norm{Q_1^{-\alpha}K_{21}(P_m^V v-P_n^V v)}_V^2\,\mu_2(\mathrm{d}v)\\
		&\leq c_{\alpha}\norm{Q_1^{\alpha}D\Phi}_{L^{\infty}(\mu_1)}^2\int_U  \tilde{a} \norm{P_m^V v-P_n^V v}_V^2\,\mu_2(\mathrm{d}u)
	\end{aligned}
	\]
	for some $\tilde{a}\in (0,\infty)$ independent of $m$ and $n$.
	Invoking the theorem of dominated convergence, we see that $((D\Phi(u),K_{21}P_n^V v)_V)_{n\in\N}$ is a Cauchy-sequence in $L^2(W;\mu^{\Phi})$.
	In particular, $(D\Phi(u),K_{21}P_n^V v)_V$ converges in $L^1(W;\mu^{\Phi})$.
	In summary this implies, when setting
	\[
		g_n(u,v)\defeq 2\left( \rho_n(u,v)+h_n(v)+ \overline{h}_n(v)+(D\Phi(u),K_{21}P_n^V v)_U+\frac{1}{2} N_n(u,v)\right),
	\]
	where $(\overline{h_n})_{n\in\N}$ is defined in Assumption \nameref{ass:inf-dim-diff-bound}, 
	that
	\begin{equation}\label{eq:inf-dim-level-set-bound}
		(\operatorname{Id}-L)N_n(u,v)\leq g_n(u,v)
	\end{equation}
	for all $(u,v)\in W$ and $n\in\N$ with $n=m^K(n)$. Clearly $(g_n)_{n\in\N}$ converges in $L^1(W;\mu^{\Phi})$. Now we can finish the 		proof as in \cite[Lemma~3.3]{EG21}.
\end{proof}

\begin{remark}\label{rem:inf-dim-bounding-rho}
	To verify the first item of Assumption \nameref{ass:inf-dim-diff-bound} the following considerations are useful.
	The idea is to give a condition under which $(P_{m^K(n)}^U u,Q_1^{-1}K_{21}P_{m^K(n)}^V v)_U$ and $(P_{m^K(n)}^V v,Q_2^{-1}K_{12}			P_{m^K(n)}^U u)_U$ define Cauchy sequences in $L^1(\mu^{\Phi})$. For that let $n\geq m$, where we assume w.l.o.g.~$m=m_K(m)$ and 			$n=m_K(n)$, it holds
	\[
	\begin{aligned}
		&\quad\int_W \abs{(P_n^U u,Q_1^{-1}K_{21}P_n^V v)_U-(P_m^U u,Q_1^{-1}K_{21}P_m^V v)_U}\mu^{\Phi}(\mathrm{d}(u,v))\\
		&\leq \int_W  \abs{\sum_{i,j=m+1}^n \lambda_{1,i}^{-1}(u,e_i)_U(v,e_j)_V(d_i,K_{21}e_j)_U}e^{-\Phi(u)}\,\mu(\mathrm{d}(u,v))\\
		&\leq \sum_{i,j=m+1}^n \lambda_{1,i}^{-1}\abs{(d_i,K_{21}e_j)_U} \int_U \abs{(u,e_i)_U}\dm_1\int_V\abs{(v,e_i)_V}\,							\mu_2(\mathrm{d}v)\\
		&= \frac{{2}}{{\pi}}\sum_{i,j=m+1}^n (\sqrt{ \lambda_{1,i}})^{-1}\sqrt{\lambda_{2,i}}\abs{(d_i,K_{21}e_j)_U}\\
		&=\frac{{2}}{{\pi}}\sum_{i,j=m+1}^n \abs{(Q_1^{-\frac{1}{2}}d_i,K_{21}Q_2^{\frac{1}{2}}e_j)_U}.
		\end{aligned}
		\]
		Since the same calculation applies to $(P_n^V v,Q_2^{-1}K_{12}P_n^U u)_V$ one has to check that 
		\[
		\sum_{i,j=1}^{\infty} \abs{(Q_1^{-\frac{1}{2}}d_i,K_{21}Q_2^{\frac{1}{2}}e_j)_U}<\infty\quad\text{and}\quad \sum_{i,j=1}^{\infty} 			\abs{(Q_2^{-\frac{1}{2}}e_i,K_{12}Q_1^{\frac{1}{2}}d_j)_V}<\infty.
		\]
\end{remark}
For the rest of this section we assume that Assumptions \nameref{ass:inf-dim-diff-bound} holds. 
\begin{prop}\label{prop:inf-dim-path-prop}
	Let $\mathbf{M}=(\Omega,\mathcal{F},(\mathcal{F}_t)_{t\geq 0}, (X_t,Y_t)_{t\geq 0},(P_w)_{w\in W})$ be the right process constructed 		at the beginning of this section.
	Then $(X_t,Y_t)_{t\geq 0}$ has weakly continuous paths $P_{\mu^{\Phi}}$-a.s.
	In particular, $\mathbf{M}$ is a $\mu^{\Phi}$-invariant Hunt process.
\end{prop}
%\cite[Theorem 4.3]{EG21_Pr}
\begin{proof}
	As in \cite[Theorem 4.3]{EG21_Pr} we can construct a countable $\Q$-algebra $\mathcal{A}\subseteq \fbs{B_W}$ of $\mathcal{T}$-continuous 	functions that separates the points of $W$ and defines a core for $(L^{\Phi},D(L^{\Phi}))$.
	Since $(F_n)_{m\in\N}$ provides a $\mu^{\Phi}$-nest of $\mathcal{T}$-compact sets we can apply use \cite[Theorem 1.1]{BBR06} and 			\cite[Remark 1.2]{BBR06} to show that $(X_t,Y_t)_{t\geq 0}$ is c\`{a}dl\`{a}g w.r.t. the weak topology $\mathcal{T}$, ${P}					_{\mu^{\Phi}}$-a.e.. By \cite[Lemma 2.1.10 and Corollary 2.1.11]{Conrad2011}, the process $(X_t,Y_t)_{t\geq 0}$ has $P_{\mu^{\Phi}}$-		almost surely $\mathcal{T}$-continuous paths up to its life time. Recall that $(X_t,Y_t)_{t\geq 0}$ has $P_{\mu^{\Phi}}$-almost surely infinite life time, by the conservativity of $(L^{\Phi},D(L^{\Phi}))$.
\end{proof}
The following considerations are important to construct a solution to \eqref{eq:sde}.

\begin{lem}\label{lem:inf-dim-eval-func}
	For any $i\in \N$, define $f_i,g_i$ via
	\[
	\begin{aligned}
		W\ni (u,v)&\mapsto f_i(u,v)\defeq (u,d_i)_U\in\R\quad \text{and} \\
		W\ni (u,v)&\mapsto g_i(u,v)\defeq (v,e_i)_V\in\R.
	\end{aligned}
	\]
	Then for $i,j\in\N$, $f_i,g_i,f_if_j,g_ig_j\in D(L^{\Phi})$ and $L^{\Phi}(f_i^2),L^{\Phi}(g_ig_j)\in L^4(\mu^{\Phi})$. Moreover, for all $(u,v)\in W$
	\[
	\begin{aligned}
		L^{\Phi}f_i(u,v) &= (v,Q_2^{-1}K_{12}d_i), \\
		L^{\Phi}(f_i^2)(u,v) &= 2f_i(u,v)L^{\Phi}f_i(u,v), \\
		L^{\Phi}g_i(u,v) &= (\partial_i K_{22}(v)e_i,e_i)_V-(v,Q_2^{-1}K_{22}(v)e_i)_V-(u,Q_1^{-1}K_{21}e_i)_U\\
		&\qquad-(D\Phi(u),K_{21}e_i)_U, \\
		L^{\Phi}(g_ig_j)(u,v) &= 2(e_i,K_{22}(v)e_j)_V+g_i(u,v)L^{\Phi}g_j(u,v)+g_jL^{\Phi}g_i(u,v).
	\end{aligned}
	\]
\end{lem}
\begin{proof}
	This follows by using a sequence of cutoff functions for each $f_i$ or $g_i$, the boundedness properties of $Q_1^{\alpha}D\Phi$ by 			\nameref{ass:pert-pot} and the assumptions on the coefficient operators $K_{21}$ and $K_{22}$ described in \Cref{def:inf-dim-operators}.
\end{proof}

\begin{prop}\label{prop:inf-dim-component-solutions}
	The stochastic process $(X_t,Y_t)_{t\geq 0}$ solves the martingale problem for $(L^{\Phi},D(L^{\Phi}))$. Moreover, for any $i\in\N$, 		we have that the real-valued processes $(X_t^i)_{t\geq 0}$ and $(Y_t^i)_{t\geq 0}$ defined by $X_t^i=(X_t,d_i)_U$ and 						$Y_t^i=(Y_t,e_i)_V$ satisfy $P_{\mu^{\Phi}}$ a.s.
	\begin{equation}\label{eq:inf-dim-sde-components}
		\begin{aligned}
			X_t^i-X_0^i &= \int_0^t (Y_s,Q_2^{-1}K_{12}d_i)_V\,\mathrm{d}s \qquad\text{ and } \\
			Y_t^i-Y_0^i &= \int_0^t (\partial_i K_{22}(Y_s)e_i,e_i)_V - (Y_s,Q_2^{-1}K_{22}(Y_s)e_i)_V - (X_s,Q_1^{-1}K_{21}e_i)_U \\
			&\qquad-(D\Phi(X_s),K_{21}e_i)_U\,\mathrm{d}s
			+ M_t^{[g_i],L^{\Phi}}
		\end{aligned}
	\end{equation}
	with $(M_t^{[g_i],L^{\Phi}})_{t\geq 0}$ being a continuous $(\mathcal{F}_t)_{t\geq 0}$-martingale such that for $i,j\in\N$ we have \[[M^{[g_i],L^{\Phi}},M^{[g_j],L^{\Phi}}]_t=2\int_0^t(e_i,K_{22}(Y_s)e_j)_V\,\mathrm{d}s.\]
\end{prop}
\begin{proof}
	The statement about the martingale problem was already mentioned in \Cref{rem:inf-dim-mart-prob}.	
	 We see that $X_t^i=f_i(X_t,Y_t)$ and $Y_t^i=g_i(X_t,Y_t)$, where $f_i$ and $g_i$ are as in the definitions from \Cref{lem:inf-dim-eval-func}.
	From the statement of the Lemma itself, we see that
	\[
		M_t^{[f_i],L^{\Phi}}= X_t^i-X_0^i-\int_0^ t (Y_s,Q_2^{-1}K_{12}d_i)\,\mathrm{d}s
	\]
	and $N_t^{[f_i],L^{\Phi}}= (M_t^{[f_i],L^{\Phi}})^2$, which implies $[M^{[f_i],L^{\Phi}}]_t=0$, hence $M_t^{[f_i],L^{\Phi}}=0$.
	This proves the first line in \eqref{eq:inf-dim-sde-components}. The second line follows analogously,
	and the representation of the quadratic covariations follows by evaluating $\frac12(N_t^{[g_i+g_j],L^{\Phi}}-N_t^{[g_i],L^{\Phi}}-			N_t^{[g_j],L^{\Phi}})$.
\end{proof}
Finally, we want to prove that the process is also a stochastically and analytically weak solution of the infinite-dimensional stochastic differential equation \eqref{eq:sde}. For this, we need to construct a suitable cylindrical Brownian motion on $V$, such that we can express the process described by $M^{V}_t\defeq \sum_{i\in\N} M_t^{[g_i],L^{\Phi}}e_i$ as a stochastic integral of $\sqrt{K_{22}}$.
In the following we set for $n\in\N$ such that $n=m^K(n)$
\[M^{(n)}_t\defeq (M_t^{[g_1],L^{\Phi}},\dots, M_t^{[g_n],L^{\Phi}})\quad\text{and}\quad\Sigma^{(n)}_t\defeq \left((K^{-\frac12}_{22}(P_n^V Y_t)e_i,e_j)_V\right)_{i,j=1}^n.\]By L\'{e}vy's characterization of Brownian motion, we get:
\begin{lem}
	For each $n\in\N$ such that $n=m^K(n)$ 
	\[
		B^{(n)}_t\defeq \int_0^t \Sigma^{(n)}_s\,\mathrm{d}M^{(n)}_s
	\]
	is an $n$-dimensional Brownian motion. Moreover, let $\beta^{(k)}$ be the $k$-th component of $B^{(m^K(k))}$. Then $(\beta^{(k)})_{k\in\N}$ is an independent sequence of one-dimensional Brownian motions.
\end{lem}
\begin{proof}
	Let $k\in\N$.
	By L\'{e}vy's characterization, it follows that
	\[
		B^{(m^K(k))}_t\defeq \int_0^t \Sigma^{(m^K(k))}_s\,\mathrm{d}M^{(m^K(k))}_s
	\]
	is an $m^K(k)$-dimensional Brownian motion.
	As $\beta^{(k)}_t$ is the $k$-th component of $B^{(m^K(k))}_t$ we know that
	$\{\beta^{(1)},\dots,\beta^{(m^K(k))}\}$ is independent for any $k\in\N$. Due to the block diagonal structure of  $\Sigma^{(m^K(k))}$ we 			know that $\beta^{(j)}$ does not depend on the $k\in\N$ such that $j\leq m^K(k)$.
\end{proof}
Now we fix some $T\in(0,\infty)$ and define the process $(B_t)_{t\in [0,T]}$ on $V$ via
\[
	B_t\defeq \sum_{k=1}^\infty \beta^{(k)}_t e_k, \qquad t\in [0,T].
\]
This is a cylindrical Brownian motion on $V$ as defined in \cite[Proposition~2.5.2]{PR07},
as can be seen by choosing $\mathcal{J}:V\to V$, $\mathcal{J}\defeq Q_2^{\frac12}$, since then
\[
	B_t^{Q_2} \defeq \sum_{k=1}^\infty \beta^{(k)}_t \mathcal{J}e_k
\]
defines a $Q_2$-Wiener process on $V$.
In the following we set $V_0\defeq Q_2^{\frac12}V$ and equip it with the inner product
\[
	(a,b)_{V_0}\defeq (Q_2^{-\frac12}a,Q_2^{-\frac12}b)_V
	\qquad\text{ for all }a,b\in V_0,
\]
which makes $V_0$ a separable Hilbert space with orthonormal basis $(Q_2^{\frac12}e_i)_{i\in\N}$.
Define $\mathcal{L}_2^0\defeq \lophs{V_0;V}$ as the Hilbert space of Hilbert-Schmidt operators from $V_0$ to $V$.
Further let
\[
	\mathcal{A}_T\defeq \sigma\{Y:[0,T]\times\Omega\to \R\mid Y\text{ is left-continuous and }(\mathcal{F}_t)_{t\in[0,T]}\text{-adapted} \}.
\]
\begin{cond}{K5}\label{ass:weak_sol}
	There is a a non-negative function $k_{22}$ in $L^1(V;\mu_2)$ such that for all $v\in V$
	\[
	\tr[K_{22}(v)]\leq k_{22}(v).
	\]
\end{cond}
\begin{lem}
	The stochastic process $\sqrt{K_{22}(Y_t)}\mathcal{J}^{-1}$, $t\in[0,T]$, is $\mathcal{L}_2^0$-valued and predictable, i.e.~$\mathcal{A}_T$-$\bs{\mathcal{L}_2^0}$-measurable.
\end{lem}
\begin{proof}
	For each $v\in V$, we have by \nameref{ass:weak_sol}
	\[
		\sum_{i\in\N} (\sqrt{K_{22}(v)}\mathcal{J}^{-1}Q_2^{\frac12}e_i,\sqrt{K_{22}(v)}\mathcal{J}^{-1}Q_2^{\frac12}e_i)_V\leq
		k_{22}(v).
	\]
	This implies that $\sqrt{K_{22}(v)}\mathcal{J}^{-1}\in\mathcal{L}_2^0$ for any $v\in V$. Moreover, for all $t\in (0,\infty)$ we have
	\[
	\begin{aligned}
		\int_{\Omega}\int_0^t \|\sqrt{K_{22}}(Y_s(\omega))\mathcal{J}^{-1}\|_{\mathcal{L}_2^0}\,\mathrm{d}s\,P_{\mu^{\Phi}}(\mathrm{d}				\omega)&\leq \int_0^t \int_{\Omega} k_{22}(Y_s(\omega))\,P_{\mu^{\Phi}}(\mathrm{d}\omega)\,\mathrm{d}s\\&=\int_0^t\int_W 					(T_sk_{22})(v) \,\mu^{\Phi}(\mathrm{d}(u,v))\,\mathrm{d}s\\&\leq t \norm{k_{22}}_{L^1(V,\mu_2)}.
	\end{aligned}
	\]
 	Hence for all $t\in (0,\infty)$: \[\int_0^t\|\sqrt{K_{22}}(Y_s(\omega))\mathcal{J}^{-1}\|_{\mathcal{L}_2^0}\,\mathrm{d}s<\infty\quad 		\mathbb{P}_{\mu^{\Phi}} a.s..\] Moreover, 
	\[
		A_i\defeq(\sqrt{K_{22}}(Y_t)\mathcal{J}^{-1}Q_2^{\frac12}e_i,\sqrt{K_{22}}(Y_t)\mathcal{J}^{-1}Q_2^{\frac12}e_i)_{V}
		=(K_{22}(P_{m^K(i)}^VY_t)e_i,e_i)_V
	\]
	is continuous and $(\mathcal{F}_t)_{t\in[0,T]}$-adapted for any $i\in\N$.
	Fix some $\eps>0$ and set
	\[
		B\defeq \{(t,\omega)\in[0,T]\times\Omega\mid \|\sqrt{K_{22}}(Y_t(\omega))\mathcal{J}^{-1}\|_{\mathcal{L}_2^0}\leq\eps \}
	\]
	as well as $B_k\defeq \{\sum_{i=1}^k A_i \leq \eps \}\in \mathcal{A}_T$ for each $k\in\N$.
	Then $B=\bigcap_{k\in\N} B_k\in\mathcal{A}_T$ as well. It is easily seen that similarly, all pre-images of closed $\eps$-balls in $			\mathcal{L}_2^0$ under $\sqrt{K_{22}(Y_\cdot)}\mathcal{J}^{-1}$ are in $\mathcal{A}_T$, so that the process is indeed predictable, 			since $\mathcal{L}_2^0$ is separable.
\end{proof}
By \cite[Section~2.3]{PR07}, the previous results imply that $\sqrt{K_{22}(Y_t)}\mathcal{J}^{-1}$ is integrable with respect to the $Q_2$-Wiener process $(B_t^{Q_2})_{t\in[0,T]}$, which shows that $\sqrt{K_{22}(Y_t)}$ is integrable with respect to $(B_t)_{t\in[0,T]}$ with
\[
	I(\sqrt{K_{22}(Y_t)})\defeq \int_0^t \sqrt{K_{22}(Y_s)}\,\mathrm{d}B_s\defeq \int_0^t \sqrt{K_{22}(Y_s)}\mathcal{J}^{-1}\,\mathrm{d} B_s^{Q_2}
\]
for all $t\in[0,T]$.
By applying \cite[Lemma~2.4.1]{PR07} for the operators $(\cdot,e_i)_V:V\to\R$, we see that
\[
	\sqrt{2}(I(\sqrt{K_{22}(Y_t)}),e_i)_V
	= 	\sqrt{2}\int_0^t (\sqrt{K_{22}(Y_s)}\mathcal{J}^{-1} \cdot, e_i)_V \,\mathrm{d} B_s^{Q_2}
	= 	M_t^{[g_i],L^{\Phi}}
\]
due to the block invariance properties of $K_{22}$. We obtain
\begin{equation}\label{eq:inf-dim-martingale-transform}
	M_t^{V}\defeq \sum_{i\in\N} M_t^{[g_i],L^{\Phi}} e_i = \int_0^t \sqrt{K_{22}(Y_s)}\,\mathrm{d}B_s.
\end{equation}
As a result, we get the following:
\begin{thm}\label{thm:sto_ana_weak_sol}
	Let \nameref{ass:inf-dim-elliptic}--\nameref{ass:weak_sol}, then $\mathbf{M}=(\Omega,\mathcal{F},(\mathcal{F}_t)_{t\geq 0}, (X_t,Y_t)_{t\geq 0},				(P_w)_{w\in W})$ provides a stochastically and analytically weak solution of \eqref{eq:sde} in the sense that there is a cylindrical 				Brownian motion $(B_t)_{t\geq 0}$ on $V$ such that $P_{\mu^{\Phi}}$ a.s.~for $\nu_1\in D(Q_2^{-1}K_{12})$ and $\nu_2\in D(Q_1^{-1}			K_{21})\cap D(Q_2^{-1}K_{22}(v))$ for all $v\in V$ we have
	\begin{equation}\label{eq:inf-dim-sde}
		\begin{aligned}
			(X_t-X_0,\nu_1)_U &= \int_0^t (Y_s,Q_2^{-1}K_{12}\nu_1)_U\,\mathrm{d}s\quad\text{and} \\
			(Y_t-Y_0,\nu_2)_V &=\int_0^t \left(\sum_{i=1}^\infty \partial_iK_{22}(Y_s)e_i,\nu_2\right)_V
			- (Y_s,Q_2^{-1}K_{22}(Y_s)\nu_2)_V\\&\qquad-(X_s,Q_1^{-1}K_{21}\nu_2)_U-(K_{12}D\Phi(X_s),\nu_2)_V\,\mathrm{d}s
			\\&\qquad+\left(\int_0^t \sqrt{2K_{22}(Y_s)}\,\mathrm{d}B_s,\nu_2\right)_V.
		\end{aligned}
	\end{equation}
\end{thm}
\begin{proof}
	Given $\nu_1\in D(Q_2^{-1}K_{12})$ and $\nu_2\in D(Q_1^{-1}K_{21})\cap D(Q_2^{-1}K_{22}(v))$ for all $v\in V$. For $n\in\N$ with 			$n=m^K(n)$ it holds by \Cref{prop:inf-dim-component-solutions}
		\[
		\begin{aligned}
			(X_t-X_0,P_n^U\nu_1)_U&=\int_0^t (Y_s,P_n^VQ_2^{-1}K_{12}\nu_1)_V\,\mathrm{d}s\\
			(Y_t-Y_0,P_n^V\nu_2)_V &= \int_0^t \left(\sum_{i=1}^\infty\partial_i K_{22}(Y_s)e_i,P_n^V\nu_2\right)_V- (Y_s,P_n^VQ_2^{-1}K_{22}(Y_s)					\nu_2)_V\\
			& \qquad- (X_s,P_n^UQ_1^{-1}K_{21}\nu_2)_U -(D\Phi(X_s),K_{21}P_n^V\nu_2)_U\,\mathrm{d}s\\
			&\qquad +( M_t^{V},P_n^V\nu_2)_V.
		\end{aligned}
		\]
	Note that $P_n^U\nu_1$, $P_n^VQ_2^{-1}K_{12}\nu_1$, $P_n^V\nu_2$, $P_n^VQ_2^{-1}K_{22}(Y_s)\nu_2$ and  $P_n^UQ_1^{-1}K_{21}\nu_2$ converge pointwisely to $\nu_1$, $Q_2^{-1}K_{12}\nu_1$, $\nu_2$, $Q_2^{-1}K_{22}(Y_s)\nu_2$ and  $Q_1^{-1}			K_{21}\nu_2$. Moreover, using \nameref{ass:inf-dim-diff-bound} and \nameref{ass:pert-pot} we can estimate for some $a,b,\tilde{a}\in\N$, compare \Cref{prop:inf-dim-assoc-process},
	\[
	\begin{aligned}
		&\abs{(\sum_{i=1}^\infty\partial_i K_{22}(Y_s)e_i,P_n^V\nu_2)_V}\leq a(1+ \norm{Y_s}^b_V) \norm{\nu_2}_V\quad\text{and}\\
		&\abs{(D\Phi(X_s),K_{21}P_n^V\nu_2)_U}\leq \sqrt{c_{\alpha}}\norm{Q_1^{\alpha}D\Phi}_{L^{\infty}(\mu_1)}\tilde{a}\norm{\nu_2}						_V.
	\end{aligned}
	\]
	Using the fact that weakly continuous functions are norm bounded we are able to apply the theorem of dominated convergence for 			$n\rightarrow  \infty$ with $n=m^K(n)$ to obtain the statement.
\end{proof}

\section{Abstract hypocoercivity framework}\label{sec:abstract-hypoc}
We include here the slight reformulation from \cite[Section 2.2]{Bertram2023} of the abstract Hilbert space hypocoercivity method described in \cite{GS14}.

Let $H$ be a separable Hilbert space with inner product $\spdx$ and induced norm $\normx$,
which has an orthogonal decomposition $H=H_1\oplus H_2$ with corresponding orthogonal projections $P:H\to H_1$, $(\operatorname{Id}-P):H\to H_2$.
Let further $(L,D(L))$ be a densely defined linear operator that generates a strongly continuous contraction semigroup $\sccs$ on $H$.
We assume the following structure on $L$:

\begin{cond}{D1}\label{ass:data-op-decomposition}
	$L=S-A$ on $\acore$, where $S$ is symmetric, $A$ is antisymmetric, and $\acore\subseteq D(L)$ is a core for $(L,D(L))$.
\end{cond}
Then both $(S,\acore)$ and $(A,\acore)$ are closable, and we denote their closures by $(S,D(S))$ and $(A,D(A))$, respectively.
These two operators are linked to the decomposition of $H$ in the following way:
\begin{cond}{D2}\label{ass:data-symmetric-compliance}
	$H_1\subseteq D(S)$ and $S\equiv 0$ on $H_1$.
\end{cond}
\begin{cond}{D3}\label{ass:data-anti-compliance}
	$P(\acore)\subseteq D(A)$, $AP(\acore)\subseteq D((AP)^*)$ and $PAP\equiv 0$ on $\acore$.
	Here $(AP)^*$ is the adjoint of the densely defined closed operator $(AP,D(AP))$ with
	\[ 
		D(AP)=\{ x\in H\mid Px\in D(A)\}.
	\]
\end{cond}

\begin{defn}\label{def:paap-operator}
	We define the operator $(G,D(G))$ by
	\[
		G\defeq -(AP)^*AP,
			\qquad
		D(G)\defeq \{ x\in D(AP)\mid APx\in D((AP)^*) \}.
	\]
\end{defn}
\begin{remark}
	Due to von Neumann's theorem (\cite[Theorem~5.1.9]{P89}), $(G,D(G))$ is self-adjoint
	and $\operatorname{Id}-G:D(G)\to H$ is bijective with bounded inverse. Since $G$ is dissipative, it generates an operator semigroup on $H$.
	
	Note that due to \nameref{ass:data-anti-compliance},
	we have $\acore\subseteq D(G)$. If additionally, $AP(\acore)\subseteq D(A)$, then $G=PA^2P$ on $\acore$.
\end{remark}

This allows us to define the following operator, which is bounded with operator norm less than 1 due again to \cite[Theorem~5.1.9]{P89}:
\begin{defn}\label{def:b-operator}
	Define the operator $(B,D(B))$ as
	\[
		B\defeq (\operatorname{Id}-G)^{-1}(AP)^*,
			\qquad
		D(B)\defeq D((AP)^*).
	\] 
	Due to boundedness, it extends uniquely to a bounded operator $B:H\to H$.
\end{defn}

\begin{cond}{H1}\label{ass:aux-bound}
	\emph{Boundedness of auxiliary operators:}
	The operators $(BS,\acore)$ and $(BA(I-P),\acore)$ are bounded
	and there exist constants $c_1,c_2<\infty$ such that
	\[
		\|BSx\|\leq c_1\|(\operatorname{Id}-P)x\| \quad\text{ and }\quad
		\|BA(\operatorname{Id}-P)x\| \leq c_2 \|(\operatorname{Id}-P)x\|
	\]
	hold for all $x\in \acore$. 
\end{cond}
\begin{cond}{H2}\label{ass:micro-coerc}
	\emph{Microscopic coercivity:} There exists some $\Lambda_m>0$ such that
	\[
		-(Sx,x) \geq \Lambda_m\|(\operatorname{Id}-P)x \|^2\qquad \text{ for all }x\in \acore.
	\]
\end{cond}
\begin{cond}{H3}\label{ass:macro-coerc}
	\emph{Macroscopic coercivity:}
	There is some $\Lambda_M>0$ such that
	\begin{equation}\label{eq:macro-coerc-ineq}
		\|APx\|^2 \geq \Lambda_M \|Px\|^2 \qquad\text{ for all }x\in D(G).
	\end{equation}
\end{cond}
\begin{remark}\label{remark:core_if_ess_self}
	If $(G,\acore)$ is already essentially self-adjoint, then \nameref{ass:macro-coerc} is satisfied if \eqref{eq:macro-coerc-ineq} holds 		for all $x\in \acore$.
	For the proof compare \cite[Corollary 2.13]{GS14}.
\end{remark}

The following hypocoercivity theorem is the central statement in this section, compare also \cite[Theorem~2.2]{GS16}.
\begin{thm}\label{thm:abstract-hypoc}
	Assume that \nameref{ass:data-op-decomposition}--\nameref{ass:data-anti-compliance} and 
	\nameref{ass:aux-bound}--\nameref{ass:macro-coerc} hold. Then there exist strictly positive
	constants $\kappa_1,\kappa_2<\infty$ which are explicitly computable
	in terms of $\Lambda_m$, $\Lambda_M$, $c_1$ and $c_2$ such that
	for all $x\in H$ we have
	\[
		\|T_tx\| \leq \kappa_1\mathrm{e}^{-\kappa_2t} \|x\|\quad\text{ for all }t\geq0.
	\]
\end{thm}
In order to verify H1, we state \cite[Lemma 3.1]{BG21-wp}:
\begin{lem}\label{lem:auxiliary-bound}
	Let $\acore$ be a core for $(G,D(G))$. Let $(T,D(T))$ be a linear operator with $\acore\subseteq D(T)$ and assume $AP(\acore)\subseteq D(T^*)$.
	Then
	\[
		(\operatorname{Id}-G)(\acore) \subseteq D((BT)^*)
		\quad\text{ with }\quad
		(BT)^*(\operatorname{Id}-G)x = T^*APx,\quad x\in \acore.
	\]
	If there exists some $C<\infty$ such that
	\begin{equation}\label{eq:adjoint_bound}
		\|(BT)^* y\| \leq C\|y\|\qquad\text{ for all } y=(\operatorname{Id}-G)x,\quad x\in \acore,
	\end{equation}
	then $(BT,D(T))$ is bounded and its closure $(\overline{BT})$ is a continuous operator on $H$ with $\| \overline{BT}\|=\| (BT)^*\|\leq C$.
	In particular, if $(S,D(S))$ and $(A,D(A))$ satisfy these assumptions with constant $C_S$ and $C_A$, respectively,
	then \nameref{ass:aux-bound} is satisfied with $c_1=C_S$ and $c_2=C_A$.
\end{lem}

\section{Hypocoercivity}\label{sec:hypo_appl}

Throughout this section, we assume \nameref{ass:inf-dim-elliptic}, \nameref{ass:inf-dim-growth} and \nameref{ass:Phi_general} for some $\alpha\in [0,\infty)$ unless specifically stated otherwise. If we additionally assume \nameref{ass:pert-pot}, \Cref{thm:ess_diss_Lphi} applies.

We now restrict the setting to the Hilbert space \[H\defeq \left\{f\in L^2(W;\mu^{\Phi}): \mu^{\Phi}(f)=0  \right\}\] and operator domain $\acore\defeq\fbs{B_W}\cap H$.

\begin{prop}\label{prop:inf-dim-sub-markov}
	Assume \nameref{ass:pert-pot}, then the operator $(L^{\Phi},\acore)$ is essentially m-dissipative on $H$ and its closure $(L^{\Phi}_0,D(L^{\Phi}_0))$ generates a strongly continuous contraction semigroup $\sccsz$ on $H$.
\end{prop}
\begin{proof}
	By $\mu^{\Phi}$-invariance of $L^{\Phi}$ and $\sccs$, $(L^{\Phi},\acore)$ and the restriction $\sccsz$ of $\sccs$ to $H$ are well-			defined as operators on $H$.
	Dissipativity of $(L^{\Phi},\acore)$ is inherited from $(L^{\Phi},\fbs{B_W})$, and the dense range condition can be verified as 			follows: Let $f\in H$, then there is a sequence $\s{f}$ in ${L^2(\mu^{\Phi})}$ such that $(\operatorname{Id}-L^{\Phi})f_n\to f$ in ${L^2(\mu^{\Phi})}		$. In particular, $\mu^{\Phi}(f_n)\to\mu^{\Phi}(f)=0$, so by setting $g_n\defeq f_n-\mu^{\Phi}(f_n)$, it follows that $(\operatorname{Id}-					L^{\Phi})g_n\to f$ since $L^{\Phi}$ acts trivially on constants.
	Since $g_n\in \acore$ for all $n\in\N$, it follows that $(L^{\Phi},\acore)$ is essentially m-dissipative, and its closure $(L^{\Phi}		_0,D(L^{\Phi}_0))$ is the generator of $\sccsz$.
\end{proof}

\begin{defn}
	Let $H=H_1\oplus H_2$, where $H_1$ is provided by the orthogonal projection
	\[
		P:H\to H_1, \quad
		f\mapsto Pf\defeq \int_{V}f(\cdot,v)\,\mu_2(\mathrm{d}v).
	\]
	For any $f\in \acore$, we can interpret $Pf$ as an element of $\fbs{B_U}$, in which case we denote it by $f_P$.
	Further let $(S_0,D(S_0))$ and $(A^{\Phi}_0,D(A^{\Phi}_0))$ be the closures in $H$ of $(S,\acore)$ and $(A^{\Phi},\acore)$, respectively.
\end{defn}

As in \cite{EG21}, we also define the following operators:
\begin{defn}
	The operators $(C,D(C))$ and $(Q_1^{-1}C,D(Q_1^{-1}C))$ on $U$ are defined by
	\[
	\begin{aligned}
		C&\defeq K_{21}Q_2^{-1}K_{12}, &D(C)&\defeq \{u\in U\mid K_{12}u\in D(Q_2^{-1}) \}\\
		Q_1^{-1}C&\defeq Q_1^{-1}K_{21}Q_2^{-1}K_{12}, &D(Q_1^{-1}C)&\defeq \{u\in D(C)\mid Cu\in D(Q_1^{-1}) \},
	\end{aligned}
	\]
	respectively. 
	For all $f\in\fbs{B_U}$, define $Nf(u)\defeq \tr[CD^2f(u)]-(u,Q_1^{-1}CDf(u))_U$ for all $u\in U$.
\end{defn}

In order to gain useful properties of $N$, we need the following:
\begin{cond}{K6}\label{ass:qc-negative-type}
	The operator $K_{21}K_{12}=K_{12}^*K_{12}$ is positive-definite on $U$.
\end{cond}
For the rest of this section we assume that \nameref{ass:qc-negative-type} holds. We continue with the collection of a few properties of the newly defined operators:
\begin{prop}\label{prop:ess_self_N}
	Let $(C,D(C))$ and $(Q_1^{-1}C,D(Q_1^{-1}C))$ be defined as above. Then
	\begin{enumerate}[(i)]
		\item $(C,D(C))$ is symmetric and positive definite on $U$.
		\item For all $n\in\N$ there is $m_k$ with $n\leq m_k$ such that $(C,D(C))$ maps $U_n$ into $U_{m_k}$.
		\item The operator $(N,\fbs{B_U})$ is essentially self-adjoint on $L^2(U;\mu_1)$ with closure $(N,D(N))$ and resolvent $						(R_{\lambda}^{N})_{\lambda>0}$. Fix $\lambda\in (0,\infty)$. It holds 
			\begin{equation}\label{eq:regularity-estimate-one_no_pot}
				\int_U(CDR_{\lambda}^Nf,DR_{\lambda}^Nf)_U\,\mathrm{d}\mu_1\leq \frac{1}{4\lambda}\int_U f^2\,\mathrm{d}\mu_1
			\end{equation}
			for all $f\in L^2(U;\mu_1)$. Moreover for all $f\in \fbs{B_U}$ we have
			\begin{equation}\label{eq:regularity-estimate-two_no_pot}
				\int_U \tr[(CD^2R_{\lambda}^Nf)^2]+\norm{Q_1^{-\frac{1}{2}}CDR_{\lambda}^Nf}^2_U\,\mathrm{d}\mu_1=\int_U(NR_{\lambda}						^Nf)^2\,\mathrm{d}\mu_1\leq 4 \int_U f^2\,\mathrm{d}\mu_1.
			\end{equation}
			In particular $D(N)\subseteq W^{1,2}_{Q_1^{-\frac{1}{2}}C}(U;\mu_1) \subseteq W^{1,2}_{C}(U;\mu_1)$ with 
			\begin{equation}\label{eq:regularity-estimate-three_no_pot}
				\int_U \norm{Q_1^{-\frac{1}{2}}CDR_{\lambda}^Nf}^2_U\,\mathrm{d}\mu_1\leq 4 \int_U f^2\,\mathrm{d}\mu_1
			\end{equation}
			for all $f\in L^2(U;\mu_1)$.
	\end{enumerate}
\end{prop}
\begin{proof}
	\begin{enumerate}[(i)]
		\item Let $u_1,u_2\in D(C)$. Then
			\[
				(Cu_1,u_2)_U= (Q_2^{-1}K_{12}u_1,K_{12}u_2)_V
				= (K_{12}u_1,Q_2^{-1}K_{12}u_2)_V
				= (u_1,Cu_2)
			\]
			due to symmetry of $Q_2^{-1}$ and definition of $K_{21}$. Since $Q_2^{-1}$ is positive-definite, positive definiteness of 			$C$ follows immediately by Assumption \nameref{ass:qc-negative-type}.
		\item Since there is some $m_k$ with $n\leq m_k$ such that $K_{12}$ maps $U_n$ to $V_{m_k}$, $K_{21}$ maps $V_n$ to $U_{m_k}$, 
			and $Q_2^{-1}$ leaves $V_{m_k}$ invariant, the claim follows directly.
		\item This can be proven analogously to \Cref{thm:inf-dim-ess-m-diss}, since the matrices in $\R^{n\times n}$ induced by
			$C$ are constant with positive eigenvalues, which allows usage of \cite[Thm.~4.5]{BG21} by point (i).
			Even though $N$ might be unbounded, it is well-defined on $\fbs{B_U}$ within $L^2(U;\mu_1)$, and dissipativity is implied by 				the integration by parts formula from \Cref{lem:inf-dim-ibp-fbs}. The regularity estimates can be shown as in 								\Cref{lem:reg_K_22} together with \cite[Theorem 3.3]{EG21}
	\end{enumerate}
\end{proof}
In analogy to \Cref{Essential m-dissipativity} we want to consider non trivial Potentials $\Phi$ and show essential self-adjointness of $(N^{\Phi},\fbs{B_U})$ defined by
\begin{equation}\label{eq:N_Phi}
	\fbs{B_U}\ni f\mapsto N^{\Phi}f\defeq\mathrm{tr}[CD^2f]-(u,Q_1^{-1}CD f)_U-(D\Phi,CDf)_U\in L^2(U;\mu_1^{\Phi}).
\end{equation}
In order to do that the next assumption, dealing with the regularity of $\Phi$, is important.

\begin{cond}{$\Phi1$}\label{ass:ess-N_Phi}Assume either
	\begin{enumerate}[(i)]
	\item $\Phi\in W_{Q_1^{\frac{1}{2}}}^{1,2}(U;\mu_1)$ with $\norm{Q_1^{\frac{1}{2}}D\Phi}_{L^{\infty}(\mu_1)}< \frac{1}{2}$ or
	\item $C\in \mathcal{L}(U)$ and $\Phi\in W_{C^{\frac{1}{2}}}^{1,2}(U;\mu_1)$ with  $\norm{C^{\frac{1}{2}}D\Phi}_{L^{\infty}(\mu_1)}<\infty$.		
	\end{enumerate}
\end{cond}
\begin{thm}\label{essdissN}
 	The operator $(N^{\Phi},\fbs{B_U})$ defined by \eqref{eq:N_Phi}, fulfills 
	\begin{equation}\label{eq:intbp_Phi}
		(N^{\Phi}f,g)_{L^2(U;\mu_1^{\Phi})}=\int_U -(CDf,Dg)_U\dm_1^{\Phi}
	\end{equation}for all $f,g\in \fbs{B_U}$. Assume that \nameref{ass:ess-N_Phi} is satisfied, then it holds $D(N)\subseteq D(N^{\Phi})$ 		with
	\[
		N^{\Phi}f=Nf-(D\Phi,CDf)_U\in L^2(U;\mu_1^{\Phi})
	\]
	for all $f\in D(N)$.
	Moreover, $(N^{\Phi},\fbs{B_U})$ is essentially self-adjoint in $L^2(U;\mu_1^{\Phi})$. The resolvent in $\lambda\in (0,\infty)$ of the 	closure $(N^{\Phi},D(N^{\Phi}))$ is denoted by $R_{\lambda}^{N^{\Phi}}$.
\end{thm}
	\begin{proof}
	The first item of the statement follows by the integration by parts formula together with the invariance properties of the involved 		operators.
	Assume that Item $(i)$ of \nameref{ass:ess-N_Phi} is valid.
	For $f\in L^2(U;\mu_1)$ set
	\[
		Tf=-(D\Phi,CDR_{1}^Nf)_U.
	\]
	Since $D(N)\subseteq W^{1,2}_{Q_1^{-\frac{1}{2}}C}(U;\mu_1) \subseteq W^{1,2}_{C}(U;\mu_1) $ the definition above is reasonable. Using the Cauchy-	Schwarz inequality, Inequality \eqref{eq:regularity-estimate-three_no_pot} and the assumption on $\Phi$ we observe
	\[
		\norm{Tf}_{L^2(\mu_1)}^2\leq \norm{Q_1^{\frac{1}{2}}D\Phi}^2_{L^{\infty}(\mu_1)}\int_U \norm{Q_1^{-\frac{1}{2}}CDR_{1}^Nf}^2_U\dm_1< \norm{f}_{L^2(\mu_1)}^2.
	\]
	Therefore the linear operator $T:L^2(U;\mu_1)\rightarrow L^2(U;\mu_1)$ is well-defined with operator norm less than one. Hence by the 		Neumann-Series theorem we obtain that $(\operatorname{Id}-T)^{-1}$ exists in $\mathcal{L}(L^2(U;\mu_1))$. In particular for a given $g\in L^2(U;\mu_1)$ we find $f\in L^2(U;\mu_1)$ with $f-Tf=g$ in $L^2(U;\mu_1)$. Since $(N,D(N))$ is m-dissipative, there is $h\in D(N)$ with $(\operatorname{Id}-N)h=f$. This yields
	\[
		(\operatorname{Id}-N)h+(D\Phi,CDh)_U=f+(D\Phi,CDR_{1}^Nf)_U=f-Tf=g.
	\]
	If we can show that $D(N)\subseteq D(N^{\Phi})$ with $N^{\Phi}f=Nf-(D\Phi,CDf)_U$ for all $f\in D(N)$ the proof is finish by the Lumer-		Philipps theorem. Indeed, this implies 
	\[
		\fbs{B_U}\subseteq L^2(U;\mu_1)\subseteq(\operatorname{Id}-N^{\Phi})(D(N)) \subseteq(\operatorname{Id}-N^{\Phi})(D(N^{\Phi})),
	\]i.e.~the dense range condition, as $\fbs{B_U}$ is dense in $L^2(U;\mu_1^{\Phi})$. So let $f\in D(N)$ be given. There is a sequence $		(f_n)_{n\in\N}\subseteq \fbs{B_U}$ s.t.~$f_n\rightarrow f$ and $Nf_n\rightarrow Nf$ in $L^2(U;\mu_1)$. As $\abs{e^{-\Phi}}\leq 1$, it is easy to see that $f_n\rightarrow f$ in $L^2(U;\mu_1^{\Phi})$. 
	In view of the assumptions on $Q_1^{\frac{1}{2}}D\Phi$ and the Inequality \eqref{eq:regularity-estimate-two_no_pot} we can estimate
	\[
	\begin{aligned}
		&\quad\norm{Nf-(D\Phi,CDf)_U-N^{\Phi}f_n}_{L^2(\mu_1^{\Phi})}\\
		&\leq \norm{N(f-f_n)}_{L^2(\mu_1)}+\norm{Q_1^{\frac{1}{2}}D\Phi}_{L^{\infty}(\mu_1)} \norm{Q_1^{-\frac{1}{2}}CD(f-f_n)}						_{L^2(\mu_1)}\\
		 &\leq \norm{N(f-f_n)}_{L^2(\mu_1)}+\frac{1}{2}\norm{(f-f_n)+N(f-f_n)}_{L^2(\mu_1)}.
	\end{aligned}
	\]
	Therefore, $N^{\Phi}f_n\rightarrow Nf-(D\Phi,CDf)_U$ in $L^2(U;\mu_1^{\Phi})$. Since $(N^{\Phi},D(N^{\Phi}))$ is closed by 						construction we obtain $D(N)\subseteq D(N^{\Phi})$ with $N^{\Phi}f=Nf-(D\Phi,CDf)_U$ for all $f\in D(N)$ as desired.
	If Item $(ii)$ from \nameref{ass:ess-N_Phi} is valid we proceed as in Theorem \ref{thm:ess_diss_Lphi} and use Inequality 					\eqref{eq:regularity-estimate-one_no_pot}.
\end{proof}

\begin{remark}
	Note that if $C\in \mathcal{L}(U)$ there are different assumption yielding essential self-adjointness of $(N^{\Phi},\fbs{B_U})$ in 			$L^2(U;\mu_1^{\Phi})$, compare \cite{BigF22}. In particular boundedness of $C^{\frac{1}{2}}D\Phi$ is not required. But as we need to 		assume \nameref{ass:pert-pot} to obtain essential m-dissipativity of $(L^{\Phi},\fbs{B_H})$ in $L^2(U;\mu^{\Phi})$ anyway, we 				formulated Assumption \nameref{ass:ess-N_Phi} as above.
\end{remark}
The next assumption is made to generalize the regularity estimates from \Cref{prop:ess_self_N}. The proof of the second order regularity estimate is similar to \cite[Theorem~2]{EG21}, where only convex potentials were considered. 
\begin{cond}{$\Phi2$}\label{ass:L_Phi-convex-reg-est}
	\begin{enumerate}[(i)]
		\item $\Phi=\Phi_1+\Phi_2$. $\Phi_1$ is in $W_{Q_1^{\alpha}}^{1,2}(U;\mu_1)$ for some $\alpha\in [0,\infty)$ and fulfills 						Assumption \nameref{hypoapproximationPotential}. $\Phi_2$ is bounded and two times continuously Fr\'{e}chet-differentiable with 				bounded first order derivative and second order derivative in $L^1(\mu_1^{\Phi})$.
		\item There is a constant $c_{\Phi_2}\in [0,\infty)$ such that for all $n\in\N$ 
			\[
				(D^2\Phi_2(\tilde{u})Cu,Cu)_U\geq -c_{\Phi_2}(Cu,u)_U\quad\text{for\;all}\quad \tilde{u}\in U\quad\text{and}\quad u\in 						U_n. 
			\]
	\end{enumerate}
\end{cond}
Note that Item $(ii)$ above is satisfied if $D^2\Phi_2$ is bounded and $C\in \mathcal{L}(U)$.
\begin{lem}\label{lem:inf-dim-regularity-estimates}
	For all $f\in  L^2(U;\mu^{\Phi}_1)$ and $\lambda\in (0,\infty)$ it holds
	\[
	\begin{aligned}\label{eq:regularity-estimate-one}
		\int_U(CDR_{\lambda}^{N^{\Phi}}f,DR_{\lambda}^{N^{\Phi}}f)_U\,\mathrm{d}\mu_1^{\Phi}
		\leq \frac{1}{4\lambda}\int_U f^2\,\mathrm{d}\mu_1^{\Phi}.
	\end{aligned}
	\]
	Assume Assumption \nameref{ass:L_Phi-convex-reg-est} holds true. Then for all $f\in \fbs{B_U}$ and $g=f-N^{\Phi}f$
	\[
	\begin{aligned}\label{eq:regularity-estimate-two}
		\int_U \tr[(CD^2f)^2]+\norm{Q_1^{-\frac{1}{2}}CDf}^2_U\,\mathrm{d}\mu_1^{\Phi}
		&\leq \left(4+\frac{c_{\Phi_2}}{4}\right) \int_U g^2\,\mathrm{d}\mu_1^{\Phi}.
	\end{aligned}	
	\]
\end{lem}
\begin{proof}
	The first inequality follows by \eqref{eq:intbp_Phi} and a reasoning similar to \Cref{lem:reg_K_22}. For the second we approximate $		\Phi_1$ with $(\Phi_{1,(n,m)})_{(n,m)\in \N^2}$, provided by \nameref{hypoapproximationPotential} and set $\Phi_{(n,m)}\defeq \Phi_{1,		(n,m)}+\Phi_2$. As in \cite[Theorem~2]{EG21} we obtain for $f\in \fbs{B_U}$ and $g_{(n,m)}=f-N^{\Phi_{(n,m)}}f$, 
	\[
	\begin{aligned}
		&\quad\int_U \tr[(CD^2f)^2]+\norm{Q_1^{-\frac{1}{2}}CDf}^2_U+((D^2\Phi_{1,(n,m)}+D^2\Phi_{2})CDf,CDf)_U\,\mathrm{d}\mu_1^{\Phi_{(n,m)}}\\
		&\leq  4 \int_U g_{(n,m)}^2\,\mathrm{d}\mu_1^{\Phi_{n,m}}.
		\end{aligned}
	\]
	Note that $\partial_{ij}\Phi_{1,(n,m)}$ and $\partial_{ij}\Phi_{2}$ exists in $L^1(\mu_1^{\Phi}),$ since $D\Phi_{1,(n,m)}$ is Lipschitz continuous and $D^2\Phi_2$ is $\mu^{\Phi}_1$ integrable by Assumption \nameref{ass:L_Phi-convex-reg-est}. 	
	Moreover, using that $\Phi_{1,(n,m)}$ is convex and Item $(ii)$ from Assumption \nameref{ass:L_Phi-convex-reg-est} we get
	\[
		\int_U \tr[(CD^2f)^2]+\norm{Q_1^{-\frac{1}{2}}CDf}^2_U\,\mathrm{d}\mu_1^{\Phi_{(n,m)}}\leq \left(4+\frac{c_{\Phi_2}}{4}\right) \int_U g_{(n,m)}^2\,\mathrm{d}\mu_1^{\Phi_{n,m}}.
	\]
	By \Cref{rem:weak_conv_measure} the left-hand side from the inequality above converges to 
	\[
	\int_U \tr[(CD^2f)^2]+\norm{Q_1^{-\frac{1}{2}}CDf}^2_U\,\mathrm{d}\mu_1^{\Phi}.
	\]
	Now observe
	\[
	\begin{aligned}
		\int_U g_{(n,m)}^2\,\mathrm{d}\mu_1^{\Phi_{n,m}}-\int_U g^2\,\mathrm{d}\mu_1^{\Phi}	&=\int_U (g_{(n,m)}^2-g^2)\frac{e^{-\Phi_{(n,m)}}}{\int_U e^{-\Phi_{(n,m)}}\,\mathrm{d}\mu_1}\,\mathrm{d}\mu_1\\&\qquad+(\mu_1^{\Phi_{n,m}}(g^2)-\mu_1^{\Phi}(g^2)).
	\end{aligned}
	\]
	By \Cref{rem:weak_conv_measure} the second term of the right-hand side of the inequality above converges to zero. By Assumption 			\nameref{hypoapproximationPotential} Item (iii) we know that 
	\[
	\lim_{(n,m)\rightarrow\infty}\norm{(Q^{\alpha}D\Phi_{1,(n,m)}-Q^{\alpha}D\Phi_{1},Q^{-\alpha}CDf)_U}_{L^2(U;\mu_1)}=0,
	\]
	as $Q^{-\alpha}CDf\in U_k$ for some $k\in\N$. Hence also $\lim_{(n,m)			\rightarrow\infty}g_{(n,m)}^2=g^2$ in $L^2(U;\mu_1)$. Using that $\frac{e^{-\Phi_{(n,m)}}}{\int_U e^{-\Phi_{(n,m)}}}$ is bounded 			independent of $(m,n)$ we can conclude that also the first term of the right-hand side of the inequality above converges to zero, 			which finishes the proof.
\end{proof}

The following now follows analogously to \cite[Lemmas~6 and 7]{EG21}:
\begin{prop}\label{prop:inf-dim-anti-collection}
	Let \nameref{ass:qc-negative-type} and \nameref{ass:ess-N_Phi} be satisfied. Then
	\begin{enumerate}[(i)]
		\item $H_1\subseteq D(S_0)$ with $S_0 P=0$.
		\item $P(\acore)\subseteq D(A^{\Phi}_0)$ and $A^{\Phi}_0Pf(u,v)=(-v,Q_2^{-1}K_{12}D_1(Pf)(u,v))$
			for all $f\in\acore$, $(u,v)\in W$.
		\item $PA^{\Phi}_0Pf=0$ for all $f\in\acore$.
		\item $A^{\Phi}_0P(\acore)\subseteq D(A^{\Phi}_0)$ and 
\[Gf\defeq P(A^{\Phi}_0)^2Pf=\tr[CD_1^2f_P]-(u,Q_1^{-1}CD_1f_P)_U-(D\Phi(u),CD_1f_P)_U\] for all $f\in\acore$. Moreover $(G,\acore)$ is essentially self-adjoint on $H$.
	\end{enumerate}
	In particular, the conditions \nameref{ass:data-op-decomposition}-\nameref{ass:data-anti-compliance} are satisfied, and $\acore$ is a 		core for the operator $(G,D(G))$ as defined in \Cref{def:paap-operator}.
\end{prop}
In particular, we can define the bounded operator $B$ on $H$ as in \Cref{def:b-operator}, which acts as $(\operatorname{Id}-G)^{-1}(A^{\Phi}_0P^*)$ on $D(A^{\Phi}_0P^*)$.
Now we can start verifying the hypocoercivity conditions, where we start with boundedness of the auxiliary operators $BA^{\Phi}_0(I-P)$ and $BS_0$.
The first part, works as in \cite[Proposition~6]{EG21}:
\begin{prop}
	The operator $(BA^{\Phi}_0(\operatorname{Id}-P),\acore)$ is bounded, so that the second inequality in \nameref{ass:aux-bound} holds with $c_2=\sqrt{8+\frac{c_{\Phi_2}}{4}}$.
\end{prop}
Now we introduce a new assumption ensuring boundedness of $BS_0$.
\begin{cond}{K7}\label{ass:inf-dim-matrix-bounded-ev}
	Let $K_{1},K_{2}(v)\in \lop{V}$ for each $v\in V$ and assume that they share the same invariance properties as $K_{22}$. Assume that $K_{22}(v)=K_{1} + K_{2}(v)$ for each $v\in V$.
	Further, let the following hold
	\begin{enumerate}[(i)]
		\item There is some $C_{1}\in(0,\infty)$ such that
			$\|Q_2^{-\frac12}K_{1} Q_2^{-\frac12}\|_{\lop{V}} \leq C_1$.
		\item There exists a measurable function $C_2:V\to [0,\infty)$ with
			\[
				\|Q_2^{-1}K_2(v)Q_2^{-\frac12}\|_{\lop{V}} \leq C_2(v)
			\]
			for almost all $v\in V$ and such that
			\[
				\overline{C}_{2}\defeq \int_V (C_2(v))^2\norm{v}_V^2\dm_2<\infty.
			\]
		\item 
			Assume that for all $v\in V$ the sequence $(\alpha_n^{22}(v))_{n\in\N}$ defined by
			\[\alpha_n^{22}(v)\defeq \sum_{k=1}^{\infty}(Q_2^{-\frac{1}{2}}\partial_kK_{22}(v)e_k,e_n)_V\] is 			in $\ell^2(\N)$
			and that $M_{22}\defeq \int_V\| ({\alpha_n^{22}}(v))_{n\in\N}\|^2_{\ell^2}\,\mu_2(\mathrm{d}v)<\infty$.
	\end{enumerate}
\end{cond}
\begin{remark}Recall $K_{22}^0$ from \nameref{ass:inf-dim-elliptic}, then $K_{22}(v)=K_{22}^0+K_{22}(v)-K_{22}^0$ is a possible decomposition of $K_{22}$ desired in \nameref{ass:inf-dim-matrix-bounded-ev}. But note, that in general $K_1$ has not to be positive definite, which is assumed for $K_{22}^0$.
\end{remark}
%with polynomial growth (i.e.~there is some $k_1,k_2\in\N$ such that $C_2(v)\leq k_1(1+\norm{v}_V^{k_2})$ for all $v\in V$),
\begin{prop}\label{prop:est_T}
	Let \nameref{ass:inf-dim-matrix-bounded-ev} hold.
	Then $(BS_0,\acore)$ is a bounded operator on $H$ and the first inequality in \nameref{ass:aux-bound} is satisfied for
	\[
		c_1\defeq \frac{1}{2}\big(\sqrt{C_{1}^2+\overline{C}_{2}}+\sqrt{M_{22}}\big).
	\]
\end{prop}
\begin{proof}
	We prove this via \Cref{lem:auxiliary-bound}, so let $f\in\acore$ and $h\in D(S_0)$ be arbitrary.
	Then by definition of $D(S_0)$, there is a sequence $\s{h}$ in $\acore$ such that $h_n\to h$ and $S_0h_n\to S_0 h$ in $H$ as 				$n\to\infty$.
	Fix some $n\in\N$, then
	\begin{equation}\label{eq:inf-dim-s-adj-dom}
	\begin{aligned}
		&\quad(S_0h_n,A^{\Phi}_0Pf)_{H}\\
		&= \int_W (D_2h_n(u,v),K_{22}(v)D_2 A^{\Phi}_0Pf(u,v))_V\,\mu^{\Phi}(\mathrm{d}(u,v)) \\
		&= -\int_U \int_V (D_2h_n(u,v),K_{22}(v)Q_2^{-1}K_{12}D_1f_P(u))_V\,\mu^{\Phi}(\mathrm{d}(u,v)) \\
		&= -\sum_{k\in\N} \int_U \int_V \partial_{2,k}h_n(u,v) (K_{22}(v)e_k,Q_2^{-1}K_{12}D_1f_P(u))_V \,\mu^{\Phi}(\mathrm{d}(u,v)).
	\end{aligned}
	\end{equation}
	Here, the first equality follows from the representation of $S$ in \Cref{lem:inf-dim-op-decomp},
	the second equality follows from \Cref{prop:inf-dim-anti-collection}~(ii),
	and the last line is due to symmetry of $K_{22}(v)$ for any $v\in V$.
	Note that the sum there is finite due to invariance properties of $K_{22}$ and $K_{12}$.
	Applying integration by parts (see \Cref{lem:inf-dim-ibp-fbs}), we obtain
	\[
	\begin{aligned}
		&\quad(S_0h_n,A^{\Phi}_0Pf)_{H}\\
		&= \sum_{k\in\N} \int_U \int_V h_n (\partial_kK_{22}(v)e_k,Q_2^{-1}K_{12}D_1f_P(u))_V \,\mu^{\Phi}(\mathrm{d}(u,v)) \\
		&\qquad - \sum_{k\in\N} \int_U \int_V (v,Q_2^{-1}e_k)_V h_n (e_k,K_{22}(v)Q_2^{-1}K_{12}D_1f_P(u))_V \,\mu^{\Phi}(\mathrm{d}(u,v))\\
		&= (h_n,Tf)_{{L^2(\mu^{\Phi})}},
	\end{aligned}
	\]
	where $T:\acore\to {L^2(\mu^{\Phi})}$ is defined by
	\[
		Tf(u,v)\defeq \sum_{k\in\N}(\partial_kK_{22}(v)e_k,Q_2^{-1}K_{12}D_1f_P(u))_V - (v,Q_2^{-1}K_{22}(v)Q_2^{-1}K_{12}D_1f_P(u))_V.
	\]
	Note that $Tf$ is indeed in ${L^2(\mu^{\Phi})}$, since all appearing sums are finite, $\|v\|_V\in L^2(\mu_2)$ and the properties of $K_{22}$.
	Moreover, since $1\in\fbs{B_W}$, it follows analogously to \eqref{eq:inf-dim-s-adj-dom} that
	$\mu^{\Phi}(Tf)=(1,Tf)_{L^2(\mu^{\Phi})}=(S1,A^{\Phi}_0Pf)_{L^2(\mu^{\Phi})}=0$, so $Tf\in H$. Now letting $n\to\infty$, we see that $A^{\Phi}_0Pf\in D(S_0^*)$ with $S_0^*A^{\Phi}_0Pf=Tf$.	
	This means that we are able to apply \Cref{lem:auxiliary-bound}, so we set $g\defeq (\operatorname{Id}-G)f$.
	We need to show that there is some $C_T<\infty$ such that
	\begin{equation}\label{eq:inf-dim-aux-bound-please}
		\|(BS_0)^*g\|_{L^2(\mu^{\Phi})}=\|S_0^*A^{\Phi}_0Pf\|_{L^2(\mu^{\Phi})}=\|Tf\|_{L^2(\mu^{\Phi})} \leq C_T \|g\|_{L^2(\mu^{\Phi})}
	\end{equation}
	holds for any choice of $f\in\acore$.
	Due to \nameref{ass:inf-dim-matrix-bounded-ev} Item (i), we have that
	\[
	\begin{aligned}
		&\quad\| (v,Q_2^{-1}K_{1}Q_2^{-1}K_{12}D_1 f_P)_V \|_H^2 \\
		&= \int_U (K_{1}Q_2^{-1}K_{12}D_1f_P(u), Q_2^{-1}K_{1}Q_2^{-1}K_{12}D_1f_P(u))_V\,\mathrm{d}\mu_1^\Phi \\
		&\leq C_1^2 \int_U (CD_1f_P,D_1 f_P)_V \,\mathrm{d}\mu_1^\Phi
		\leq \frac{C_1^2}{4} \int_U ((\operatorname{Id}-N)f_P)^2\,\mathrm{d}\mu_1^{\Phi}\\
		&=\frac{C_1^2}{4} \int_U \left(\int_V (\operatorname{Id}-G)f\,\mathrm{d}\mu_2\right)^2\,\mathrm{d}\mu_1^{\Phi}
		\leq \frac{C_1^2}{4} \|g\|_H^2, \\
	\end{aligned}
	\]
	where we applied \Cref{cor:gaussian-integral} and the estimate from \eqref{eq:regularity-estimate-one}. On the other hand, by Item (ii) from \nameref{ass:inf-dim-matrix-bounded-ev}. 
	\[
	\begin{aligned}
		&\quad\| (v,Q_2^{-1}K_{2}(v)Q_2^{-1}K_{12}D_1 f_P(u))_V \|_H^2 \\
		&\leq\int_W \|v\|_V^2 \|(Q_2^{-1}K_{2}(v)Q_2^{-\frac12}) Q_2^{-\frac12}K_{12}D_1 f_P(u)\|_V^2 \,\mathrm{d}\mu^\Phi\\
		&\leq\int_W \|v\|_V^2 C_2(v)^2\,
			\|Q_2^{-\frac12}K_{12}D_1f_P(u)\|_V^2\,\mathrm{d}\mu^\Phi
		\leq \overline{C}_{2} \int_U (CD_1f_P,D_1 f_P)_V\,\mathrm{d}\mu_1^{\Phi}\\
		&\leq\frac{\overline{C}_{2}}{4} \|g\|_H^2.
	\end{aligned}
	\]
	This shows that the second summand of $Tf$ can be bounded relatively to $g$. To deal with the first one, note that by Item (iii) from \nameref{ass:inf-dim-matrix-bounded-ev} we can estimate
	\[
	\begin{aligned}
		&\quad\left(\sum_{k\in\N}(\partial_kK_{22}(v)e_k,Q_2^{-1}K_{12}D_1f_P(u))_V  \right)^2\\
		&=\left(\sum_{k\in\N}\sum_{j\in\N}(Q_2^{-\frac{1}{2}}\partial_kK_{22}(v)e_k,e_j)_V(e_j,Q_2^{-\frac{1}{2}}K_{12}D_1f_P(u))_V  \right)^2\\
		&=\left(\sum_{j\in\N}\alpha_j^{22}(v)(e_j,Q_2^{-\frac{1}{2}}K_{12}D_1f_P(u))_V\right)^2\\
		&\leq\| ({\alpha_j^{22}}(v))_{j\in\N}\|^2_{\ell^2}(D_1f_P(u),CD_1f_P(u))_U.
	\end{aligned}
	\]
	Note that the right hand side factorizes into the $u$- and $v$-dependent components, so integration over $\mu^{\Phi}$
	yields a product of integrals with respect to $\mu_1^{\Phi}$ and $\mu_2$.	
	We obtain
	\[
		\int_W \left(\sum_{k\in\N}(\partial_kK_{22}(v)e_k,Q_2^{-1}K_{12}D_1f_P(u))_V  \right)^2\,\mathrm{d}\mu^{\Phi}\leq {M_{22}}\int_U (D_1f_P,CD_1f_P)_U\,\mathrm{d}\mu_1^{\Phi}.
	\]
	As above, this shows that the first summand of $Tf$ can also be bounded relative to $g$.
	Overall, we can see
	\[
		\|Tf\|_{L^2(\mu^{\Phi})} \leq \frac{1}{2}\left(\sqrt{C_1^2+\overline{C}_{2}}+\sqrt{M_{22}}\right)\|g\|_{L^2(\mu^{\Phi})}.
	\]
	This means that \eqref{eq:inf-dim-aux-bound-please} holds for $C_T=c_1$ as claimed, and \Cref{lem:auxiliary-bound}
	proves that $c_1$ is indeed an upper bound for the operator $BS_0$.	
\end{proof}

\begin{remark}\label{remark:alter_K_seven}
	In the proof, we have always used $(CD_1f_P,D_1f_P)_U$ as a bounding term, in order to apply the first inequality from \Cref{lem:inf-dim-regularity-estimates}.
	It seems clear that by involving eigenvalues of $Q_1$ into the assumptions \nameref{ass:inf-dim-matrix-bounded-ev},
	we can leverage all the invariance properties across finite-dimensional subspaces to instead use $(Q_1^{-1}CDf_P,CDf_P)_U$ as a bound.
	In that case, the aforementioned lemma enables us to the bound all terms relatively to $g$ if we assume the following alternative condition.
\end{remark}
\begin{cond}{K7*}\label{ass:inf-dim-matrix-bounded-ev_star}
	Recall \nameref{ass:inf-dim-elliptic} and let $K_{22}$ be such that $K_{22}(v)=K_{1} + K_{2}(v)$ for each $v\in V$. Assume 					\nameref{ass:qc-negative-type} such that $K_{12}$ is injective.
	Further, let the following hold
	\begin{enumerate}[(i)]
		\item There is some $C_{1}\in(0,\infty)$ such that
			$\|Q_2^{-\frac{1}{2}}K_1K_{12}^{-1} Q_1^{\frac{1}{2}}\|_{\lop{V}} \leq C_1$.
		\item There exists a measurable function $C_2:V\mapsto [0,\infty)$ with
			\[
				\|Q_1^{\frac{1}{2}}K_{12}^{-1} K_2(v)Q_2^{-1}\|_{\lop{V}} \leq C_2(v)
			\]
			for almost all $v\in V$ and such that
			\[
				\overline{C}_{2}\defeq \int_V (C_2(v))^2\norm{v}_V^2\dm_2<\infty.
			\]
		\item 
			Assume that the sequence $\s{\alpha^{22}(v)}$ defined by
			\[\alpha_n^{22}(v)\defeq \sum_{k=1}^{\infty}(Q_1^{\frac{1}{2}}K_{12}^{-1}\partial_kK_{22}(v)e_k,d_n)_V\] is 				an element of $\ell^2(\N)$
			and that $M_{22}\defeq \int_V\| ({\alpha_n^{22}}(v))_{n\in\N}\|^2_{\ell^2}\,\mu_2(\mathrm{d}v)<\infty$.
	\end{enumerate}
\end{cond}
Now that the first hypocoercivity condition is proven, we are left to show \nameref{ass:micro-coerc} and \nameref{ass:macro-coerc}.
For this, we assume modified Poincar\'{e} inequalities based on $K_{22}$ and $K_{12}$.
\begin{cond}{K8}\label{ass:inf-dim-v-poincare}
	Assume that there is some $c_S\in(0,\infty)$ such that
	\[
		\int_V (K_{22}D_2f,D_2f)_V\,\mathrm{d}\mu_2 \geq c_S \int_V \left(f-\mu_2(f)\right)^2\,\mathrm{d}\mu_2\quad \text{for all}\quad f\in\fbs{B_V}.
	\]
	
\end{cond}
\begin{cond}{K9}\label{ass:inf-dim-u-poincare}
	Assume that there is some $c_A\in(0,\infty)$ such that
	\[
		\int_U (Q_2^{-1}K_{12}D_1f,K_{12}D_1f)_V\,\mathrm{d}\mu^{\Phi}_1 \geq c_A \int_U \left(f-\mu_1(f)\right)^2\,\mathrm{d}\mu^{\Phi}_1\quad \text{for all}\quad f\in\fbs{B_U}.
	\]
\end{cond}

\begin{remark}\label{rem:inf-dim-eigenvalue-estimates}
	\begin{enumerate}[(i)]
		\item Recall $K_{22}^0$ from \nameref{ass:inf-dim-elliptic} and assume that there is a sequence of eigenvalues $(\lambda_k^0)_{k\in\N}$ of $K_{22}^0$ with respect to the basis $B_V$. Let $\lambda_{2,i}$ denote the $i$-th eigenvalue of $Q_2$,
			then due to \Cref{lem:inf-dim-poincare}, we have
			\[
				\int_V (Q_2D_2f,D_2f)_V\,\mathrm{d}\mu_2 \geq \lambda_{2,1} \int_V \left(f-\mu_2(f)\right)^2\,\mathrm{d}\mu_2\quad \text{for all}\quad f\in\fbs{B_V}.
			\]
			So if there is some $\omega_{22}\in(0,\infty)$ such that $\lambda_k^0\geq \omega_{22}\lambda_{2,k}$ 				for each $k\in\N$, then \nameref{ass:inf-dim-v-poincare} holds with
			$c_S=\frac{\lambda_{2,1}}{\omega_{22}}$.
		\item Similarly, if $\Phi=\Phi_1+\Phi_2$ is as described in Assumptions \nameref{ass:L_Phi-convex-reg-est}   and there is some $				\omega_{12}\in(0,\infty)$ such that $(Q_2^{-1}K_{12}d_k,K_{12}d_k)\geq \omega_{12}\lambda_{1,k}$
			for all $k\in\N$, then \nameref{ass:inf-dim-u-poincare} holds with $c_A= \frac{\lambda_{1,1}}{\omega_{12}e^{\norm{\Phi_2}_{osc}}}$ by \Cref{lem:inf-dim-poincare}.
	\end{enumerate}
\end{remark}

Under these conditions, we can easily verify macroscopic and microscopic coercivity.
\begin{prop}
	Let \nameref{ass:inf-dim-v-poincare} hold. Then $S_0$ satisfies \nameref{ass:micro-coerc} with $\Lambda_m=c_S$.
\end{prop}
\begin{proof}
	Let $f\in\acore$ and set $f_u\defeq f(u,\cdot)-Pf(u)\in\fbs{B_V}$ for any $u\in U$.
	Then $\mu_2(f_u)=0$ and $D_2f_u(v)=D_2f(u,v)$ for all $u\in U$, $v\in V$.
	By \nameref{ass:inf-dim-v-poincare} and \Cref{lem:inf-dim-op-decomp}, it then holds that
	\[
	\begin{aligned}
		c_S\|(\operatorname{Id}-P)f\|_H^2
		&= c_S\int_U \int_V f_u^2\,\mathrm{d}\mu_2\,\mathrm{d}\mu^{\Phi}_1
		\leq \int_U\int_V(K_{22}(v)D_2f_u,D_2f_u)_V\,\mathrm{d}\mu_2\,\mathrm{d}\mu^{\Phi}_1 \\
		&= \int_W (K_{22}D_2f,D_2f)_V\,\mathrm{d}\mu^{\Phi}
		=-(S_0f,f)_{H}.
	\end{aligned}
	\]
\end{proof}

\begin{prop}
	Let \nameref{ass:inf-dim-u-poincare} hold. Then $A^{\Phi}_0$ satisfies \nameref{ass:macro-coerc} with $\Lambda_M=c_A$.
\end{prop}
\begin{proof}
	Let $f\in\acore$, then $f_P\in\fbs{B_U}$ with $\mu^{\Phi}_1(f_P)=0$.
	By using \nameref{ass:inf-dim-u-poincare}, then \Cref{cor:gaussian-integral}, and finally \Cref{prop:inf-dim-anti-collection}~(ii),
	it follows that
	\[
	\begin{aligned}
		c_A\|Pf\|_H^2
		&= c_A\int_U f_P^2\,\mathrm{d}\mu_1^{\Phi}
		\leq \int_U (Q_2^{-1}K_{12}D_1f_P,K_{12}D_1f_P)_V\,\mathrm{d}\mu_1^{\Phi} \\
		&= \int_U\int_V (v,Q_2^{-1}K_{12}D_1f_P)_V^2\,\mathrm{d}\mu_2\,\mathrm{d}\mu_1^{\Phi}
		=\|A^{\Phi}_0Pf\|_{H}^2.
	\end{aligned}
	\]Hence the claim follows by \Cref{prop:inf-dim-anti-collection} and \Cref{remark:core_if_ess_self}.
\end{proof}

The main result of this section is now immediate.
\begin{thm}\label{thm:inf-dim-hypoc-applied}
	Let the conditions \nameref{ass:inf-dim-elliptic}-\nameref{ass:pert-pot}, \nameref{ass:qc-negative-type}-\nameref{ass:inf-dim-u-poincare} (either \nameref{ass:inf-dim-matrix-bounded-ev} or \nameref{ass:inf-dim-matrix-bounded-ev_star}), \nameref{ass:ess-N_Phi} and 		\nameref{ass:L_Phi-convex-reg-est} hold.
	Then the semigroup $\sccs$ on $L^2(\mu^{\Phi})$ generated by the closure $(L^{\Phi},D(L^{\Phi}))$ of $(L^{\Phi},\fbs{B_W})$
	is hypocoercive in the sense that for each $\theta_1\in (1,\infty)$,
	there is some $\theta_2\in (0,\infty)$ such that
	\[
		\left\| T_tf- \mu^{\Phi}(f) \right\|_{L^2(\mu^{\Phi})} \leq \theta_1 \mathrm{e}^{-\theta_2 t}\left\| f-\mu^{\Phi}(f)  \right\|_{L^2(\mu^{\Phi})}
	\]
	for all $f\in {L^2(\mu^{\Phi})}$ and all $t\geq 0$.
	For $\theta_1\in (1,\infty)$ the constant $\theta_2$ determining the speed of convergence can be explicitly computed in terms of 			$c_S$, $c_A$, $c_{\Phi_2}$, and $c_1$ as
	\[
		\theta_2=\frac{1}{2}\frac{\theta_1-1}{\theta_1}\frac{\min\{c_S,c_1\}}{(1+2\sqrt{2}+c_1+\frac{c_{\Phi_2}}{4})\Big(1+\frac{1+c_A}				{2c_A}(1+c_1+\sqrt{8+\frac{c_{\Phi_2}}{4}})\Big)+ \frac{1}{2}\frac{c_A}{1+c_A}}\frac{c_A}{1+c_A}.
	\]	
	 Moreover, the transition semigroup $(p_t)_{t\geq 0}$ associated to the stochastic process $\mathbf{M}=(\Omega,\mathcal{F},(\mathcal{F}_t)_{t\geq 0}, (X_t,Y_t)_{t\geq 0},(P_w)_{w\in W})$ is hypocoercive with the rate computed in \Cref{thm:inf-dim-hypoc-applied}, which shows exponential convergence of the 	weak solution to the equilibrium described by the invariant measure $\mu^{\Phi}$.
\end{thm}
\begin{proof}
	By the prior considerations, \Cref{thm:abstract-hypoc} can be applied to $\sccsz$ to yield
	\[
		\| T_t^0 f\|_H \leq \kappa_1\mathrm{e}^{-\kappa_2 t} \|f\|_H
		\qquad\text{ for all }f\in H,\ t\geq 0
	\]
	for suitable $\kappa_1,\,\kappa_2\in (0,\infty)$. By conservativity and $\mu^{\Phi}$-invariance of $\sccs$ (see \Cref{thm:inf-dim-ess-m-diss}), this implies
	\[
		\| T_t f - \mu^{\Phi}(T_t f)\|_{L^2(\mu^{\Phi})}
		= \| T_t (f-\mu^{\Phi}(f))\|_H
		\leq \kappa_1\mathrm{e}^{-\kappa_2 t} \|f-\mu^{\Phi}(f)\|_{L^2(\mu^{\Phi})}
	\]
	for all $f\in {L^2(\mu^{\Phi})}$ and $t\geq 0$. The computation of $\theta_1 $ and $\theta_2$ can then be done as in the proof of 			\cite[Theorem 5.2]{EG21_Pr}, see also \cite[Theorem 1.1]{GS16}.
	By association, $(p_t)_{t\geq 0}$ is a $\mu^{\Phi}$-version of $\sccs$. Therefore $(p_t)_{t\geq 0}$ satisfies the estimate from 			\Cref{thm:inf-dim-hypoc-applied} and the second statement is shown.
\end{proof}

We end this section with an $L^2$-exponential ergodicity result for the associated right process\[\mathbf{M}=(\Omega,\mathcal{F},(\mathcal{F}_t)_{t\geq 0}, (X_t,Y_t)_{t\geq 0},(P_w)_{w\in W}),\] which was already shown in \cite[Corollary~6.12]{EG21_Pr}. We want to highlight that Assumptions \nameref{ass:inf-dim-diff-bound} and \nameref{ass:weak_sol} are not involved to study the longtime behavior of $\mathbf{M}$ in the sense of the following corollary.

\begin{corollary}\label{ergodic}
	Assume we are in the situation as in the \Cref{thm:inf-dim-hypoc-applied} with constants $\theta_1\in (1,\infty)$ and $\theta_2\in (0,\infty)$ determined in \Cref{thm:inf-dim-hypoc-applied}.
	Then for all $t\in (0,\infty)$ and $f\in L^2(W;\mu^{\Phi})$ it holds
	\[
	\begin{aligned}
		&\quad\left\Vert\frac{1}{t}\int_{[0,t)}f(X_s,Y_s)\,\lambda(\mathrm{d}s)-(f,1)_{L^2(\mu^{\Phi})}\right\Vert_{L^2({P}_{\mu^{\Phi}})}\\
		&\leq \frac{1}{\sqrt{t}}\sqrt{\frac{2\theta_1}{\theta_2}\big(1-\frac{1}{t\theta_2}(1-e^{-t\theta_2})\big)}\norm{f-(f,						1)_{L^2(\mu^{\Phi})}}_{L^2(\mu^{\Phi})}.
	\end{aligned}
	\]
	We call $\mathbf{M}=(\Omega,\mathcal{F},(\mathcal{F}_t)_{t\geq 0}, (X_t,Y_t)_{t\geq 0},(P_w)_{w\in W})$ with this property $L^2$-			exponentially ergodic, i.e.~ergodic with a rate that corresponds to exponential convergence of the corresponding semigroup.
\end{corollary}

\begin{remark}
	Roughly speaking, the result above means that time average converges to space average in $L^2({P}_{\mu^{\Phi}})$ with rate $t^{-\frac{1}{2}}$. If the spectrum of $(L^{\Phi},D(L^{\Phi}))$ contains a negative eigenvalue $-\kappa$ with corresponding eigenvector $f$, 	then this rate is optimal. Indeed, in this case one can show, compare \cite[Remark~6.13]{EG21_Pr}
	\[
	\begin{aligned}
		\left\Vert\frac{1}{t}\int_{[0,t)}f(X_s)\,\lambda(\mathrm{d}s)\right\Vert_{L^2({P}_{\mu^{\Phi}})}
		=\frac{1}{\sqrt{t}}\sqrt{\frac{2}{\kappa}\big(1-\frac{1}{t\kappa}(1-e^{-t\kappa})\big)}\norm{f}_{L^2(\mu^{\Phi})}
	\end{aligned}
	\]
	for all $t\in (0,\infty)$. 
\end{remark}

\section{Examples}\label{sec:examples}
In this section we have a look at certain examples, where the results we derived above are applicable. The examples are inspired by the ones in \cite[Section~5~and~6]{Lunardi14}.
\subsection{Stochastic reaction-diffusion equations}

Let $\mathrm{d}\xi$ be the standard Lebesgue measure on $((0,1),\mathcal{B}(0,1))$. Set $U=L^2(0,1)\defeq L^2((0,1),\mathrm{d}\xi)$ and $W=U\times U$. Moreover, let
$(-\frac{\partial^2}{\partial\xi^2}, D(\frac{\partial^2}{\partial\xi^2}))$ be the negative second order derivative with Dirichlet boundary conditions, i.e. 
\[
	D(\frac{\partial^2}{\partial\xi^2})=W^{1,2}_0(0,1)\cap W^{2,2}(0,1)\subseteq L^2(0,1).
\]Consider the linear continuous positive self-adjoint operator
\[
	Q=(-\frac{\partial^2}{\partial\xi^2})^{-1}:L^2(0,1)\rightarrow D(\frac{\partial^2}{\partial\xi^2}).
\]
Denote by $B_U=(e_k)_{k\in\N}=(\sqrt{2}\sin(k\pi\cdot))_{k\in\N}$ the orthonormal basis of $L^2(0,1)$ diagonalizing $Q$ with corresponding eigenvalues $(\lambda_k)_{k\in\N}=(\frac{1}{(k\pi)^2})_{k\in\N}$.
For $\alpha_1,\alpha_2\in (\frac{1}{2},\infty)$ we consider two centered non-degenerate infinite-dimensional Gaussian measures $\mu_1$ and $\mu_2$ on $(U,\mathcal{B}(U))$, with covariance operators
\[
	Q_1=Q^{\alpha_1}\quad\text{and}\quad Q_2=Q^{\alpha_2},
\]
respectively.
Since $(\lambda_k)_{k\in\N}\in \ell^{\theta}(\N)$ for $\theta\in (\frac{1}{2},\infty)$, $Q_1$ and $Q_2$ are indeed trace class.
The eigenvalues of $Q_1$ and $Q_2$ are given by
\[
	\lambda_{1,k}\defeq \lambda_k^{\alpha_1}\quad\text{and}\quad\lambda_{2,k}\defeq \lambda_k^{\alpha_2},\quad k\in\N,
\]
respectively.
For $\sigma_1\in[0,\infty)$, we choose $K_{12}=Q^{\sigma_1}$ and since $K_{21}=K_{12}^*$ also $K_{21}=Q^{\sigma_1}$. Moreover, we define $K_{22}$ by specifying its eigenvalue functions $\lambda_{22,k}:V\to\R$. Let $c,\tilde{c}\in (0,\infty)$ be constant.
For each $k\in\N$, let $\sigma_2\in[0,\infty)$, $\beta_k\in(0,1)$, $\varphi_k\in C^1(\R;[0,\infty))$ and $\psi_k\in C^1(\R^k;[0,\infty))$, both with bounded derivative.

Define
\[
	\lambda_{22,k}(v)\defeq c\lambda_{k}^{\sigma_2}+ \gamma_k\left(\varphi_k(|p_k^U v|^{\beta_k+1})+\psi_k(p_k^U v)\right),
\]
where
\[
	\gamma_k\defeq \frac{\tilde{c}\lambda_{k}^{\sigma_2}}{\|\varphi'_k\|_{\infty}+\|D\psi_k\|_{\infty}+\abs{\varphi_k(0)}+\abs{\psi_k(0)}}.
\]
One can check that $ c\lambda_{k}^{\sigma_2}\leq \lambda_{22,k}(v)=\lambda_{22,k}(P_k^U v)\leq A \lambda_{k}^{\sigma_2}(1+|p_k^U v|^{2})\leq A \lambda_{k}^{\sigma_2}(1+\norm{v}_U^{2})$ for all $k\in\N$, $v\in U$ and some $A\in [0,\infty)$ independent of $k$.
Moreover, for $i\geq k$, we have $\partial_i\lambda_{22,k}(v)=0$, and for $1\leq i\leq k$, it holds that
\[
\begin{aligned}
	|\partial_i\lambda_{22,k}(v)|
	&= \gamma_k \left|\varphi_k'(|p_k^U v|^{\beta_k+1})(\beta_k+1)|p_k^U v|^{\beta_k-1}(v,e_i)_U + \partial_i \psi(p_k^U v)\right| \\
	&\leq \gamma_k(\beta_k+1)(\|\varphi'_k\|_{C^0}+\|D\psi_k\|_{C^0})(1+|p_k^U v|^{\beta_k}) \\
	&\leq 2\tilde{c}\lambda_{k}^{\sigma_2}(1+\|P_k^U v\|_U^{\beta_k})
\end{aligned}
\]
for all $v\in U$. Now we simply set $K_{22}(v)e_i\defeq \lambda_{22,k}(v)e_i$, which describes a symmetric positive-definite bounded linear operator on $U$
as required for \Cref{def:inf-dim-operators}.
Moreover, \nameref{ass:inf-dim-elliptic} holds for $K_{22}^0=cQ^{\sigma_2}$, \nameref{ass:inf-dim-growth} is satisfied for $N_k\defeq 2\tilde{c}\lambda_{k}^{\sigma_2}$.

We fix a continuous differentiable (convex) function $\phi:\R\rightarrow \R$, which is bounded from below an with bounded derivative.
	For such $\phi$ we consider potentials $\Phi:L^2(0,1)\rightarrow \R$ defined by
	\[
		\Phi\defeq \Phi_1+\Phi_2\quad\text{where}\quad \Phi_1(u)=\int_{(0,1)}\phi( u(\xi))\,\mathrm{d}\xi 
	\]and $\Phi_2:U\rightarrow \R$ is bounded and two times continuously Fr\'{e}chet-differentiable with bounded first and second order derivative.
\begin{remark}\label{der}
	Note the boundedness of $\phi'$ implies that $\phi$ grows at most linear.
	Moreover, $\Phi_1$ is lower semicontinuous by Fatou's lemma, bounded from below and in $L^p(U;\mu_1)$ for all $p\in [1,\infty)$. If $		\phi$ is convex the same holds true for $\Phi_1$. Using \cite[Proposition 5.2]{Lunardi14} we know that $\Phi$  in $W^{1,2}(U;\mu_1)$ 		with $D\Phi(u)=\phi'\circ u+D\Phi_2(u)$ for $u\in L^2(0,1)$ and that $\Phi$ fulfills Assumption \nameref{hypoapproximationPotential}.
	In particular, we obtain 
	\[
		\norm{D\Phi}_{L^{\infty}(\mu_1)}\leq\sup_{t\in\R}\abs{\phi'(t)}+\norm{D\Phi_2}_{L^{\infty}(\mu_1)}<\infty.
	\]
\end{remark}
If we assume
\[
	\sigma_2\leq 2\sigma_1,
\]
then \nameref{ass:pert-pot} is valid with $\alpha=0$. Therefore all assumptions for essential m-dissipativity are fulfilled and we obtain essential m-dissipativity on $L^2(L^{2}(0,1)\times L^{2}(0,1);\mu^{\Phi})$ of $(L_{\Phi},\fbs{B_W})$. As explained in Section \ref{sec:process} there is an associated  right process 
\[
	\mathbf{M}=(\Omega,\mathcal{F},(\mathcal{F}_t)_{t\geq 0}, (X_t,Y_t)_{t\geq 0},(P_w)_{w\in W})
\]with infinite life time, such that its transition semigroup and resolvent coincides in $L^2(W;\mu^{\Phi})$ with the ones generated by $(L_{\Phi},\fbs{B_W})$.
To show that the process has weak continuous paths we invoke \Cref{rem:inf-dim-bounding-rho} and assume
\[
	-\frac{\alpha_1}{2}+\sigma_1+\frac{\alpha_2}{2}>\frac{1}{2}\quad\text{and}\quad -\frac{\alpha_2}{2}+\sigma_1+\frac{\alpha_1}{2}				>\frac{1}{2}
\] 
to verify Item (i) from \nameref{ass:inf-dim-diff-bound} as well as $\sigma_2\in (\frac{1}{2},\infty)$ and $\sigma_2-\alpha_2\in (\frac{1}{2},\infty)$ to verify Item (ii) and (iii) from \nameref{ass:inf-dim-diff-bound}, respectively.
If we additionally want to construct a stochastically and analytically weak solution with weakly continuous paths in the sense of \Cref{thm:sto_ana_weak_sol} we need to assume $\sigma_2\in(\frac{1}{2},\infty)$ in order to verify Assumption \nameref{ass:weak_sol}. This is redundant since we already assume that $\sigma_2-\alpha_2\in (\frac{1}{2},\infty)$.

Concerning Hypocoercivity, we note that Assumption \nameref{ass:qc-negative-type} is obviously valid. But to obtain the final hypocoercivity result we distinguish to cases. First suppose $\sigma_2=\frac{3}{2}\alpha_2$.
Since we assume $\sigma_2\leq 2\sigma_1$ and $\sigma_2=\frac{3}{2}\alpha_2$ we know that $C=Q^{2\sigma_1-\alpha_2}$ is bounded. Hence Assumption \nameref{ass:ess-N_Phi} and \nameref{ass:L_Phi-convex-reg-est} are valid.

Item (i) and (ii) from assumption
\nameref{ass:inf-dim-matrix-bounded-ev} hold for 
\[
	C_1=0\quad\text{and}\quad C_{2}(v)=A 	(1+\norm{v}_U^{2})\quad v\in U.
\]
For Item (iii) note that $\alpha_{n}^{22}(v)\leq 2\tilde{c}(2+\norm{v}_U)\lambda_{n}^{\frac{3\alpha_2}{2}-\frac{\alpha_2}{2}}$, for all $v\in U$, which describes an $\ell^2$-sequence, since
\[
	\frac{3\alpha_2}{2}-\frac{\alpha_2}{2}= {\alpha_2}> \frac{1}{2}.
\]
Moreover
\[
	\int_U  \|({\alpha_n^{22}}(v))_{n\in\N}\|^2_{\ell^2}\,\mu_2(\mathrm{d}v)\leq 4\tilde{c}^2\|(\lambda_{n}^{\alpha_2})_{n\in\N}\|			^2_{\ell^2}\int_U(2+\norm{v}_U)^2\,\mu_2(\mathrm{d}v)<\infty.
\]
  
Second, as described in \Cref{remark:alter_K_seven} we can also verify Assumption \nameref{ass:inf-dim-matrix-bounded-ev_star}. 

Item (i) and (ii) from assumption
\nameref{ass:inf-dim-matrix-bounded-ev_star} hold for 
\[
	-{\alpha_2}+\sigma_2-\sigma_1+\frac{\alpha_1}{2}>0\quad\text{with}\quad C_1=0\quad\text{and}\quad C_{2}(v)=A 	(1+\norm{v}_U^{2})\quad v\in U.
	\]
For Item (iii) note that $\alpha_{n}^{22}(v)\leq 2\tilde{c}(2+\norm{v}_U)\lambda_{n}^{\frac{\alpha_1}{2}-\sigma_1+\sigma_2}$, for all $v\in U$, which describes an $\ell^2$-sequence, since by the inequality above
\[
	\frac{\alpha_1}{2}-\sigma_1+\sigma_2> \alpha_{2}>\frac{1}{2}.
\]
Moreover, 
\[
	\int_U  \|({\alpha_n^{22}}(v))_{n\in\N}\|^2_{\ell^2}\,\mu_2(\mathrm{d}v)\leq 4\tilde{c}^2\|(\lambda_{n}^{\frac{\alpha_1}{2}-\sigma_1+\sigma_2})_{n\in\N}\|^2_{\ell^2}\int_U(2+\norm{v}_U)^2\,\mu_2(\mathrm{d}v)<\infty.
\]
Since we want to allow unbounded $(C,D(C))$ in this situation, we assume $\Phi_2=0$ and scale $\Phi_1$ such that $\norm{Q_1^{\frac{1}{2}}D\Phi}_{L^{\infty}(\mu_1)}<\frac{1}{2}$ to verify assumption \nameref{ass:ess-N_Phi} and \nameref{ass:L_Phi-convex-reg-est}.
Finally, in both cases, using \Cref{rem:inf-dim-eigenvalue-estimates} we obtain \nameref{ass:inf-dim-v-poincare}
and \nameref{ass:inf-dim-u-poincare}, if
\[
	2\sigma_1-\alpha_2\leq \alpha_1.
\]

\begin{summary}\label{rem:table_rde_one}
For parameters $\alpha_1,\alpha_2\in (\frac{1}{2},\infty)$, $\sigma_1,\sigma_2\in [0,\infty)$ with $\sigma_2\leq 2\sigma_1$ and $\Phi$ with $\norm{D\Phi}_{L^{\infty}(\mu_1)}<\infty$ (i.e.~\nameref{ass:inf-dim-elliptic}-\nameref{ass:pert-pot} hold), we know that the generator $L^{\Phi}$ of the degenerate stochastic reaction-diffusion type equation
\[
\begin{aligned}
	\mathrm{d}X_t&=(-\frac{\partial^2}{\partial\xi^2})^{-\sigma_1+\alpha_2}Y_t\,\mathrm{d}t\\
	\mathrm{d}Y_t&=\sum_{i=1}^\infty \partial_iK_{22}(Y_t)f_i-K_{22}(Y_t)(-\frac{\partial^2}{\partial\xi^2})^{\alpha_2}Y_t-(-					\frac{\partial^2}{\partial\xi^2})^{-\sigma_1+\alpha_1} X_t\\
	&\qquad-(-\frac{\partial^2}{\partial\xi^2})^{-\sigma_1}\phi'(X_t)-(-\frac{\partial^2}{\partial\xi^2})^{-\sigma_1}D\Phi_2(X_t)\,\mathrm{d}t
	+\sqrt{2K_{22}(Y_t)}\,\mathrm{d}B_t
\end{aligned}
\]
is essentially m-dissipative. Moreover, as $L^{\Phi}$ is sub-Markovian, conservative and $\mu^{\Phi}$ invariant there exists a corresponding right process 
\[
	\mathbf{M}=(\Omega,\mathcal{F},(\mathcal{F}_t)_{t\geq 0}, (X_t,Y_t)_{t\geq 0},(P_w)_{w\in W})
\]
with infinite life time $P_{\mu^{\Phi}}$ a.s., solving the martingale problem $(L^{\Phi},D(L^{\Phi})$ under $P_{\mu^{\Phi}}$. We want to highlight that $K_{22}:L^2(0,1)\to \lop{L^2(0,1)}$ is not finitely based. Further, note that $\sigma_2=0$ is a valid choice. In this case the variable diffusion matrix $K_{22}$ is not trace class.
The table below summarize the conditions on the parameters and the potential such that $\mathbf{M}$ provides a stochastically and analytically weak solution with weak continuous paths to the equation above and that the semigroup $\sccs$ generated by $(L^{\Phi},D(L^{\Phi}))$ is hypocoercive.

\smallskip

\begin{center}
\begin{tabular}{ | m{5,4em} | m{4,4cm}| m{5,25cm} | } 
\hline
 & $\mathbf{M}$ stochastically and analytically weak solution with weak continuous paths & $\sccs$ is hypocoercive \\
\hline
\nameref{ass:inf-dim-diff-bound} & $\sigma_2-\alpha_2$, $\pm\nicefrac{\alpha_1}{2}+\sigma_1\mp\nicefrac{\alpha_2}{2}>\nicefrac{1}{2}$&  \\
\hline 
\nameref{ass:weak_sol} & follows by \nameref{ass:inf-dim-diff-bound} &\\
\hline
\nameref{ass:qc-negative-type} &  by construction &\\
\hline
\nameref{ass:ess-N_Phi}, \nameref{ass:L_Phi-convex-reg-est}, \nameref{ass:inf-dim-matrix-bounded-ev} &  & $\sigma_2=\nicefrac{3}{2}\alpha_2$\\
\hline
\nameref{ass:ess-N_Phi}, \nameref{ass:L_Phi-convex-reg-est}, \nameref{ass:inf-dim-matrix-bounded-ev_star}&  & $-{\alpha_2}+\sigma_2-\sigma_1+\nicefrac{\alpha_1}{2}>0$, $\Phi_2=0$ and $\norm{D\Phi}_{L^{\infty}(\mu_1)}<\nicefrac{1}{2}$\\
\hline
\nameref{ass:inf-dim-v-poincare}, \nameref{ass:inf-dim-u-poincare}&  & $2\sigma_1-\alpha_2\leq \alpha_1$\\
\hline
\end{tabular}
\end{center}

\smallskip

As explained in \Cref{ergodic}, we can combine the results described above to verify that $\textbf{M}$ is $L^2$-exponentially ergodic.
Moreover, choosing $\sigma_1=\alpha_2$ and the other parameters accordingly, we are able rewrite the system of coupled equations above, as an infinite-dimensional second order in time stochastic differential equation.
\end{summary}
\subsection{Stochastic Cahn-Hilliard equations}
Set
\[
	X\defeq\lbrace x\in W^{1,2}(0,1)\mid \int_0^1x(\xi)\,\mathrm{d}\xi=0\rbrace
\]
and equip it with the inner product $(\cdot,\cdot)_X$ given by
\[
	(x,y)_X\defeq\int_0^1\frac{\partial}{\partial\xi}x(\xi)\frac{\partial}{\partial\xi}y(\xi)\,\mathrm{d}\xi.
\]
Now let $U$ be the dual space of $X$, endowed with the canonical dual inner product and norm. Further, set $W=U\times U$ and for $p\in[1,\infty)$, $\tilde{L}^p(0,1)=\lbrace x\in {L}^p(0,1)\mid \int_0^1x(\xi)\,\mathrm{d}\xi=0\rbrace$. In the following we consider $\tilde{L}^p(0,1)$ as a subspace of $U$ by identifying an element $x\in \tilde{L}^p(0,1)$ with the continuous linear functional $y\mapsto \int_0^1x(\xi)y(\xi)\,\mathrm{d}\xi$ in $U$.
Next we define the map 
\[
	B:X\to X'=U,\quad Bx(y)=\int_0^1\frac{\partial}{\partial\xi}x(\xi)\frac{\partial}{\partial\xi}y(\xi)\,\mathrm{d}\xi,\quad y\in X.
\]
Note that for $x\in W^{2,2}(0,1)\cap X$ with Neumann boundary conditions we have
\[
	Bx(y)=-\int_0^1\frac{\partial^2}{\partial\xi^2}x(\xi)y(\xi)\,\mathrm{d}\xi,\quad y\in X.
\]
I.e.~$B$ can be identified with the extension of minus the second order derivative with Neumann boundary conditions. One can show that $B$ is isometric and we have for $z\in \tilde{L}^2(0,1)$ and $x\in X$, $(z,Bx)_U=(z,x)_{L^2(0,1)}$. It is well known that $(e_k)_{k\in\N}=(\sqrt{2}\cos(k\pi\cdot))_{k\in\N}$ is an orthonormal basis of $\tilde{L}^2(0,1)$ with $Be_k=(k\pi)^2e_k$. Therefore $(f_k)_{k\in\N}$ defined by $f_k=k\pi e_k$ is an orthonormal basis of $U$.
The operator $-B^2:D(B^2)\to U$ is a self-adjoint realization of the negative fourth order derivative with zero boundary conditions for the first and third order derivative. Moreover we have
\[
	(-B^2)^{-1}f_k=-\frac{1}{(\pi k)^4}f_k.
\]
Therefore $Q=-(-B^2)^{-1}:U\mapsto U$ is a positive self-adjoint operator, with corresponding orthonormal basis $(f_k)_{k\in\N}$ and eigenvalues $(\lambda_k)_{k\in\N}=(\frac{1}{(\pi k)^4})_{k\in\N}$. Let $\alpha_1,\alpha_2\in(\frac{1}{4},\infty)$. Since the sequence of eigenvalues of $Q$ is in $l^{\theta}(\N)$ for all $\theta \in (\frac{1}{4},\infty)$ it is reasonable to consider 
\[
	Q_1=Q^{\alpha_1}\quad \text{and}\quad Q_2=Q^{\alpha_2}
\] 
as covariance operators for the infinite-dimensional Gaussian measures $\mu_1$ and $\mu_2$, respectively. Moreover,
\[
	\lambda_{1,k}\defeq \lambda_k^{\alpha_1}\quad \quad\text{and}\quad \lambda_{2,k}\defeq \lambda_k^{\alpha_2},\quad k\in\N
\]
are the Eigenvalues of $Q_1$ and $Q_2$, respectively. As in the reaction-diffusion setting we choose
$K_{12}=Q^{\sigma_1}$ for some $\sigma_1\in [0,\infty)$. Since $K_{21}=K_{12}^*$ also $K_{21}=Q^{\sigma_1}$. We define $K_{22}$ by specifying its eigenvalue functions $\lambda_{22,k}:U\to\R$. Let $c,\tilde{c}\in (0,\infty)$ be constant.
Fix $\sigma_2\in[0,\infty)$ and for each $k\in\N$,  $\beta_k\in(0,1)$, $\varphi_k\in C^1(\R;[0,\infty))$ and $\psi_k\in C^1(\R^k;[0,\infty))$, both with bounded derivative.
Define $\lambda_{22,k}(v)$ and $K_{22}(v)$ as in the previous example. Therefore by the exact same reasoning as before, \nameref{ass:inf-dim-elliptic} holds for $K_{22}^0=cQ^{\sigma_2}$ and \nameref{ass:inf-dim-growth} is satisfied for $N_k\defeq 2\tilde{c}\lambda_{k}^{\sigma_2}$.
We fix a (convex) function $\phi\in C^1(\R)$, which is bounded from below. We assume that there are constants $A_1\in (0,\infty)$ and $p\in [1,\infty)$ such that
\[
	\abs{\phi'(t)}\leq A_1(1+\abs{t}^{p-1}),\quad t\in\R.
\]
I.e.~$\phi$ and its derivative have at most polynomial growth of order $p$ and $p-1$, respectively. For such $\phi$ we consider potentials $\Phi_1:U\rightarrow (-\infty,\infty]$ defined by

\begin{align*}
	\Phi_1(u)=\begin{cases}\int_0^1\phi(u(\xi))\,\mathrm{d}\xi,& u\in \tilde{L}^p(0,1), \\ \infty, & u\notin \tilde{L}^p(0,1). \end{cases}
\end{align*}

Before we investigate the Sobolev regularity of $\Phi_1$ we need some auxiliary results from \cite[Section 6]{Lunardi14}.Below, \Cref{lemma_phi_bound} includes a refinement of \cite[Proposition 6.5]{Lunardi14}, where only the case $\tilde{\alpha}=\frac{1}{2}$ was considered.
\begin{lem}\label{propo_bounded_Cahn_poly}
	For all $q\in [0,\infty)$  there is a constant $C_q\in (0,\infty)$ such that
	\[
		\int_U\int_0^1\abs{{P}^U_n(u)(\xi)}^q\,\mathrm{d}\xi\dm_1\leq C_q\left(\sum_{i=1}^n\frac{1}{(\pi k)^2}\right)^{\frac{q}{2}}
	\]
	and $\mu_1(\tilde{L}^q(0,1))=1$.
\end{lem}
\begin{lem}\label{lemma_phi_bound}
	$\Phi_1\in W^{1,2}_{Q^{\tilde{\alpha}}}(U;\mu_1)$ for all $\tilde{\alpha}\in (\frac{3}{8},\infty)$ and it satisfies Assumption \nameref{hypoapproximationPotential}. Moreover, if $p=1$, i.e.~if $\phi'$ is bounded we additionally have
\[
		\norm{Q^{\tilde{\alpha}}D\Phi_1}_{L^{\infty}(\mu_1)}<\infty.
\]
\end{lem}
\begin{proof}
	For $n\in\N$ and $u\in U$ we define
	\[
		(\Phi_1)_n(u)=\int_0^1\phi(P_n^U(u)(\xi))\,\mathrm{d}\xi.
	\]
	As in \cite[Proposition 6.5]{Lunardi14} one can show that $((\Phi_1)_n)_{n\in\N}$ is a sequence of continuously differentiable 				functions from $U$ to $\R$ converging to $\Phi_1$ in $L^q(U;\mu_1)$ for all $q\in [1,\infty)$. 		To conclude that $\Phi_1\in W^{1,2}_{Q^{\tilde{\alpha}}}(U;\mu_1)$ we use \cite[Lemma 5.4.4]{Bogachev}, i.e.~it is enough to show that $((\Phi_1)_n)_{n\in\N}		$ is bounded in $ W^{1,2}_{Q^{\tilde{\alpha}}}(U;\mu_1)$. Boundedness of $((\Phi_1)_n)_{n\in\N}$ in $L^q(U;\mu_1)$ for all $q\in [1,\infty)		$ follows by \Cref{propo_bounded_Cahn_poly} and the polynomial growth of $\phi$. As in \cite[Proposition 6.5]{Lunardi14} we obtain
	\[
		\partial_k(\Phi_1)_n(u)=\begin{cases}\int_0^1\phi'(P_n^U(u)(\xi))f_k(\xi)\,\mathrm{d}\xi& k\leq n \\ 0 & k>n \end{cases}.
	\]
	Hence for $k\leq n$
	\[
	\begin{aligned}
		\abs{\partial_k(\Phi_1)_n(u)}&=\abs{\int_0^1\phi'(P_n^U(u)(\xi))f_k(\xi)\,\mathrm{d}\xi}\leq A_1\sqrt{2}\pi k\int_0^1(1+\abs{P_n^U(u)(\xi)}^{p-1})\,\mathrm{d}\xi\\
		&\leq A_1\sqrt{2}\pi k\norm{1+\abs{P_n^U(u)}}^{p-1}_{L^{2(p-1)}(0,1)}.
		\end{aligned}
	\]
	This yields
	\[
	\begin{aligned}
		\norm{Q^{\tilde{\alpha}}D(\Phi_1)_n(u)}_U^2&=\sum_{k=1}^n\frac{1}{(\pi k)^{8\tilde{\alpha}}}\abs{\partial_k(\Phi_1)_n(u)}^2\\&\leq 2A^2_1\sum_{k=1}^n \frac{ (\pi k)^2}{(\pi 			k)^{8\tilde{\alpha}}}\norm{1+\abs{P_n^U(u)}}^{2(p-1)}_{L^{2(p-1)}(0,1)}\\
		&\leq 2A^2_1 \norm{1+\abs{P_n^U(u)}}^{2(p-1)}_{L^{2(p-1)}(0,1)} \sum_{k=1}^{\infty} \frac{1}{(\pi k)^{8\tilde{\alpha}-2}}.
	\end{aligned}
	\]
	Hence, by Proposition \ref{propo_bounded_Cahn_poly} we get boundedness of $(Q^{\tilde{\alpha}}D(\Phi_1)_n)_{n\in\N}$ in $L^2(U;\mu_1)$ 	as desired. A closer look in \cite[Lemma 5.4.4]{Bogachev}, i.e.~using the Banach-Saks property of the Hilbert space $W^{1,2}_{Q^{\tilde{\alpha}}}(U;\mu_1)$ shows that there is a subsequence $(n_k)_{k\in\N}$ such that the Cesaro mean $\psi_N\defeq\frac{1}{N}\sum_{k=1}^N {(\Phi_1)}_{n_k}$ converges to $\Phi_1$ in $W^{1,2}_{Q^{\tilde{\alpha}}}(U;\mu_1)$. W.l.o.g.~we can assume that $\lim_{N\rightarrow\infty}\psi_N(x)=\Phi_1(x)$ and $\lim_{N\rightarrow\infty}Q^{\tilde{\alpha}}D\psi_N(x)=Q^{\tilde{\alpha}}D\Phi_1(x)$ for $\mu_1$-a.e.~$x\in U$. Since  $((\Phi_1)_n)_{n\in\N}$ is a sequence of convex, continuously differentiable functions from $U$ to $\R$ which are bounded from below, the same is true for $(\psi_N)_{N\in\N}$. Now for a sequence $(t_M)_{M\in\N}\subseteq (0,\infty)$ converging to zero denote by $\psi_{(M,N)}$ the Moreau-Yosida Approximation of $\psi_N$, compare \cite[Section 2.3]{Lunardi14}, of 			order $t_M$. It is well known that $\psi_{(M,N)}$ is a convex function from $U$ to $\R$ and that for all  for all $u\in U$ and $N,M\in\N$
	\begin{enumerate}[(i)]
		\item  $-\infty<\inf_{u\in U}\Phi(u)\leq\inf_{u\in U}\psi_N(u)\leq\psi_{M,N}(u)\leq \psi_N(u)$ as well as 
		$\lim_{M\rightarrow \infty}\psi_{M,N}(u)=\psi_N(u)$
		\item $\psi_{M,N}$ is Fr\'{e}chet-differentiable with Lipschitz continuous gradient.
		\item  $\lim_{M\rightarrow \infty}\norm{(Q^{\tilde{\alpha}}D\psi_{M,N}-Q^{\tilde{\alpha}}D\psi_{N},f_i) _U}_{L^2(U;\mu_1)}=0$ for all 					$i\in\N$.
	\end{enumerate}
	Using the approximation properties of $(\psi_N)_{N\in\N}$, we obtain Assumption \nameref{hypoapproximationPotential}.\\
	To show the last statement assume $p=1$. This shows that $Q^{\tilde{\alpha}}D(\Phi_1)_{n_k}$ is bounded in $U$ independent of $k$. 			Hence the same holds true for $Q^{\tilde{\alpha}}D\psi_N$. Since  $Q^{\tilde{\alpha}}D\psi_N$ converges pointwisely to 						$Q^{\tilde{\alpha}}D(\Phi_1)$, $\mu_1$-almost everywhere, also the last statement is shown.
\end{proof}
\begin{remark}
	As in the previous example we are also able to consider potentials of the form
	\[
		\Phi\defeq \Phi_1+\Phi_2,
	\]
	where $\Phi_2:U\rightarrow \R$ is bounded and two times continuously Fr\'{e}chet-differentiable with bounded first and second order 			derivative.
\end{remark}

In order to satisfy Assumption \nameref{ass:pert-pot} we need to find $\alpha\in (0,\infty)$ such that $\alpha\alpha_1\in (\frac{3}{8},\infty)$ and
	\begin{equation}\label{cahn_ess}
		\sigma_2\leq -2{\alpha\alpha_1}+2\sigma_1.
	\end{equation}
In this case all assumptions for essential m-dissipativity are fulfilled and we obtain essential m-dissipativity in $L^2(W;\mu^{\Phi})$ of $(L_{\Phi},\fbs{B_W})$.
As in the stochastic reaction-diffusion setting there is an associated right process 
\[
	\mathbf{M}=(\Omega,\mathcal{F},(\mathcal{F}_t)_{t\geq 0}, (X_t,Y_t)_{t\geq 0},(P_w)_{w\in W})
\]with infinite life time, such that its transition semigroup and resolvent coincides in $L^2(W;\mu^{\Phi})$ with the ones generated by $(L_{\Phi},\fbs{B_W})$.
To show that the process has weak continuous paths we again invoke \Cref{rem:inf-dim-bounding-rho} and assume
\[
	-\frac{\alpha_1}{2}+\sigma_1+\frac{\alpha_2}{2}>\frac{1}{4}\quad\text{and}\quad -\frac{\alpha_2}{2}+\sigma_1+\frac{\alpha_1}{2}				>\frac{1}{4}
\] to verify Item (i) from \nameref{ass:inf-dim-diff-bound} as well as $\sigma_2\in (\frac{1}{4},\infty)$ and $\sigma_2-\alpha_2\in (\frac{1}{4},\infty)$ to verify Item (ii) and (iii) from \nameref{ass:inf-dim-diff-bound}, respectively.
Since we already assume that $\sigma_2-\alpha_2\in (\frac{1}{4},\infty)$ Assumption \nameref{ass:weak_sol} holds true and we get a stochastically and analytically weak solution in the sense of \Cref{thm:sto_ana_weak_sol}.
Concerning hypocoercivity, we note that Assumption \nameref{ass:qc-negative-type} is obviously valid. But to obtain the final hypocoercivity result we again distinguish to cases. First suppose $\sigma_2=\frac{3}{2}\alpha_2$.
In order to show that $\Phi$ fulfills the Assumptions \nameref{ass:ess-N_Phi} and \nameref{ass:L_Phi-convex-reg-est} it is required that $\alpha_1>\frac{3}{4}$ and $\phi$ is scaled such that $\norm{Q_1^{\frac{1}{2}}D\Phi}_{L^{\infty}(\mu_1)}<\frac{1}{2}$ or if $2\sigma_1-\alpha_2\geq 0$ that $\frac{2\sigma_1-\alpha_2}{2}\in (\frac{3}{8},\infty)$.
Item (i) and (ii) from assumption
\nameref{ass:inf-dim-matrix-bounded-ev} hold for 
\[
	C_1=0\quad\text{and}\quad C_{2}(v)=A 	(1+\norm{v}_U^{2}),\quad v\in U.
\]
For Item (iii) note that $\alpha_{n}^{22}(v)\leq 2\tilde{c}(2+\norm{v}_U)\lambda_{n}^{\frac{3\alpha_2}{2}-\frac{\alpha_2}{2}}$, for all $v\in U$, which describes an $\ell^2$-sequence, since
\[
	\frac{3\alpha_2}{2}-\frac{\alpha_2}{2}= {\alpha_2}> \frac{1}{4}.
\]
Moreover
\[
	\int_U  \|({\alpha_n^{22}}(v))_{n\in\N}\|^2_{\ell^2}\,\mu_2(\mathrm{d}v)\leq 4\tilde{c}^2\|(\lambda_{n}^{\alpha_2})_{n\in\N}\|			^2_{\ell^2}\int_U(2+\norm{v}_U)^2\,\mu_2(\mathrm{d}v)<\infty
\]
The second option is to verify Assumption \nameref{ass:inf-dim-matrix-bounded-ev_star} using \Cref{remark:alter_K_seven}. 
Indeed, Item (i) and (ii) from Assumption
\nameref{ass:inf-dim-matrix-bounded-ev_star} hold for 
\[
	-{\alpha_2}+\sigma_2-\sigma_1+\frac{\alpha_1}{2}>0\quad\text{with}\quad C_1=0\quad\text{and}\quad C_{2}(v)=A 	(1+\norm{v}_U^{2})			\quad v\in U.
\]
For Item (iii) note that $\alpha_{n}^{22}(v)\leq 2\tilde{c}(2+\norm{v}_U)\lambda_{n}^{\frac{\alpha_1}{2}-\sigma_1+\sigma_2}$, for all $v\in U$, which describes an $\ell^2$-sequence, since by the inequality above
\[
	\frac{\alpha_1}{2}-\sigma_1+\sigma_2> \alpha_{2}>\frac{1}{4}.
\]
Moreover, 
\[
	\int_U  \|({\alpha_n^{22}}(v))_{n\in\N}\|^2_{\ell^2}\,\mu_2(\mathrm{d}v)\leq 4\tilde{c}^2\|(\lambda_{n}^{\frac{\alpha_1}{2}-			\sigma_1+\sigma_2})_{n\in\N}\|^2_{\ell^2}\int_U(2+\norm{v}_U)^2\,\mu_2(\mathrm{d}v)<\infty.
\]
Since we want to allow unbounded $(C,D(C))$ in this situation, we assume $\Phi_2=0$, $\alpha_1\in (\frac{3}{4},\infty)$ and scale $\Phi_1$ such that $\norm{Q_1^{\frac{1}{2}}D\Phi}_{L^{\infty}(\mu_1)}<\frac{1}{2}$ to verify assumption \nameref{ass:ess-N_Phi} and \nameref{ass:L_Phi-convex-reg-est}.
In analogy to the previous section, we obtain \nameref{ass:inf-dim-v-poincare}
and \nameref{ass:inf-dim-u-poincare}, if
\[
	2\sigma_1-\alpha_2\leq \alpha_1.
\]

As already mentioned we can interpret $-B^2$ as a realization of the negative fourth order derivative with zero boundary conditions for the first and third order derivative, in the following denoted by $-\frac{\partial^4}{\partial\xi^4}$. In particular $Q=(\frac{\partial^4}{\partial\xi^4})^{-1}$. The arguments from \cite[Section 6]{Lunardi14} justify to interpret $D\Phi_1(u)$ as $\frac{\partial^2}{\partial\xi^2}\phi'(u)$, $u\in U$.

\begin{summary}\label{rem:table_rde_two}
For parameters $\alpha_1,\alpha_2\in (\frac{1}{4},\infty)$, $\sigma_1,\sigma_2,\alpha\in [0,\infty)$ assume that $\sigma_2\leq -2\alpha\alpha_1+2\sigma_1$ and $\alpha\alpha_1\in(\frac{3}{8},\infty)$. Let $\Phi$ be as above with $\norm{Q_1^{\alpha}D\Phi}_{L^{\infty}(\mu_1)}<\infty$. Since then Assumption \nameref{ass:inf-dim-elliptic}-\nameref{ass:pert-pot} hold, we know that the generator $L^{\Phi}$ of the degenerate stochastic Cahn-Hilliard type equation
\[
\begin{aligned}
	\mathrm{d}X_t&=(\frac{\partial^4}{\partial\xi^4})^{-\sigma_1+\alpha_2}Y_t\,\mathrm{d}t\\
	\mathrm{d}Y_t&=\sum_{i=1}^\infty \partial_iK_{22}(Y_t)f_i-K_{22}(Y_t)(\frac{\partial^4}{\partial\xi^4})^{\alpha_2}Y_t-						(\frac{\partial^4}{\partial\xi^4})^{-\sigma_1+\alpha_1} X_t\\
	&\qquad-(\frac{\partial^4}{\partial\xi^4})^{-\sigma_1}\frac{\partial^2}{\partial\xi^2}\phi'(X_t)-(\frac{\partial^4}{\partial\xi^4})^{-			\sigma_1}D\Phi_2(X_t)\mathrm{d}t
	+\sqrt{2K_{22}(Y_t)}\,\mathrm{d}B_t
\end{aligned}
\]
is essentially m-dissipative.
The existence of a corresponding right process 
\[
	\mathbf{M}=(\Omega,\mathcal{F},(\mathcal{F}_t)_{t\geq 0}, (X_t,Y_t)_{t\geq 0},(P_w)_{w\in W})
\]
with infinite life time $P_{\mu^{\Phi}}$ a.s., solving the martingale problem $(L^{\Phi},D(L^{\Phi}))$ under $P_{\mu^{\Phi}}$ follows. As in \Cref{rem:table_rde_one}, $K_{22}:U\to \lop{U}$ is not finitely based and not necessarily mapping into the space of trace class operators.
We summarize the conditions on the parameters and the potential such that $\mathbf{M}$ provides a stochastically and analytically weak solution with weak continuous paths to the equation above and that the semigroup generated by $(L^{\Phi},D(L^{\Phi}))$ is hypocoercive.

\smallskip

\begin{center}
\begin{tabular}{ | m{5,4em} | m{4,4cm}| m{5,25cm} | } 
\hline
 & $\mathbf{M}$ stochastically and analytically weak solution with weak continuous paths  & $\sccs$ is hypocoercive \\
\hline
\nameref{ass:inf-dim-diff-bound} & $\sigma_2-\alpha_2$, $\pm\nicefrac{\alpha_1}{2}+\sigma_1\mp\nicefrac{\alpha_2}{2} >\nicefrac{1}{4}$&  \\
\hline 
\nameref{ass:weak_sol} & follows by \nameref{ass:inf-dim-diff-bound} &\\
\hline
\nameref{ass:qc-negative-type} &  by construction &\\
\hline
\nameref{ass:ess-N_Phi}, \nameref{ass:L_Phi-convex-reg-est}, \nameref{ass:inf-dim-matrix-bounded-ev} &  & $\sigma_2=\nicefrac{3}{2}\alpha_2$ and either ($\alpha_1>\nicefrac{3}{4}$ and $\norm{Q_1^{\frac{1}{2}} D\Phi}_{L^{\infty}(\mu_1)}<\nicefrac{1}{2}$) or ($2\sigma_1\geq \alpha_2$ and $\nicefrac{(2\sigma_1-\alpha_2)}{2}>\nicefrac{3}{8}$)\\
\hline
\nameref{ass:ess-N_Phi}, \nameref{ass:L_Phi-convex-reg-est}, \nameref{ass:inf-dim-matrix-bounded-ev_star}&  & $-{\alpha_2}+\sigma_2-\sigma_1+\nicefrac{\alpha_1}{2}>0$, $\alpha_1 >\nicefrac{3}{4}$, $\Phi_2=0$ and $\norm{Q_1^{\frac{1}{2}}D\Phi}_{L^{\infty}(\mu_1)}<\nicefrac{1}{2}$\\
\hline
\nameref{ass:inf-dim-v-poincare}, \nameref{ass:inf-dim-u-poincare}&  & $2\sigma_1-\alpha_2\leq \alpha_1$\\
\hline
\end{tabular}
\end{center}

\smallskip

As in the previous example, we can combine the results described above to verify that $\textbf{M}$ is $L^2$-exponentially ergodic. Choosing $\sigma_1=\alpha_2$ and the corresponding parameters accordingly, we are able rewrite the system of coupled equations above, as an infinite-dimensional second order in time stochastic differential equation.
\end{summary}

\backmatter

\bmhead{Acknowledgments}

The second author gratefully acknowledges financial support in the form of a
fellowship from the ``Studienstiftung des deutschen Volkes''.

\section*{Declarations}

%Some journals require declarations to be submitted in a standardised format. Please check the Instructions for Authors of the journal to which you are submitting to see if you need to complete this section. If yes, your manuscript must contain the following sections under the heading `Declarations':
\section{Funding}
Not applicable
\section{Conflict of interest/Competing interests}
Not applicable
\section{Ethics approval}
Not applicable

\section{Consent to participate}
Not applicable

\section{Consent for publication}
Not applicable

\section{Availability of data and materials}
Not applicable
\section{Code availability}
Not applicable
\section{Authors' contributions}
Not applicable

\bibliography{sources}
%


\end{document}